\documentclass[12pt,a4paper]{article}
\usepackage[utf8]{inputenc}
\pdfoutput=1

\newcommand\Tstrut{\rule{0pt}{2.9ex}}       
\newcommand\Bstrut{\rule[-1.3ex]{0pt}{0pt}} 
\newcommand\TBstrut{\Tstrut\Bstrut}         

\usepackage{authblk} 
\usepackage{amssymb,amsthm}
\usepackage{mathtools}
\usepackage{multirow}
\usepackage{xspace}
\usepackage{amsmath,bm} 
\usepackage{algorithm}
\usepackage{algpseudocode}
\usepackage{color}
\usepackage{graphicx}
\usepackage{caption}

\usepackage{hyperref}
\usepackage[misc]{ifsym}

\providecommand{\url}[1]{{#1}}

\expandafter\ifx\csname urlstyle\endcsname\relax
  \providecommand{\doi}[1]{DOI~\discretionary{}{}{}#1}\else
  \providecommand{\doi}{DOI~\discretionary{}{}{}\begingroup
  \urlstyle{rm}\Url}\fi

\newcommand{\specfem}{{\tt SPECFEM2D}\xspace}
\newcommand{\QoI}{\ensuremath{\mathcal{Q}}}
\newcommand{\QoItext}{QoI}
\newcommand{\QoItexts}{QoIs}
\newcommand{\norm}[1]{\left\lVert#1\right\rVert}
\providecommand{\displ}{\ensuremath{\mathbf{s}}}     
\providecommand{\strain}{\ensuremath{\bm{\epsilon}}} 
\providecommand{\ParEar}{\ensuremath{\bm{\theta}}}   
\providecommand{\source}[1][~]{\ensuremath{\mathbf{x_{s}^{#1}}}}     
\providecommand{\sourcex}{\ensuremath{x_s}}          
\providecommand{\rec}[1][\mathbf{r}]{\ensuremath{\mathbf{x}_{#1}}}        
\providecommand{\Nrec}{\ensuremath{N_{rec}}}           
\providecommand{\data}{\ensuremath{\mathbf{d}}}      
\providecommand{\WhiteNoise}[1][{k,n}]{\ensuremath{\bm{\varepsilon}_{#1}}} 
\providecommand{\elasticmod}{\ensuremath{\mathbf{c}}}   

\providecommand{\stress}{\ensuremath{\mathbf{T}}}    
\providecommand{\density}{\ensuremath{\rho}}         
\providecommand{\shearmod}{\ensuremath{\mu}}         
\providecommand{\Lamefp}{\ensuremath{\lambda}}       
\providecommand{\Lamesp}{\shearmod}                  

\providecommand{\comprv}{\ensuremath{\alpha}}        
\providecommand{\shearv}{\ensuremath{\beta}}         

\providecommand{\force}{\ensuremath{\mathbf{f}}}     

\providecommand{\xv}{\ensuremath{\mathbf{x}}}        
\providecommand{\domain}{\ensuremath{D}}             
\providecommand{\boundary}{\ensuremath{\partial D}}  
\providecommand{\surface}{\ensuremath{\partial D_S}} 
\providecommand{\artbound}{\ensuremath{\partial D_A}}
\providecommand{\normal}{\ensuremath{\mathbf{\hat{n}}}} 
\providecommand{\testfun}{\ensuremath{\mathbf{w}}}   
\providecommand{\solspace}{\ensuremath{\mathbf{V_\displ}}}
\providecommand{\Imatrix}{\ensuremath{\mathbf{I}}}   

\providecommand{\samplespace}{\ensuremath{\Omega}}
\providecommand{\sample}{\ensuremath{\omega}}

\providecommand{\unifrnd}{\ensuremath{\mathcal{U}}}

\providecommand{\doubledot}{\mathbin{\vcentcolon}} 
\providecommand{\rset}{\mathbb{R}} 
\providecommand{\grad}{\boldsymbol{\nabla}} 
\providecommand{\divergence}{\boldsymbol{\nabla\cdot}} 

\DeclareMathOperator{\trace}{Tr}

\providecommand{\e}[1]{\ensuremath{\cdot 10^{#1}}}
\providecommand{\E}[1]{{\ensuremath{\mathrm{E}}\mspace{-2mu}\left[#1\right]}}
\providecommand{\prob}[1]{{\ensuremath{\mathrm{P}}\mspace{-2mu}\left[#1\right]}}
\providecommand{\var}[1]{{\ensuremath{\mathrm{Var}}\mspace{-2mu}\left[#1\right]}}
\providecommand{\tol}{\ensuremath{\mathrm{TOL}}} 
\providecommand{\work}{\ensuremath{W}}           
\providecommand{\nrs}{N}                         
\providecommand{\Est}{\mathcal{A}}               
\providecommand{\EstMLMC}{\mathcal{A}_{MLMC}}    
\providecommand{\EstMC}{\mathcal{A}_{MC}}        
\providecommand{\varMC}{\mathcal{V}}             
\providecommand{\vg}{{\ensuremath{V_0}}}         
\providecommand{\Cw}{{K_{w}}}                    
\providecommand{\Cs}{{K_{s}}}                    
\providecommand{\Ow}{{q_{w}}}                    
\providecommand{\Os}{{q_{s}}}                    
\providecommand{\splitting}{\varphi}             
\providecommand{\confidence}{\xi}                
\providecommand{\confpar}{C_\xi}                  
\providecommand{\meshparam}[1][\ell]{h_{#1}}     
\providecommand{\dx}[1][\ell]{\Delta x_{#1}}   
\providecommand{\dt}[1][\ell]{\Delta t_{#1}}     

\newtheorem{problem}{Problem}
\theoremstyle{remark}
\newtheorem{remark}{Remark}

\title{Multilevel Monte Carlo Acceleration of Seismic~Wave Propagation
  under Uncertainty} 

\author[1,2]{Marco~Ballesio~\thanks{Email: marco.ballesio@kaust.edu.sa}}
\author[1,2]{Joakim~Beck}
\author[1,2]{Anamika~Pandey~\thanks{Email: anamika.pandey@kaust.edu.sa}}
\author[3]{Laura~Parisi}
\author[1,2]{Erik~von~Schwerin}
\author[1,2,4]{Ra\'{u}l~Tempone}
\affil[1]{\small Computer, Electrical and Mathematical Sciences and Engineering (CEMSE), King Abdullah University of Science and Technology (KAUST), Thuwal 23955-6900, Kingdom of Saudi Arabia.}
\affil[2]{\small KAUST SRI Center for Uncertainty Quantification in Computational
  Science and Engineering}
\affil[3]{\small Physical Science and Engineering (PSE), King Abdullah University of Science and Technology (KAUST), Thuwal 23955-6900, Kingdom of Saudi Arabia.}
\affil[4]{\small Alexander von Humboldt Professor in Mathematics for Uncertainty
  Quantification, RWTH Aachen University, Germany}

\date{}

\begin{document}

\maketitle
\begin{abstract}  We interpret uncertainty in a model for seismic wave
  propagation by treating the model parameters as random variables,  
  and apply the Multilevel Monte Carlo (MLMC) method to reduce the
  cost of approximating expected values of selected, physically
  relevant, quantities of interest (\QoItext) with respect to the
  random variables. 

  Targeting source inversion problems, where the source of an
  earthquake is inferred from ground motion recordings on the Earth's
  surface, we consider two \QoItexts~that measure the discrepancies
  between computed seismic signals and given reference signals: 
  one \QoItext, $\QoI_E$, is defined in terms of the $L^2$-misfit,
  which is directly related to maximum likelihood estimates of the
  source parameters; the other, $\QoI_W$, is based on the quadratic
  Wasserstein distance between probability distributions, and
  represents one possible choice in a class of such misfit functions
  that have become increasingly popular to solve seismic inversion in
  recent years. 

  We simulate seismic wave propagation, including seismic attenuation,
  using a publicly available code in widespread use, based on the
  spectral element method.
  Using random coefficients and deterministic initial and boundary
  data, we present benchmark numerical experiments with synthetic data
  in a two-dimensional physical domain and a one-dimensional velocity
  model where the assumed parameter uncertainty is motivated by
  realistic Earth models. 
  Here, the computational cost of the standard Monte Carlo method
  was reduced by up to $97\%$ for $\QoI_E$, and up to $78\%$ for
  $\QoI_W$, using a relevant range of tolerances. 
  Shifting to three-dimensional domains is straight-forward and will
  further increase the relative computational work reduction.

\end{abstract}


\pagestyle{myheadings}
\thispagestyle{plain}
\markboth{MLMC Acceleration of Seismic Wave Propagation under Uncertainty}
  {MLMC Acceleration of Seismic Wave Propagation under Uncertainty}

\section{Introduction}
\label{sec:intro}

Recent large earthquakes and their devastating effects on society and
infrastructure (e.g., New Zealand, 2011; Japan, 2011; Nepal, 2015)
emphasize the urgent need for reliable and robust earthquake-parameter
estimations for subsequent risk assessment and mitigation. Seismic
source inversion is a key component of seismic hazard assessments
where the probabilities of future earthquake events in the region are
of interest.  

From ground motion recordings at the surface of the Earth, i.e.
seismograms, we are interested in efficiently computing the likelihood 
of postulated parameters describing the unknown source of an earthquake. 
A sub-problem of the seismic source inversion is to infer the location
of the source (hypocenter) and the origin time. 
We take the expected value of the quantity of interest, which, in our
case, is the misfit between the observed and predicted ground
displacements for a given seismic source location and origin time, to
find the location and time of highest likelihood from observed seismogram data.
Several mathematical and computational approaches can be used to
calculate the predicted ground motions in the source inversion problem. 
These techniques span from approximately calculating only some 
of the waveform attributes, (e.g. peak ground acceleration or 
seismic phase arrival-times), often by using simple one-dimensional
velocity models, to the simulation of the full wave propagation in a
three-dimensionally varying structure.

In this work, the mathematical model and its output are random to
account for the lack of precise knowledge about some of its
parameters. In particular, to account for uncertainties in the
material properties of the Earth, these are modeled by random variables. 
The most common approach used to compute expected values of random variables 
is to employ Monte Carlo (MC) sampling. 
MC is non-intrusive, in the sense that it doesn't
require underlying deterministic computational codes to be modified,
but only called with randomly sampled parameters. Another striking
advantage of MC is that no regularity assumptions are needed on the
quantity of interest with respect to the uncertain parameters, other
than that the variance has to be bounded to use the Central Limit
Theorem to predict the convergence rate. 
However, in situations when
generating individual samples from the computational model is highly
expensive, often due to the need for a fine, high resolution,
discretization of the physical system, MC can be too costly to
use. To avoid a large number of evaluations of the computer model with
high resolution, but still preserve the advantages of using MC, we
apply Multilevel Monte Carlo (MLMC)
sampling~\cite{giles08,heinrich01,giles_AcNum} to substantially reduce
the computational cost by distributing the sampling over computations
with different discretization sizes. 
Polynomial chaos surrogate models
have been used to exploit the regularity of the waveform
solution~\cite{gmd-2018-4}; in cases when the waveform is analytic
with respect to the random parameters, the asymptotic convergence can
be super-algebraic. 
However, if the waveform is not analytic, one typically only achieves
algebraic convergence in $L^2$ asymptotically, 
as was shown in~\cite{Motamed2013} to be the case with stochastic
collocation for the second order wave equation with discontinuous
random wave speed. 
The setup in~\cite{Motamed2013} with a stratified medium is analogous
to the situation we will treat in the numerical experiments for a
viscoelastic seismic wave propagation problem in the present paper. 
The efficiency of polynomial chaos, or stochastic collocation, methods
deteriorates as 
the number of ``effective'' random variables increases, whereas the MC
only depends weakly on this number. Here MC
type methods have an advantage for the problems we are interested in.

We are motivated by eventually solving the inverse problem; therefore 
the Quantities of Interest (\QoItexts) considered in this paper are
misfit functions quantifying the distance between observed and
predicted displacement time series at a given set of seismographs on
Earth's surface. 
One \QoItext~is based on the $L^2$-norm of the distances; another makes use of
the Wasserstein distance between probability densities, which requires
transformation of the displacement time series to be applicable.
The advantages of the Wasserstein distance over the $L^2$
difference for full-waveform inversion, for instance, that the former
circumvents the issue of cycle skipping, have been shown
in~\cite{Engquist2014application} and further studied 
in~\cite{engquist2016optimal,Engquist2018Analysis,Engquist2018Application}. 
A recent preprint~\cite{Motamed_Appelo} combines this type of
\QoItexts~with a Bayesian framework for inverse problems.

In our demonstration of MLMC for the two \QoItexts, we consider a
seismic wave propagation in a 
semi-infinite two-dimensional domain with free surface boundary
conditions on the Earth's surface, heterogeneous
viscoelastic media, and a point-body time-varying
force. 
We use \specfem~\cite{SPECFEM2D,tromp2008spectral} for the numerical 
computation where we consider an isotropic viscoelastic Earth model.
The heterogeneous media is divided into homogeneous horizontal
layers with uncertain densities and velocities. The densities and
shear wave velocities are treated as random, independent between the
subdomains, and they are uniformly distributed over the assumed
intervals of uncertainty in the respective layers. 
The compressional wave
velocities of the subdomains follow a multivariate uniform
distribution conditional on the shear wave speeds. Our choice of
probability distributions to describe the uncertainties are motivated by
results given in~\cite{Albaric2010,Roecker2017}. 

The paper is outlined as follows:
In Section~\ref{sec:prob}, the seismic wave propagation problem is
described for a viscoelastic medium with random Earth material
properties. The two \QoItexts~are described in Section~\ref{sec:QoI}.
The computational techniques, including (i) numerical approximation of the 
viscoelastic Earth material model and of the resulting initial
boundary value problem, the combination of which is taken as
``black-box'' solver by the widely used seismological software package
{\tt SPECFEM2D}~\cite{SPECFEM2D}, and (ii) MLMC approximation of
\QoItexts~depending on the random Earth material properties, are
described in Section~\ref{sec:comp}.
The configuration of the numerical tests is described in
Section~\ref{sec:num}, together with the results, showing  a considerable
decrease in computational cost 
compared to the standard Monte Carlo approximation for the same accuracy.

\section{Seismic Wave Propagation Model with \\ Random
  Parameters}
\label{sec:prob}

Here we describe the model we use for seismic wave propagation in a
heterogeneous Earth medium, given by an initial boundary value problem
(IBVP). 
We interpret the inherent uncertainty in the Earth material properties
through random parameters which define the compressional and shear
wave speed fields and the mass density.

First, we state the strong form of the IBVP in the case of a
deterministic elastic Earth model, later to be extended to a
particular anelastic model in the context of a weak form of the IBVP,
suitable for the numerical approximation methods used in
Section~\ref{sec:comp} and~\ref{sec:num}. 
Finally, we state assumptions on the random material parameter fields.

\subsection{Strong Form of Initial Boundary Value Problem}
\label{sec:strong}

We consider a heterogeneous medium occupying a
domain $\domain \subset \rset^3$ modeling the Earth.
We denote by 
$\displ\doubledot\domain \times (0,\mathcal{T}]\to\rset^3$ 
the space-time displacement field induced by a seismic event in $\domain$. 
In the deterministic setting $\displ$ is assumed to satisfy 
\begin{subequations}
 \label{eq:strong_eq_motionMain}
  \begin{align}
    \label{eq:strong_eq_motion}
    \density(\xv)\partial_t^2\displ(\xv, t) - \divergence\stress(\grad
    \displ(\xv, t)) & = \force(\xv,t), && 
    \text{~$\forall(\xv, t)\in\domain \times(0,\mathcal{T}]$}, 
      \intertext{for some finite time horizon, given by $\mathcal{T}>0$, 
        with the initial conditions}
    \label{eq:InitialCondition}
    \begin{array}{r}
      \displ(\xv,0)\!\!\!\\
      \partial_{t}\displ(\xv,0)\!\!\!
    \end{array}
    & 
    \begin{array}{l}
      = \mathbf{g}_{1}(\xv),\\
      = \mathbf{g}_{2}(\xv),
    \end{array}
    && \text{~$\forall\xv\in \domain$,}
    \intertext{and the free surface boundary condition on Earth's
       surface $\surface$}
    \label{eq:free_surface}
    \normal\cdot\stress & = 0, && \text{on \surface,}
  \end{align}
\end{subequations}
where $\stress$ denotes the stress tensor, and $\normal$ denotes the
unit outward normal to $\surface$. 
Together with a constitutive relation between stress and strain,
\eqref{eq:strong_eq_motion}--\eqref{eq:free_surface} form 
an IBVP for seismic wave propagation; two different
consitutive relations will be considered below.
In this paper $\partial_{t}$ denotes time derivative and 
$\grad$ and $\divergence$ denote spatial gradient and divergence
operators, respectively. 

With $\density$ denoting the density, $\force$ becomes a body force
that includes the force causing the seismic event. In this study, we
consider a simple point body force acting with time-varying magnitude
at a fixed point, as described in~\cite{aki_richards2002,Dahlen_Tromp1998}.

For an isotropic elastic Earth medium undergoing infinitesimal
deformations, the constitutive stress-strain relation can be described
by 
\begin{equation}
  \label{eq:Actual_HookesLaw}
  \stress(\grad \displ) = 
  \Lamefp \trace(\strain(\grad\displ))\Imatrix + 2\Lamesp\,\strain(\grad\displ),
\end{equation}
with $\trace(\strain)$ the trace of the symmetric infinitesimal strain tensor, 
$\strain(\grad \displ) = \frac{1}{2} \left[\grad \displ + \left(\grad
    \displ\right)^\intercal\right]$, and $\Imatrix$ the 
identity tensor; see~\cite{aki_richards2002,Dahlen_Tromp1998,Carcione2014}. 
In the case of isotropic heterogeneous elastic media undergoing
infinitesimal deformations, the first and second Lam\'{e} parameter,
denoted $\Lamefp$ and 
$\Lamesp$ respectively, are functions of the spatial position, but
for notational simplicity, we often omit the dependencies on $\xv$. 
These parameters, together, constitute a parametrization of the
elastic moduli for homogeneous isotropic media and together with
$\density$ also determine the compressional wave speed, $\comprv$, and shear
wave speed, $\shearv$, by  
\begin{equation}
 \comprv  = \sqrt{\dfrac{\Lamefp+2 \Lamesp}{\density}}, 
~~\quad ~~\quad  \shearv = \sqrt{\frac{\shearmod}{\density}}.
\end{equation}

Either one of the triplets
$(\density,\Lamefp,\Lamesp)$ and $(\density, \comprv, \shearv)$ defines
the Earth's material properties with varying spatial position for a general
velocity model~\cite{Virieux1986}.

\paragraph{Simplification of the full Earth model}

For the purpose of the numerical computations and 
the well-posedness of the underlying wave propagation models, 
we will later replace the whole Earth domain by a semi-infinite domain, 
which we will truncate with absorbing boundary conditions at the artificial
boundaries introduced by the truncation. 
The domain boundary is $\boundary=\surface\cup\artbound$
with $\artbound$ denoting the artificial boundary.
From now on, we will consider $\domain$ to be an open 
bounded subset of $\rset^d$, where $d=2~\mathrm{or}~3$ denotes the dimension
of the physical domain. 
We will consider numerical examples with $\domain \subset
\rset^2$. 

\subsection{Weak Form of Initial Boundary Value Problem} 
\label{sec:weak}

The numerical methods for simulating the seismic wave propagation used
in this paper are based on an alternative formulation of
IBVP~\eqref{eq:strong_eq_motionMain} that uses the weak form in
space. To obtain such form, one multiplies~\eqref{eq:strong_eq_motion}
at time $t$ by a sufficiently regular test function $\testfun$, and
integrates over the physical domain $\domain$.
Using integration by parts and imposing the traction-free boundary
condition~\eqref{eq:free_surface}, one derivative is shifted
from the unknown displacement, $\displ$, to the test function
$\testfun$, i.e.,
\begin{align*}
  \int_\domain\testfun\cdot\divergence\stress \,d\xv & = 
  - \int_\domain\grad\testfun\doubledot\stress \,d\xv,
\end{align*}
where $\doubledot$ denotes the double contraction.
In this context, ``sufficiently regular test function'', means that
$\testfun \in \mathbf{H}^1(\domain)$, where $\mathbf{H}^1(\domain)$
denotes the Sobolev space $\mathbf{W}^{1,2}(\domain)$, 
\begin{align*}
  \mathbf{H}^1(\domain) & = 
  \left\lbrace
    \mathbf{w}: \domain \rightarrow \rset^{d}
    \,\, \text{s.t.} \,
    \norm{\mathbf{w}}_{\mathbf{H}^1(\domain)} < \infty
  \right\rbrace,
\end{align*}
equipped with the usual inner product
\begin{align*}
  \left\langle \mathbf{u},\mathbf{v} \right\rangle_{\mathbf{H}^1(\domain)}
  & = 
  \int_\domain \left( 
    \mathbf{u}\cdot\mathbf{v} + 
    l^2 \grad \mathbf{u}\doubledot \grad \mathbf{v} \right),
\end{align*}
where $l$ is a characteristic length scale, 
and the corresponding induced norm
\begin{align*}
  \norm{\mathbf{u}}_{\mathbf{H}^1(\domain)}
  & = 
    \sqrt{
    \left\langle \mathbf{u},\mathbf{u}
    \right\rangle_{\mathbf{H}^1(\domain)}
    }.
\end{align*}
The weak form of the IBVP then becomes: 
\begin{problem}[Weak form of isotropic elastic IBVP]
\label{prob:weak1}
Find $\displ\in\solspace$, which both satisfies the initial
conditions~\eqref{eq:InitialCondition}, and for $\stress$
in~\eqref{eq:Actual_HookesLaw} satisfies 
\begin{align}
  \label{eq:weak_eq_motion}
  \int_\domain\density\testfun\cdot\partial_t^2\displ \,d\xv & = 
  - \int_\domain\grad\testfun\doubledot\stress\left(\grad\displ\right) \,d\xv
  + \int_\domain \force\cdot\testfun \,d\xv,  && \forall~ \testfun \in \mathbf{H}^1(\domain),
\end{align}
almost everywhere in the time interval $\left[0,\mathcal{T}\right]$,
where the trial space
\begin{align}
  \solspace & =
  \left\lbrace
    \displ: [0,\mathcal{T}] \rightarrow \mathbf{H}^1(D)
    \, \left\vert \, 
    \begin{array}{l}
      \displ \in \mathbf{L}^{2}(0, \mathcal{T}; \mathbf{H}^{1}(\domain)),
      \\
      \partial_t\displ \in 
        \mathbf{L}^{2}(0, \mathcal{T}; \mathbf{L}^{2}(\domain)), 
      ~\text{and}~ 
      \\
      \partial_{tt}\displ \in \mathbf{L}^{2}(0, \mathcal{T}; \mathbf{H}^{-1}(\domain))
    \end{array}
  \right\rbrace
  \right.,
  \label{eq:solspace_P1}
\end{align}
using the spaces defined in~\eqref{eq:L2_t_X}.
\end{problem}

Above, the time dependent functions, 
$\displ(\xv,t)$, $\partial_t\displ(\xv,t)$,
$\partial_{tt}\displ(\xv,t)$, belong to Bochner spaces
\begin{align}
  \mathbf{L}^{2}(0, \mathcal{T}; \mathcal{X}) 
  & = 
  \left\lbrace \mathbf{u}:
    [0,\mathcal{T}] \rightarrow \mathcal{X},\text{~strongly measurable~}
    \left\vert
      \int_{[0,\mathcal{T}]} \norm{\mathbf{u}}_{\mathcal{X}}^{2} \,dt < \infty
    \right\rbrace
  \right. ,
  \label{eq:L2_t_X}
\end{align}
where the appropriate choice of $\mathcal{X}$ depends on the
number of spatial derivatives needed: $\mathbf{H}^1(\domain)$,
$\mathbf{L}^2(\domain)$, and $\mathbf{H}^{-1}(\domain)$,
respectively, with the latter space being the dual space of
$\mathbf{H}^1(\domain)$. 

According to~\cite{LionsI1972,LionsII1972}, Problem~\ref{prob:weak1}
is well-posed under certain regularity assumptions and appropriate
boundary condition on $\artbound$, see page 32, equation (5.29) in~\cite{LionsII1972}. 
More precisely, assuming that the
density~$\density$ is bounded away from zero, $\density\geq\density_{min}>0$, 
the problem fits into the setting of Section~1, Chapter~5,
of~\cite{LionsII1972}, after dividing through by $\density$.
Assuming that $\Lamefp,\Lamesp$ in~\eqref{eq:Actual_HookesLaw} are also
sufficiently regular,  
according to Theorem~2.1 in the chapter it holds that, if the force 
$\force/\density \in \mathbf{L}^{2}(0, \mathcal{T}; \mathbf{L}^{2}(\domain))$, 
the initial data
$\mathbf{g}_{1} \in \mathbf{H}^{1}(\domain)$
and 
$\mathbf{g}_{2} \in \mathbf{L}^{2}(\domain)$,
and the boundary, $\surface$, is infinitely differentiable, 
there exists a unique solution to Problem~\ref{prob:weak1},
\begin{align*}
  \displ & \in \mathbf{C}^{0}([0, \mathcal{T}]; \mathbf{H}^{1}(D)) 
           \cap \mathbf{C}^{1}([0, \mathcal{T}]; \mathbf{L}^{2}(D)) 
           \cap \mathbf{H}^{2}((0, \mathcal{T}); \mathbf{H}^{-1}(D)), 
\end{align*}
which depends continuously on the initial data. 

In our numerical experiments, however, we instead use a problem with
piecewise constant material parameters, violating the smoothness
assumption, but only at material interfaces in the interior of the
domain. Furthermore, a singular source term, $\force$, common in
seismic modeling, will be used.

\subsection{Weak Form Including Seismic Attenuation}
\label{sec:attenuation}

For a more
realistic Earth model, we include seismic attenuation in the wave
propagation. 
The main cause of seismic attenuation is the relatively small but not
negligible anelasticity of the Earth. 
In the literature, e.g., Chapter~6 of~\cite{Dahlen_Tromp1998} or
Chapter~1 of~\cite{Carcione2014}, anelasticity of the Earth is
modeled by combining the mechanical properties of elastic solids and
viscous fluids. In a heterogeneous linear isotropic viscoelastic medium, the
displacement field $\displ(\xv,t)$ 
follows
the IBVP~\eqref{eq:strong_eq_motionMain}, but the stress tensor
$\stress$ depends linearly upon the entire history of the
infinitesimal strain, and the constitutive
relation~\eqref{eq:Actual_HookesLaw} will be replaced by
\begin{align}
  \label{eq:anisoHookeLaw}
  \stress\left(\xv,t;\left\{\grad\displ\right\}_0^t\right)
  & = \int_{-\infty}^{t} \elasticmod(\xv, t - t^{'}) \doubledot 
  \partial_{t}\strain(\grad\displ(\xv, t^{'})) \,dt^{'}, 
\end{align}
where $\elasticmod$ represents the anelastic fourth
order tensor which accounts for the Earth's material properties, which will 
be further discussed in Section~\ref{sec:comp_attenuation}.  

With the constitutive relation~\eqref{eq:anisoHookeLaw}
replacing~\eqref{eq:Actual_HookesLaw},  
the IBVP becomes:
\begin{problem}[Weak form of isotropic viscoelastic IBVP]
\label{prob:weak2}
Find $\displ\in\solspace$, defined in~\eqref{eq:solspace_P1}, which both
satisfies the initial conditions~\eqref{eq:InitialCondition}, and for 
$\stress$ in~\eqref{eq:anisoHookeLaw} satisfies 
\begin{multline}
  \int_\domain\density\testfun(\xv)\cdot\partial_t^2\displ(\xv,t)\,d\xv = 
  - \int_\domain\grad\testfun\doubledot
    \stress\left(\xv,t;\left\{\grad\displ\right\}_0^t\right)\,d\xv
  + \int_\domain \force(\xv,t)\cdot\testfun(\xv)\,d\xv,  
  \\ \forall~ \testfun \in \mathbf{H}^1(\domain),
  \label{eq:weak_eq_motion_attenuation}
\end{multline}
almost everywhere in the time interval $\left[0,\mathcal{T}\right]$.
\end{problem}

In the case of a bounded domain, $\domain$, with homogeneous initial conditions, theoretical well-posedness results for the viscoelastic model considered in Problem~\ref{prob:weak2}, as well as a wide range of other viscoelastic models, are given in~\cite{Sayas2018}. More precisely, the viscoelastic material tensor $\elasticmod(\xv, t - t^{'})$ should follow the hypothesis of symmetricity, positivity, and boundedness, as given on page 60 in~\cite{Sayas2018}; the boundary, $\boundary$, can be a combination of nonoverlapping Dirichlet and Neumann parts, and we can define a bounded and surjective trace operator from $\mathbf{H}^1(\domain)$ to $\mathbf{H}^{1/2}(\boundary)$. For a detailed description of $H^{1/2}(\boundary)$ in the case of mixed type of boundaries, see page 58 in~\cite{Sayas2018}. Note that a $d-1$ times differentiable $\boundary$ with locally Lipschtiz $\domain$ will always have a trace operator.

\subsection{Statement in Stochastic setting} 
\label{sec:stochastic}

Here we model the uncertain Earth material properties,
$(\density,\Lamefp,\Lamesp)$, as time-independent random
fields
$(\density,\Lamefp,\Lamesp):\domain\times\samplespace\to\rset^3$,
where $\samplespace$ is the sample space of a complete probability
space.
We assume that the random fields are bounded from above
and below, uniformly both in physical space and in sample space, and
with the lower bounds strictly positive,
\begin{subequations}
  \label{eq:RandomParCond1}
  \begin{alignat}{6}
    0 & < & \density_{min}  & \leq & \density(\xv, \sample) & 
     \leq & \density_{max} & < \infty, 
    &\quad\quad& \forall \xv \in \domain, \forall \sample \in \samplespace, \\
    0 & < & \Lamefp_{min}  & \leq & \Lamefp(\xv, \sample) & 
    \leq & \Lamefp_{max} & < \infty, 
    &\quad\quad& \forall \xv \in \domain, \forall \sample \in \samplespace, \\
    0 & < & \Lamesp_{min}  & \leq & \Lamesp(\xv, \sample)  & 
    \leq & \Lamesp_{max} & < \infty, 
    &\quad\quad& \forall \xv \in \domain, \forall \sample \in \samplespace.
  \end{alignat}
\end{subequations}
For any given sample $\omega$ the displacement field
$\displ(\cdot,\cdot,\sample)$ solves
Problem~\ref{prob:weak1} or Problem~\ref{prob:weak2}, for the
respective case.

Any known well-posedness properties of the deterministic
Problem~\ref{prob:weak1} and Problem~\ref{prob:weak2} are directly  
inherited in their stochastic form, assuming the same regularity of
realizations of the random fields as of their deterministic counterparts.

\section{Quantities of Interest}
\label{sec:QoI}

Two \QoItext~suitable for different approaches to seismic inversion
will be described. 
The common feature is that they quantify the misfit between data,
consisting of ground motion measured at the Earth's surface, at fixed equidistant 
observation times $\{t_k^d\}_{k=0}^K$, $t_k=k\dt[]^d$, and model
predictions, consisting of the corresponding model predicted ground motion. 
The displacement data, $\data$, and the model predicted displacement,
$\displ$, are given for a finite number of receivers, $\Nrec$,
at locations $\{\rec[\mathbf{r},n]\}_{n=1}^{\Nrec}$.
Let us ignore model errors and assume that the measured
data is given by the model, $\displ$, depending on two parameters,
denoted by $\source$ and $\ParEar$. Here $\source$ corresponds to
the unknown source location, which can be modeled as deterministic or
stochastic depending on the approach to the source inversion problem,
and $\ParEar$ is a random nuisance parameter, corresponding to the
uncertain Earth material parameters. 
We assume that $\data$ is given by the model up to some additive
noise:  
\begin{align}
  \data(\rec[\mathbf{r},n],t_k) & = \displ(\rec[\mathbf{r},n],t_k;\source[\ast],\ParEar^\ast)
 + \WhiteNoise, && 
    \begin{array}{l}
      n=1,2,\dots,\Nrec,\\
      k=0,1,\dots,K,
    \end{array}
  \label{eq:data}
\end{align}
where $\WhiteNoise\sim\mathcal{N}(\mathbf{0},\sigma^{2}\Imatrix)$, independent
identically distributed (i.i.d.), and $\source[\ast]$ and
$\ParEar^\ast$ denote some fixed values of $\source$ and $\ParEar$,
respectively. 
We consider the random parameter, $\ParEar$, to consist of
the material triplet, $\left(\density,\comprv,\shearv\right)$ as a
random variable or field and all other parameters than $\ParEar$ and
$\source$ as given. 

The additivity assumption on the noise, while naive, can easily be
replaced by more complex, correlated, noise models without affecting the
usefulness or implementation of the MLMC approach described in
Section~\ref{sec:MLMC}. 

Let us denote a \QoItext~by $\QoI(\source,\ParEar)$. An
example of a seismic source inversion approach is finding the
location, $\source$, that yields the lowest expected value of $\QoI$,
i.e., the solution of 
\begin{equation}
\textrm{argmin}_{\source[\ast]} \E{\QoI(\source,\ParEar) \vert \source=\source[\ast]},
\end{equation}
where $\E{\cdot\vert\cdot}$ denotes conditional expectation of 
the first argument with respect to the second argument.
Both \QoItexts~investigated in this work can be used to construct likelihood
functions for statistical inversion, see
e.g.~\cite{bissiri2016general}. For instance, finding the source
location, $\source$, by maximizing the marginal likelihood, i.e., the
solution of 
\begin{equation}
\textrm{argmax}_{\source[\ast]} \E{\mathcal{L}(\source,\ParEar) \vert \source=\source[\ast]},
\end{equation}
where $\mathcal{L}$ is the likelihood function.

The two \QoItexts~below represent two different classes, 
which have many variations. In the present work we neither aim to 
compare the two \QoItexts~to each other, nor to choose between different
representatives of the two classes. Instead, we want to show the efficiency of 
MLMC applied to \QoItext~from both classes. Furthermore, in
the two \QoItexts~below, we assume that the discrete
time observations have been extended to a continuous function,
e.g., by linear interpolation between data points; we could 
also have expressed the \QoItext~in terms of discrete time observations.

\paragraph{$L^2$-based \QoItext} The first \QoItext~
studied in this work, denoted by $\QoI_{E}$, is based on the commonly used
$L^2$ misfit between predicted data $\displ$ and measured data
$\data$: 
\begin{align}
  \QoI_{E}(\source,\ParEar) & = 
  \dfrac{1}{\mathcal{T}}
  \int_{0}^{\mathcal{T}} \sum_{n=1}^{\Nrec}
    \left\vert\displ(\rec[\mathbf{r},n],t;\source,\ParEar)
   -\data(\rec[\mathbf{r},n],t)\right\vert^{2}\mathrm{d}t, 
  \label{eq:QoI_E}
\end{align}
where $\vert\cdot\vert$ is the Euclidean norm in $\rset^d$ and
$\mathcal{T}$ is the total simulation time. 
This quantity of interest is directly related to the seismic inversion
problem through its connection to the likelihood for
normally-distributed variables: 
\begin{align} 
  \label{eq:Lik}
  \mathcal{L}\big(\source,\ParEar\big)
  & \propto\exp(-\Theta(\source,\ParEar)) = 1-\Theta(\source,\ParEar)+\mathcal{O}(\Theta^{2}(\source,\ParEar)),
\end{align}
where
$\Theta(\source,\ParEar)
:=\dfrac{1}{2\sigma^{2}}\norm{\displ(\source,\ParEar)-\data}_{L^2(0,\mathcal{T}]}^{2}
=\dfrac{\mathcal{T}}{2\sigma^{2}}\QoI_{E}(\source,\ParEar)$. 

A drawback with the $L^2$ misfit function for full-waveform seismic
inversion, see e.g.~\cite{Engquist2018Analysis}, is the well-known
cycle skipping issue
which typically leads to many local optima and raises a substantial challenge to
subsequent tasks such as optimization and Bayesian inference. 

\paragraph{$W_2^2$-based \QoItext}

An alternative \QoItext~was introduced, in the setting of seismic
inversion, and analyzed
in~\cite{Engquist2014application,engquist2016optimal,Engquist2018Analysis,Engquist2018Application}, 
where it is shown to have several desirable properties which
$\QoI_E$ is lacking; in particular, in an idealized case, 
if one of the two waveforms is shifted in time, this \QoItext~is a convex
function of the shift; see Theorem~2.1 in~\cite{Engquist2014application} and 
the discussion in~\cite{engquist2016optimal}. 
This \QoItext~is based on the quadratic Wasserstein distance between
two probability density functions (PDFs), $\psi:X\rightarrow\rset^{+}$
and $\phi:X\rightarrow \rset^{+}$, which is defined as 
\begin{align} 
  W_{2}^{2}(\psi,\phi) & =
  \inf_{\mathcal{P}\in\mathcal{M}}
    \int_{X}|x-\mathcal{P}(x)|^{2}\psi(x)\,\mathrm{d}x, 
  \label{eq:Wass1}
\end{align}
where $\mathcal{M}$ is the set of all maps that rearrange the
PDF $\psi$ into $\phi$. 
When $X$ is an interval in $\rset$, an explicit form
\begin{align}
  W_{2}^{2}(\psi,\phi) & =
  \int_{0}^{1}|\Psi^{-1}(t)-\Phi^{-1}(t)|^{2}\mathrm{d}t
  \label{eq:Wass2}
\end{align}
exists, where $\Psi(\cdot)$ and $\Phi(\cdot)$ are the cumulative
distribution functions (CDFs) of $\psi$ and $\phi$ respectively.

How to optimally construct a~\QoItext~for seismic source-inversion based
on the $W_2^2$-distance, or similar distances, is an active research topic, 
and several recent papers discuss advantages and disadvantages of various 
approaches; see e.g.~\cite{Engquist2018Analysis,Metivier2018,Motamed_Appelo}. 
Here, we use one of the earliest suggestions, proposed 
in~\cite{Engquist2014application}.

To eliminate the scaling due to the length of the time interval we
make a change of variable $\tau=t/\mathcal{T}$ so that $X=[0,1]$ below.
Typically, the waveforms will not be PDFs, even in their component
parts. 
If we assume that $\psi$ and $\phi$ are two more general
one-dimensional functions, taking both positive and negative values in
the interval, then the
non-negative parts $\psi^{+}$ and $\phi^{+}$ and non-positive parts
$\psi^{-}$ and $\phi^{-}$ can be considered separately, and one can define 
\begin{align}
  \mathcal{W}(\psi,\phi) & := 
  W_2^2\left(
    \frac{\phi^-}{\int_{0}^1\phi^-},
    \frac{\psi^-}{\int_{0}^1\psi^-}
  \right) 
  + 
  W_2^2
  \left(
    \frac{\phi^+}{\int_{0}^1\phi^+}, 
    \frac{\psi^+}{\int_{0}^1\psi^+}
  \right).
  \label{eq:W}
\end{align}
To define the \QoItext~we sum $\mathcal{W}(\cdot,\cdot)$ applied to all
spatial components of the vector-valued $\displ$ and $\data$ in all
receiver locations, i.e.
\begin{align}
  \label{eq:QoI_W}
  \QoI_W(\source,\ParEar) & = 
  \sum_{n=1}^{\Nrec} \sum_{j=1}^d 
    \mathcal{W}\big(\,\displ_j(\rec[\mathbf{r},n],t(\tau);\source,\ParEar)\, , 
                \,\data_j(\rec[\mathbf{r},n],t(\tau))\,\big).
\end{align}

\begin{remark}[Assumption on alternating signs of $\displ$ and $\data$]
  \label{rem:alternating_signs}
  Note that the definition of $\QoI_W$ above requires all components
  of $\displ$ and $\data$ in all receivers to obtain both positive and
  negative values in the time interval $[0,\mathcal{T}]$
  for~\eqref{eq:W} to be well-defined. With Gaussian noise
  in~\eqref{eq:data}, the probability of violating this assumption is
  always positive, though typically too small to observe in practice if
  the observation interval and the receivers are properly set up. To
  complete the definition of $\QoI_W$, we may extend~\eqref{eq:W} by
  replacing
  $W_2^2\left(\phi^\ast/\int_0^1\phi^\ast,\psi^\ast/\int_0^1\psi^\ast\right)$
  by its maximal possible value, 1, whenever at least one of 
  $\phi^\ast$ and $\psi^\ast$ is identically 0.
  Note that this can lead to issues due to reduced regularity 
  beyond the loss of differentiability caused by the splitting into positive 
  and negative parts.
\end{remark}

\begin{figure}
  \centering
  \includegraphics[width=60mm]{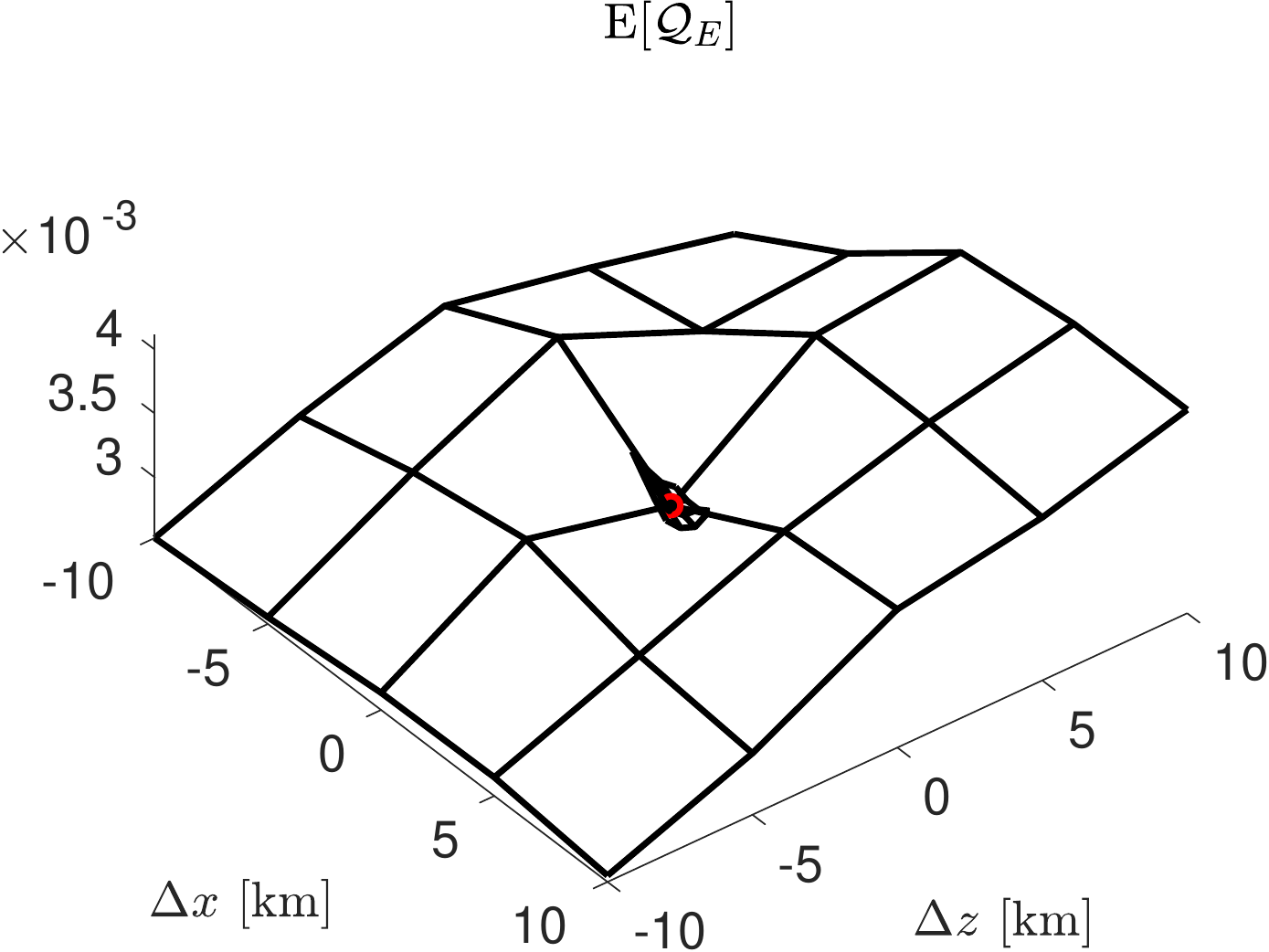}
  \includegraphics[width=60mm]{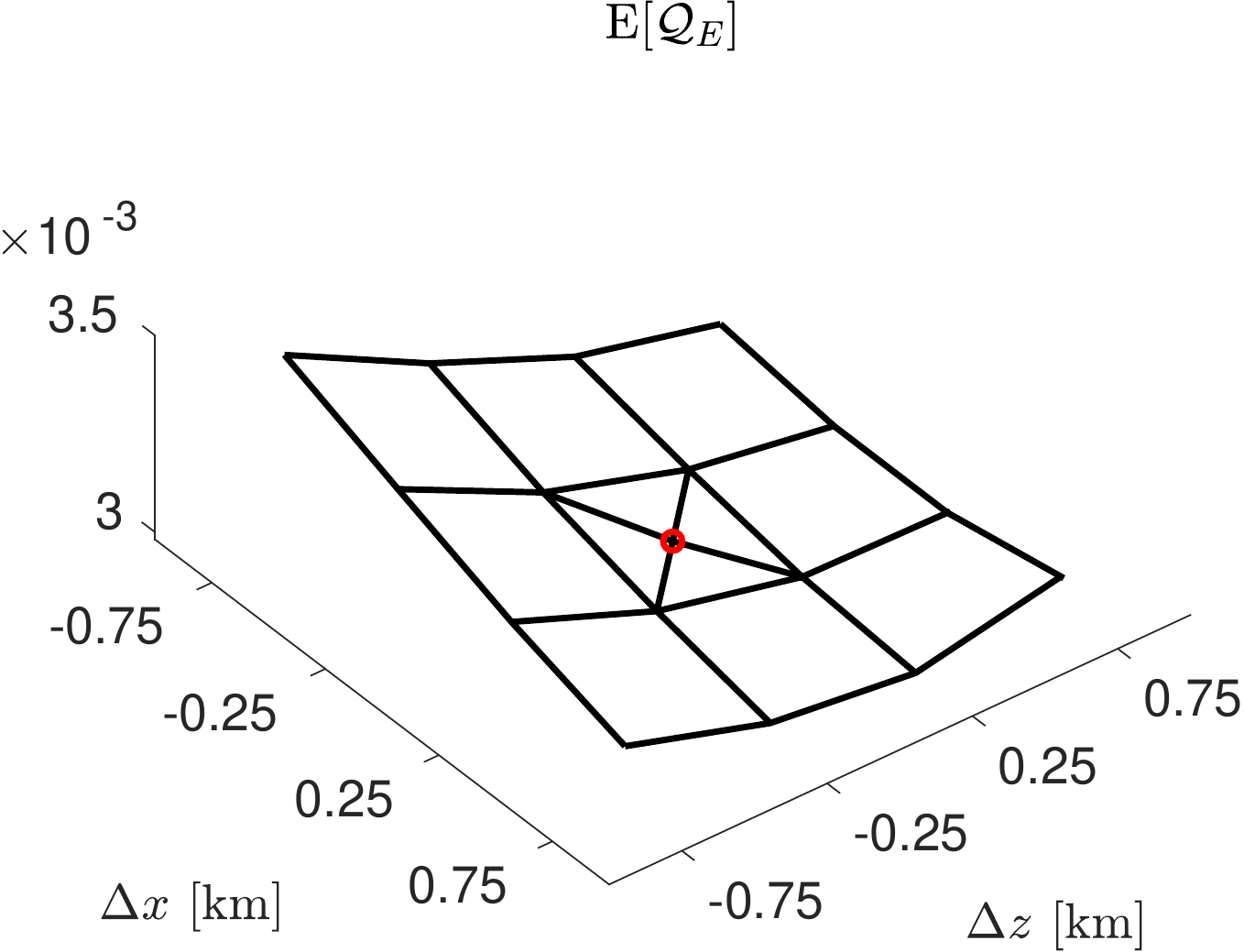}\\
  \vspace{0.05\textwidth}
  \includegraphics[width=60mm]{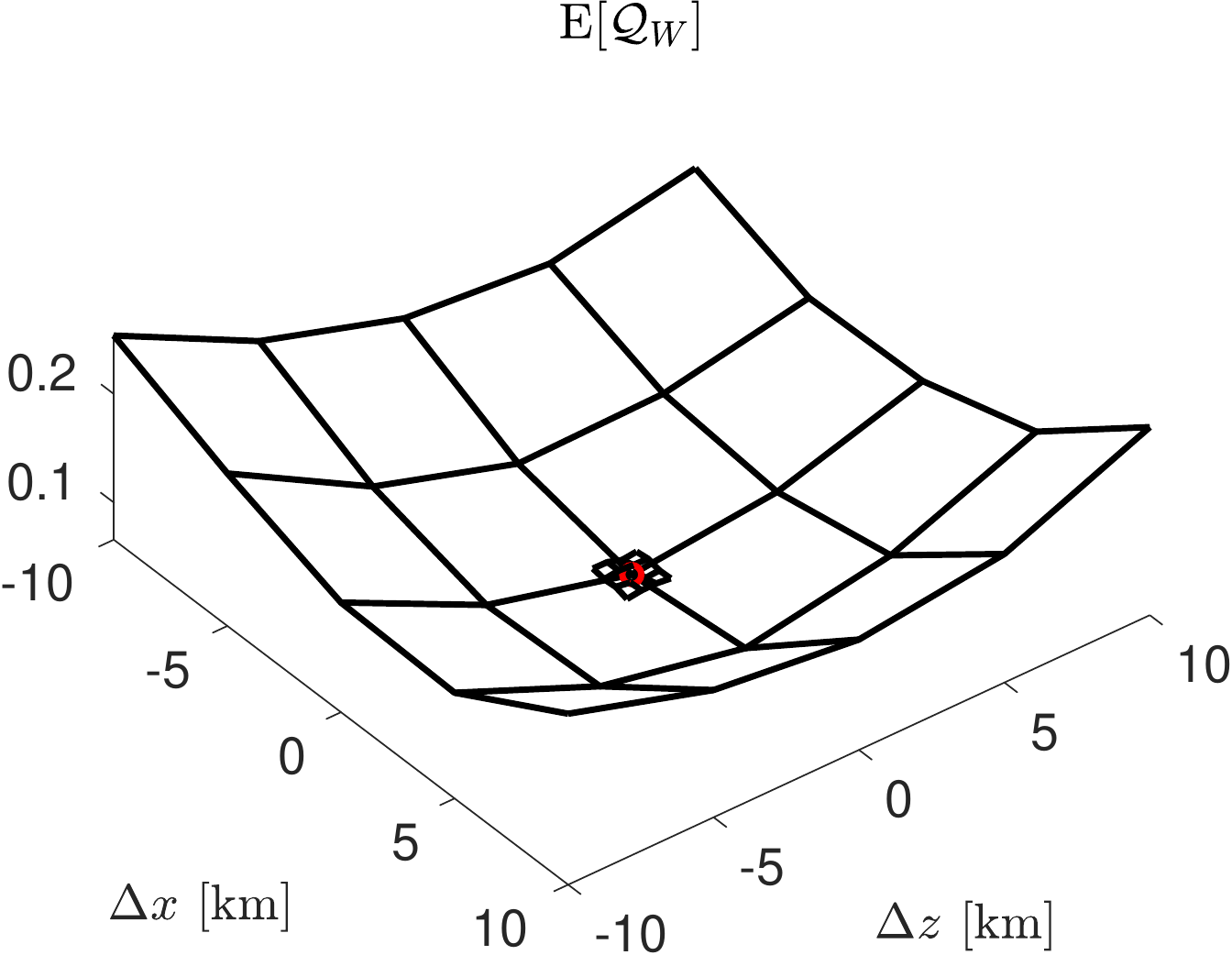}
  \includegraphics[width=60mm]{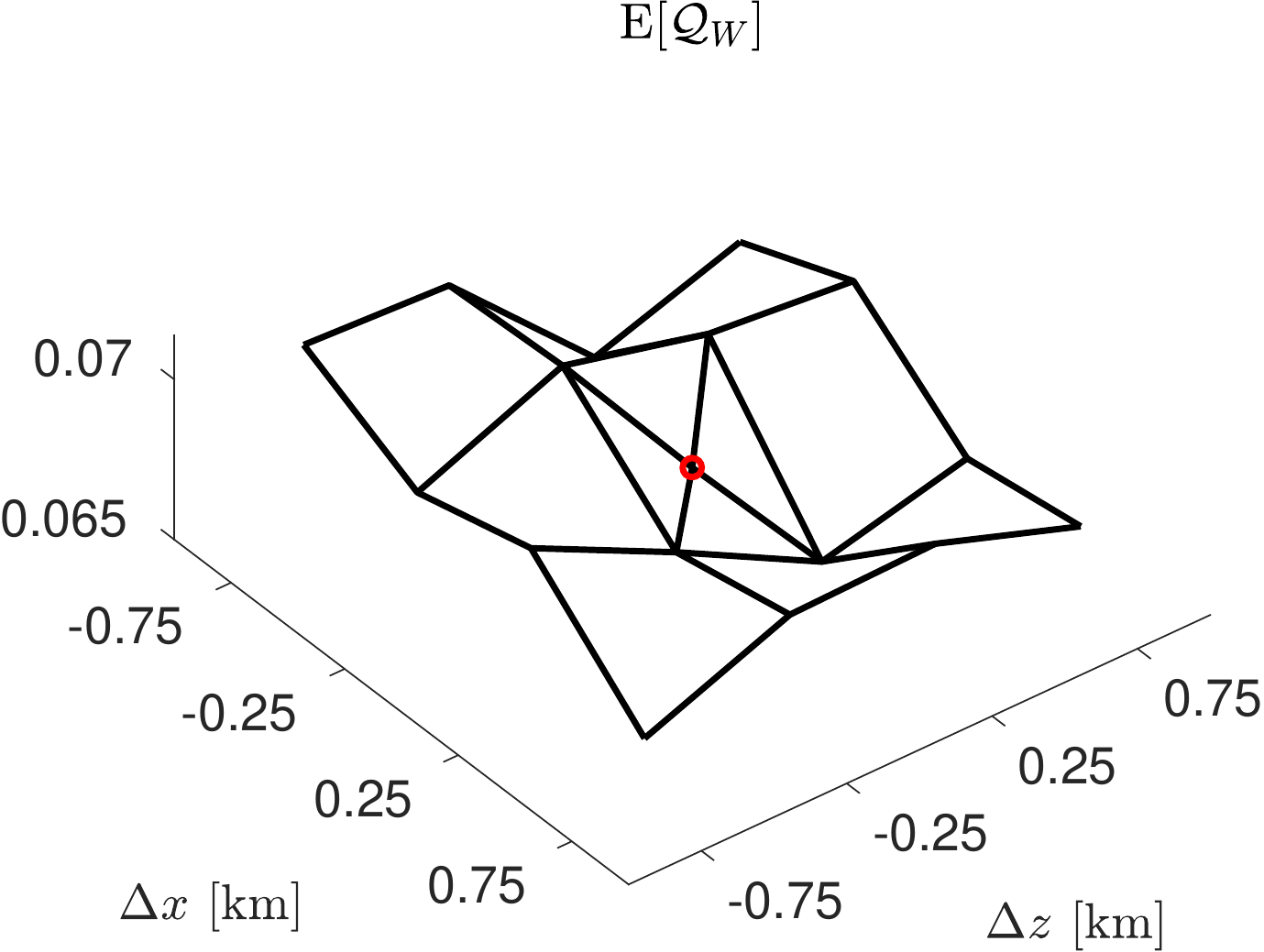}\\
  \vspace{0.05\textwidth}
  \includegraphics[width=60mm]{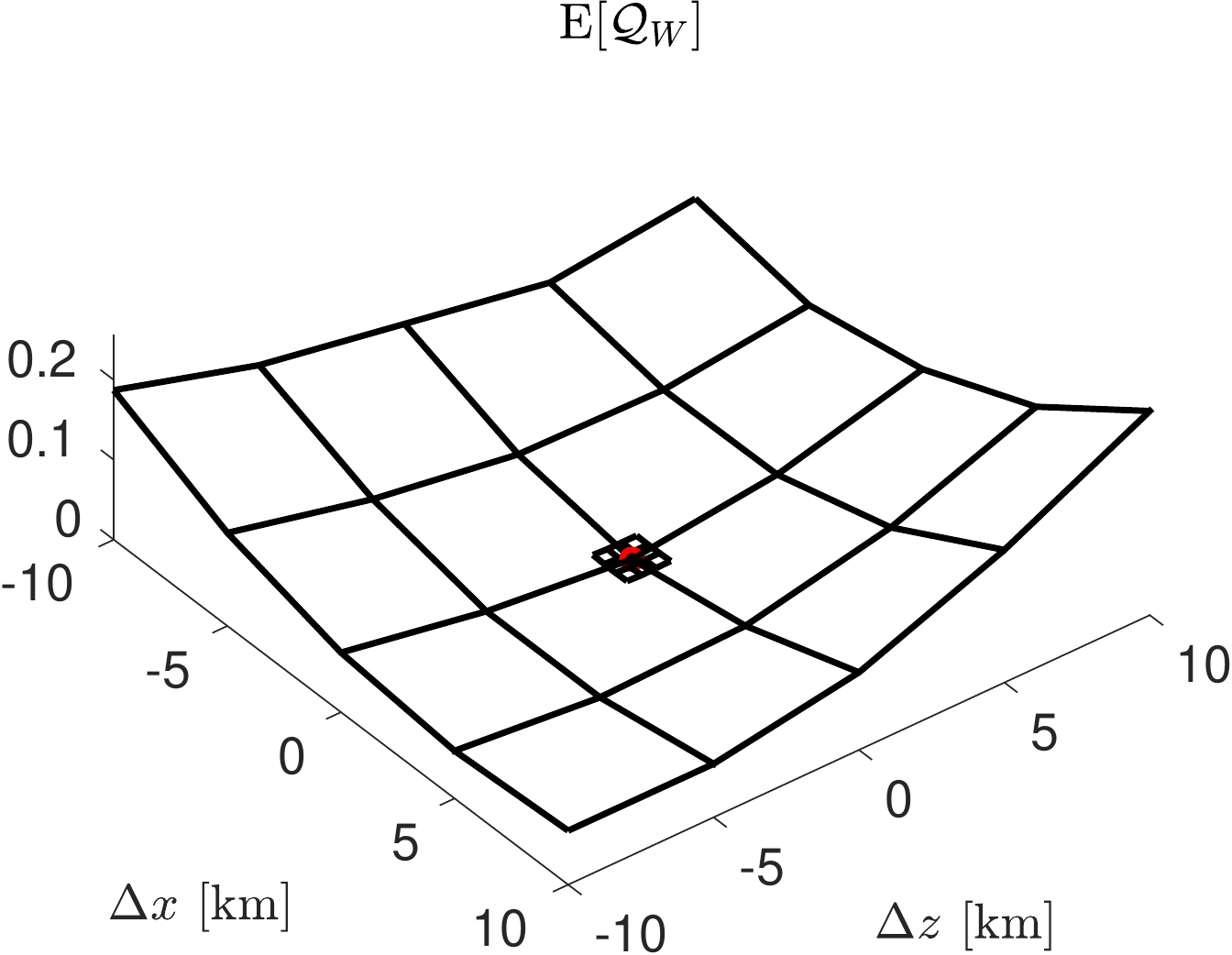}
  \includegraphics[width=60mm]{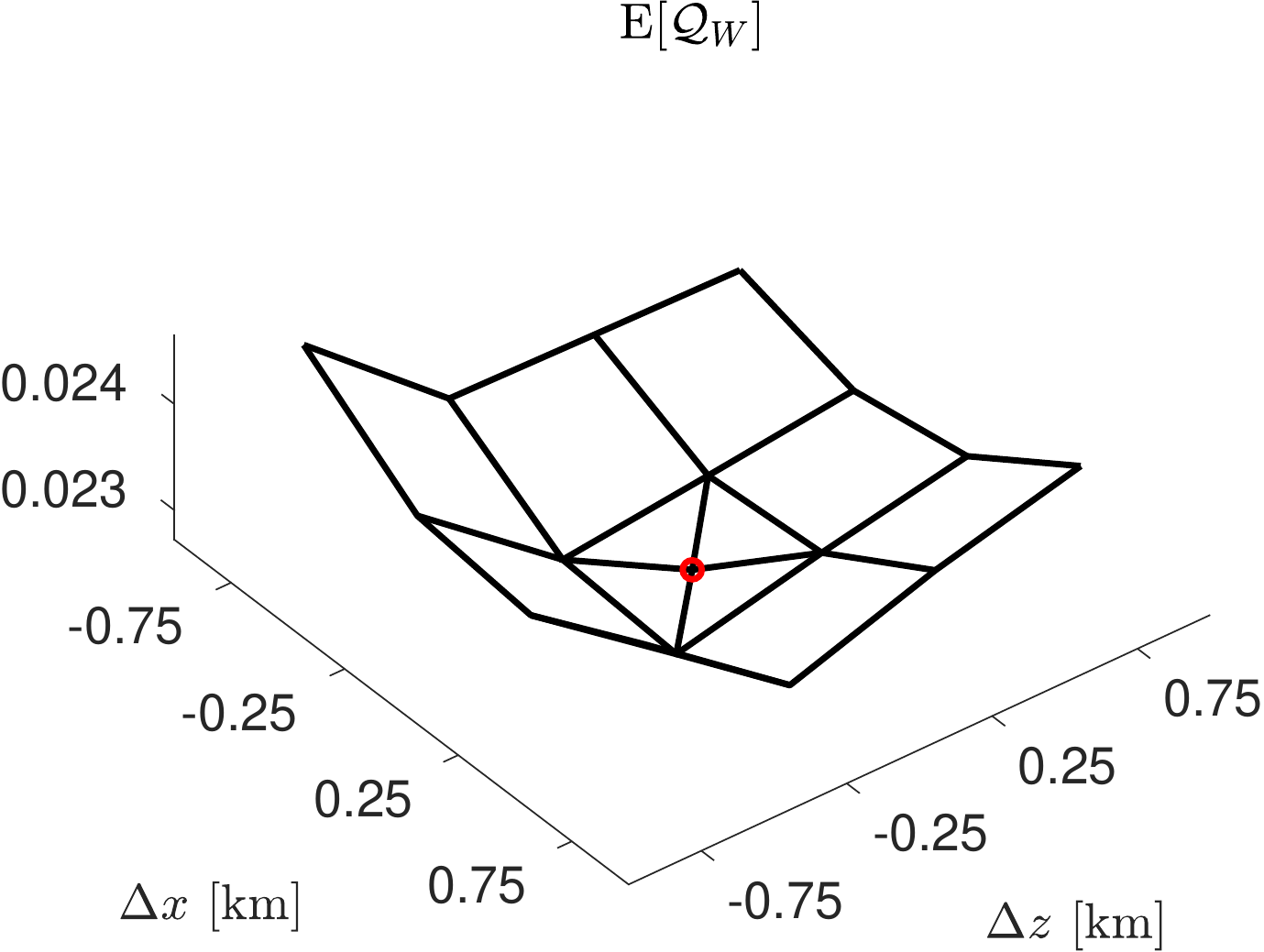}
  \caption{The expected value, $\E{\QoI(\source,\ParEar)}$, as a
    function of the shift between source location in the MLMC
    simulation and in the simulation used to generate the synthetic
    data; see Remark~\ref{rem:convexity}.
    \newline
    The top row shows $\QoI_E$, computed with synthetic data with
    additive noise, the middle row $\QoI_W$ with the same synthetic
    data, and the bottom row, $\QoI_W$ without noise added to the
    synthetic data.
    \newline
    The effect of adding or removing a noise of this level is not
    visible in $\QoI_E$, and therefore the corresponding figures of
    $\QoI_E$ without added noise are omitted.} 
  \label{fig:QoI}
\end{figure}

\begin{remark}[Use of $\QoI_W$ in source-inversion]
  \label{rem:convexity}
  The convexity of $\QoI_W$ with respect to time-shifts in signals
  is directly related to source inversion
  problems,~\cite{engquist2016optimal}, since
  perturbations of the source location approximately result in shifts 
  in the arrival times at the receivers. 
  However, convexity with respect to the source location 
  is not guaranteed.
  As an illustration, consider the case where the data $\data$ is
  synthetic data obtained from a computed approximation of
  $\displ(\cdot,\cdot;\source[\ast],\ParEar^\ast)$.
  Figure~\ref{fig:QoI} shows results from a numerical
  example very similar to the one in Section~\ref{sec:num}. In
  Figure~\ref{fig:QoI} 
  $\E{\QoI_E(\source,\ParEar)}$ and $\E{\QoI_W(\source,\ParEar)}$ are
  approximated for 41 different synthetic data, obtained by shifting
  the source location, $\source[\ast]$, while keeping the Earth
  material parameters fixed, as given in
  Table~\ref{tab:Earth_model_2D_synthetic} on page~\pageref{tab:Earth_model_2D_synthetic}. 
  In the figure, 
  $(\Delta x,\Delta z)=\Delta\mathbf{x}=\source[\ast]-\source$ is the
  shift of the source location in the synthetic data relative to the
  source location encoded in $\force$ when approximating Problem~\ref{prob:weak2}. 
  The left column shows a larger region around the point
  $\Delta\mathbf{x}=\mathbf{0}$, marked with a red circle, and the
  right gives a detailed view around
  $\Delta\mathbf{x}=\mathbf{0}$.
  The left column clearly illustrates a situation where the large-scale
  behavior of $\E{\QoI_W(\source,\ParEar)}$ is convex with respect to 
  $\Delta\mathbf{x}$, while $\E{\QoI_E(\source,\ParEar)}$ is not. In the 
  absence of noise, both \QoItexts~are convex in a small neighborhood around 
  the point $\Delta\mathbf{x}=\mathbf{0}$.
  Avoiding the non-convex behavior on the larger scale 
  significantly simplifies the source inversion problem. 
\end{remark}

\section{Computational Techniques}
\label{sec:comp}

In this section, we start with a concise description of
how Problem~\ref{prob:weak1} and Problem~\ref{prob:weak2} are
approximated numerically. The domain, $\domain$, is modified in two
steps: first, the finite Earth model is replaced by a half-plane in
two dimensions or a half-space in three dimensions, and second, this
semi-infinite domain is truncated, introducing absorbing boundary
conditions on the artificial boundaries. We also describe a
simplification of the stress tensor model~\eqref{eq:anisoHookeLaw} that
results in a viscoelastic stress tensor suitable for numerical 
implementation. Then, we proceed with providing computational
approximations of the \QoItext~given in Section~\ref{sec:QoI}. The
section ends with a summary of the MLMC algorithm for computing
the expected value of the \QoItext.

\subsection{Numerical Approximation of Initial Boundary Value Problem}
\label{sec:comp_weak}

A numerical approximation of Problem~\ref{prob:weak1}
or~\ref{prob:weak2} can either be achieved by (i)~approximating the
seismic wave propagation produced by a seismic event on the whole
Earth, or by (ii)~restricting the computational domain, $\domain$, to
a local region around the source and the receivers. In either case,
there are purpose-built software packages based on the Spectral
Element Method (SEM)~\cite{Komatitsch_Tromp1999,Dimitri1998} that
will be used in this paper. 
The MLMC method does not fundamentally depend on which
of the alternatives, (i) or (ii), that is used, or on the choice of
SEM over other approximation methods.
Indeed, an important advantage of MLMC, or more generally MC, methods
is that they are non-intrusive in the sense that they can
straightforwardly be applied by randomly sampling the Earth material
parameters and then executing any such publicly available simulation
code to compute the corresponding sample of the \QoItext.

In our numerical example, we choose alternative (ii), and proceed in
two steps: first, we approximate the Earth locally by a half-plane, in
a two-dimensional test case, or by a half-space, in the full
three-dimensional problem; second, the half-plane or half-space is
truncated to a finite domain, where absorbing boundary conditions
(ABC) are introduced on the artificial boundaries to mimic the
absorption of seismic energy as the waves leave the region around the
receivers.  
The variational equations~\eqref{eq:weak_eq_motion}
and~\eqref{eq:weak_eq_motion_attenuation} now contain a non-vanishing
boundary term, corresponding to the part of the boundary, $\artbound$, 
where absorbing boundary conditions apply.  
We use a perfectly matched layer (PML) approximation of the ABC,
introduced in~\cite{BERENGER1994185} and used in many fields; see
e.g.~\cite{xie2014PML,KomatitschPerfectlyMatched} in the context of
seismic wave propagation. 
However, in the absence of true PML, see~\cite{Xie2016_OceanAcustic},
for Problem~\ref{prob:weak2} which has attenuating Earth material
properties, in practice we choose the truncation of $\domain$ such
that $\artbound$ is far enough from all receivers to guarantee that no
reflected waves reach the receivers in 
the time interval $[0,\mathcal{T}]$, given the maximal wave speeds
allowed by the range of uncertainties~\eqref{eq:RandomParCond1}.

To apply SEM, first, a semi-discrete version of the variational
equation is introduced by discretizing space and introducing a finite
dimensional solution space where the solution at time $t$ can be
represented by a finite vector, $\mathbf{S}(t)$. Then, the time evolution of the
SEM approximation, $\mathbf{S}(t)$, of the seismic wavefield solves an initial
value problem for the second order ordinary differential equation
(ODE) in time
\begin{align}
  \label{eq:SemidisWF}
   M\ddot{\mathbf{S}} + C\dot{\mathbf{S}} +K\mathbf{S} & = F, && 0< t \leq \mathcal{T}, 
\end{align}
where $M$ is the mass matrix, $C$ is the global absorbing boundary
matrix, $K$ the global stiffness matrix, and $F$ the source term.

To get the semi-discrete form, $\domain$ is divided into
non-overlapping elements of maximal size $\dx[]$, similarly to what is
done when using a standard finite-element method. 
Quadrilateral elements are used in two space dimensions and hexahedral
in three. 
Each element is defined in terms of a number, $n_{c}$, of control
points and an equal number of shape functions which describe the
isomorphic mapping between the element and reference square or cube.
The shape functions are products of Lagrange polynomials of low
degree. 
In the remainder, we assume that no error is introduced by the
representation of the shape of the elements, which is justified by the
very simple geometry of the test problem in Section~\ref{sec:num}.
The displacement field on every element is approximated by a Lagrange
polynomial of higher degree, $N_{l}$, and the approximation of the
variational form~\eqref{eq:weak_eq_motion}
or~\eqref{eq:weak_eq_motion_attenuation}, including the artificial 
boundary term, over an element is based on the Gauss-Lobatto-Legendre
integration rule on the same $N_{l} + 1$ points used for Lagrange
interpolation;
this choice leads to a diagonal mass matrix, $M$, which is beneficial
in the numerical stepping scheme.
More details on the construction of these matrices and the source term
can be found in~\cite{Dimitri1998}.

The initial value problem for the ODE~\eqref{eq:SemidisWF} is
approximately solved by introducing a discretization of the time
interval, and by applying a time-stepping method. Among multiple available 
choices, this work uses the second-order accurate explicit
Newmark-type scheme; see for example Chapter~9
in~\cite{Hughes2000}. It is a conditionally stable scheme and the
associated condition on the time step leads to
$\dt[]\leq c\,\dx[]N_{l}^{-2}$, for uniform spatial discretizations. 
Thanks to the diagonal nature of $M$ and the sparsity of $C$ and $K$, the cost
per time step of the Newmark scheme is proportional to the number of
unknowns in $\mathbf{S}(t_j)$, and thus to $\dx[]^{-d}$, and since the number of time
steps is inversely proportional to $\dt[]\propto\dx[]$ the total work
is proportional to $\dx[]^{-(d+1)}$.

Determining $\dt[]$ by the stability constraint,
$\dt[]\propto\dx[]$, we expect the second order accuracy of the
Newmark scheme to asymptotically be the leading order error term as
$\dx[]\to0$, assuming sufficient regularity of the true solution. 

\subsection{Computational Model of Seismic Attenuation}
\label{sec:comp_attenuation}

Approximately solving Problem~\ref{prob:weak2}, as it is stated in
Section~\ref{sec:attenuation}, is very difficult since the stress~$\stress$
in~\eqref{eq:anisoHookeLaw} at time $t$ depends on the entire solution
history $\displ(\cdot,t')$ for $-\infty<t'\leq t$, or in practice for
$0\leq t'\leq t$ since the displacement is assumed to be constant up
to time 0. 
Even in a discretized form in an explicit time stepping scheme,
approximately updating the stress according
to~\eqref{eq:anisoHookeLaw} would require storing the strain history
for all previous time steps in every single discretization point where
$\stress$ must be approximated, requiring unfeasible amounts of
computer memory and computational time. 
Therefore, the model of the viscoelastic properties of the medium is
often simplified to a generalized Zener model using a series of
``standard linear solid" (SLS) mechanisms. The present work uses the 
implementation of the generalized Zener model in {\tt SPECFEM2D}.
Below, we will briefly sketch the simplification
of~\eqref{eq:anisoHookeLaw}. For a more detailed 
description, we refer readers to~\cite{LiAnKa76,Zhang2016,Komatitsch_Tromp1999,MoczoRobEisner2007,Carcione2014}.

The integral in the stress-strain 
relation~\eqref{eq:anisoHookeLaw} can be expressed as a
convolution in time by defining the relaxation tensor $\elasticmod$ to
be zero in $\domain\times\rset^-$, 
i.e. $\elasticmod(\cdot,t) = \tilde{\elasticmod}(\cdot,t)H(t)$, where
$H(t)$ is the Heaviside function. That is,
\begin{multline}
  \stress\left(\xv,t;\left\{\grad\displ\right\}_0^t\right)  = \int_{-\infty}^{\infty}
  \elasticmod(\xv, t - t^{'}) \doubledot 
  \partial_{t}\strain(\grad\displ(\xv,t^{'})) \,dt^{'} 
   = \left(\elasticmod \ast \partial_{t}\strain\right) (\xv,t)
  \\
   = \left(\partial_{t} \elasticmod \ast \strain\right) (\xv,t), 
  \label{eq:anisoHookeLawConv}
\end{multline}
which, as discussed in~\cite{Emmerich1987}, can be formulated in the
frequency domain as 
\begin{align}
  \label{eq:anisoHookeLawFreqD}
  \widehat\stress(\xv,\omega) & = 
  \widehat{\mathbf{M}}(\xv,\omega)\widehat\strain(\xv,\omega),
\end{align}   
where $\mathbf{M}(\cdot,t)=\partial_t\elasticmod(\cdot,t)$.
In seismology, it has been observed,~\cite{Carcione1988},
that the so-called quality factor
\begin{align}
  \label{qualityFactor}
  Q(\omega) & = 
  \frac{Re\left(\widehat{\mathbf{M}}(\cdot,\omega)\right)}
       {Im\left(\widehat{\mathbf{M}}(\cdot,\omega)\right)},
\end{align}
is approximately constant over a wide range of frequencies. 
This $Q$ is an intrinsic property of the Earth material that describes
the decay of amplitude in seismic waves due to the loss of energy to
heat, and its impact on $\mathbf{M}$ is explicitly given
  in equation (5) of~\cite{Emmerich1987}. 
This observation allows modeling $\widehat{\mathbf{M}}(\cdot,\omega)$ in
the frequency domain through a series of a number, $B$, of SLS.
Then the stress-strain relation~\eqref{eq:anisoHookeLawConv} can be
approximated as
\begin{align}
  \label{eq:anelasticSLSHL}
  \stress & = \elasticmod^{U}\doubledot\strain - \sum_{b=1}^{B}\mathbf{R}^b,
\end{align} 
where $\elasticmod^{U}$ is the unrelaxed viscoelastic fourth order
tensor, which for an isotropic Earth model is defined by $\Lamefp$ and
$\Lamesp$ or equivalently $\comprv$ and $\shearv$. The relaxation
functions $\mathbf{R}^b$ for each SLS satisfy initial value problems
for a damping ODE; by the non-linear optimization approach given
in~\cite{Blanc01042016}, implemented in \specfem and used in this
work, the ODE for each $\mathbf{R}^b$ is determined by two
parameters: the quality factor, $Q$, and the number of SLS, $B$. 

\subsection{Quantities of Interest}
\label{sec:comp_QoI} 

The two \QoItext, $\QoI_E$ and $\QoI_W$ defined in~\eqref{eq:QoI_E}
and~\eqref{eq:QoI_W} respectively, are approximated from discrete
time series $\{\displ(\rec,t_j^s;\cdot,\cdot)\}_{j=0}^J$ and
$\{\data(\rec,t_k^d)\}_{k=0}^K$. 
The data observation times $0=t_0^d<\ldots<t_{K}^d=\mathcal{T}$, 
are considered given and fixed, and with realistic frequencies of the
measurements, $100-200\,\mathrm{Hz}$, it is natural to take smaller
time steps in the time discretization
$0=t_0^s<\dots<t_J^s=\mathcal{T}$ of the numerical approximation of
$\displ$ and we assume that
$\left\{t_k^d\right\}_{k=0}^K\subseteq\left\{t_j^s\right\}_{j=0}^J$. Furthermore, 
the time discretization is assumed to be
characterized by one parameter $\dt[]$, e.g., the constant time step size of a 
uniform discretization.
Since we defined both \QoItext~as sums over all receivers and all
components of the vector-valued functions $\displ$ and $\data$ in the
receivers, it is sufficient to describe the approximation in the case of
two scalar functions $\phi^s$, representing the simulated
displacement, and $\phi^d$, representing the observed data. 
Here, we define the function $\phi^d$ as the piecewise linear
interpolation in time of the data points.

\paragraph{Approximation of $\QoI_E$}

The time integral in~\eqref{eq:QoI_E} is approximated by
the Trapezoidal rule on the time discretization of the numerical
simulation. With the assumption that 
$\{t_k^d\}_{k=0}^K\subseteq\{t_j^s\}_{j=0}^J$ and the definition of
$\phi^d$ as the piecewise linear interpolation in time of the data
points, only the discretization of $\phi^s$ contributes to the
error. The numerical approximation of $\QoI_E$ will then have an
$\mathcal{O}(\dt[]^2)$ asymptotic error, provided that $\phi^s$ is
sufficiently smooth.

\paragraph{Approximation of $\QoI_W$}

The approximation of~\eqref{eq:QoI_W}, through~\eqref{eq:W}, requires
both positive and negative values to be attained in all components of
$\{\displ(\rec,t_j^s;\cdot,\cdot)\}_{j=0}^J$ and 
$\{\data(\rec,t_k^d)\}_{k=0}^K$ in all receivers; see
Remark~\ref{rem:alternating_signs} on
page~\pageref{rem:alternating_signs}. Assuming that this holds we 
approximate the $W_2^2$-distance between the normalized non-negative
parts of $\phi^s$ and $\phi^d$; we treat the non-positive
analogously. To this end, the zeros of $\phi^s$ and $\phi^d$, denoted
$\{z_j^s\}_{j=0}^{\hat{J}}$ and $\{z_k^d\}_{k=0}^{\hat{K}}$ respectively, 
are approximated by linear interpolation, and they are included in the
respective time discretizations, generating
$\{t_j^s\}_{j=0}^{J^\ast}=\{t_j^s\}_{j=0}^J\bigcup\{z_j^s\}_{j=0}^{\hat{J}}$
and 
$\{t_k^d\}_{k=0}^{K^\ast}=\{t_k^d\}_{k=0}^K\bigcup\{z_k^d\}_{k=0}^{\hat{K}}$,
and thus $\{\phi_j^{s,+}\}_{j=0}^{J^\ast}$ and
$\{\phi_k^{d,+}\}_{k=0}^{K^\ast}$ are obtained. Then the corresponding
values of the CDFs, $\{\Phi_j^{s,+}\}_{j=0}^{J^\ast}$ and
$\{\Phi_k^{d,+}\}_{k=0}^{K^\ast}$ are approximated by the Trapezoidal
rule, followed by normalization. Finally, the inverse of $\Phi^{d,+}$
in $\{\Phi_j^{s,+}\}_{j=0}^{J^\ast}$, i.e. 
$\left\{\left[\Phi^{d,+}\right]^{-1}\left(\Phi_j^{s,+}\right)\right\}_{j=0}^{J^\ast}$, 
is approximated by linear interpolation, and analogously for the inverse
of $\Phi^{s,+}$ in $\{\Phi_k^{d,+}\}_{k=0}^{K^\ast}$, before the
integral in~\eqref{eq:Wass2} is approximated by the Trapezoidal rule
on the discretization 
$\{\Phi_j^{s,+}\}_{j=0}^{J^\ast}\bigcup\{\Phi_k^{d,+}\}_{k=0}^{K^\ast}$
of $[0,1]$.
These steps combined lead to an $\mathcal{O}(\dt[]^2)$ asymptotic
error, provided that $\phi^s$ is sufficiently smooth.

\begin{remark}[Cost of approximating $\mathcal{Q}_W$]
  In the present context, passive source-inversion with a small 
  number of receivers, in the order of 10, the cost of approximating 
  $\mathcal{Q}_W$ from $\displ$ and $\data$ is negligible compared 
  to that of computing the approximation of $\displ$ itself.
  This might not always be the case in other contexts where one has
  data from a large number of receivers, as could be the case e.g. in
  seismic imaging for oil and gas exploration.
\end{remark}

\begin{remark}[On expected weak and strong convergence rates]
  The deterministic rate of convergence, as $\dt[]^s\to 0$, of both 
  $\QoI$ approximations, is two. Here, $\dt[]^s$ is identical to the time
  step of the underlying approximation method of
  Problem~\ref{prob:weak1} or~\ref{prob:weak2}. 
  In the numerical approximation of $\displ$, with the second order
  Newmark scheme in time, the asymptotic convergence rate is also two at
  best, which holds if the solution is sufficiently regular.
  By our assumptions on the random fields satisfying the same
  regularity as in the deterministic case and being uniformly bounded,
  from above and away from zero from below, in the physical domain
  and with respect to random outcomes, 
  we expect both the weak and the
  strong rates of convergence to be the same as the deterministic
  convergence rate of the numerical approximation; i.e. at best two,
  asymptotically as $\dt[]\propto\dx[]$ goes to zero.
\end{remark}

\subsection{MLMC Algorithm}
\label{sec:MLMC}

Here, we summarize the MLMC algorithm introduced by
Giles~\cite{giles08} and independently, in a setting further from
the one used here, by Heinrich in~\cite{heinrich01}, which has since
become widely used~\cite{giles_AcNum}. 

MLMC is a way of reducing the computational cost of standard MC, for
achieving a given accuracy in the estimation of the expected value of some \QoItext, in
situations when the samples in the MC method are obtained by numerical
approximation methods characterized by a refinement parameter,
$\meshparam$, controlling both the accuracy and the cost.

\paragraph{Goal}
We aim to approximate the expected value of some \QoItext, $\E{\QoI}$, 
by an estimator $\Est$, with the accuracy requirement that
\begin{align}
  \vert\E{\QoI}-\Est\vert & \leq \tol, 
  && \text{with probability $1-\confidence$, for $0<\confidence\ll 1$,}
  \label{eq:goal}
\end{align}
where $\tol>0$ is a user-prescribed error tolerance.
To this end we require
\begin{subequations}
  \label{eq:splitting}
  \begin{align}
    \label{eq:bias}
    \left|\E{\QoI-\Est}\right| & \leq (1-\splitting)\tol \\
    \intertext{and}
    \label{eq:staterr}
    \prob{\left|\E{\Est}-\Est\right| > \splitting\tol} & \leq
                                                         \confidence,
  \end{align}  
\end{subequations}
for some $0<\splitting<1$, which we are free to choose.

\paragraph{Assumptions on the Numerical Approximation Model}
Consider a sequence of discretization-based approximations
of $\QoI$ characterized by a refinement parameter,
$\{\meshparam\}_{\ell=0}^\infty$. Let $\QoI_\ell(\ParEar)$ denote the
resulting  approximation of $\QoI$, using the refinement parameter
$\meshparam$ for an outcome of the 
random variable $\ParEar$. In this work, we consider successive
halvings in the refinement parameter, $\meshparam=2^{-\ell}$, which in 
Section~\ref{sec:comp_weak} corresponds to the spatial mesh size and
the temporal step size, $\dx=\dx[0]\meshparam$ and
$\dt=\dt[0]\meshparam$, respectively. 
We then make the following assumptions on how the cost and accuracy of
the numerical approximations depend on $\meshparam$.
We assume that the work per sample of $\QoI_\ell$, denoted $\work_\ell$,
depends on $\meshparam$ as
\begin{subequations}
  \label{eq:MC_models}
  \begin{align}
    \label{eq:work_model}
    \work_\ell & \propto\meshparam^{-\gamma}, 
    && \gamma>0,\\
    \intertext{and the weak order of convergence is $\Ow$, so that we 
     can model,}
    \label{eq:bias_model}
    \vert\E{\QoI-\QoI_\ell}\vert & = \Cw\meshparam[\ell]^\Ow, 
    && \Cw,\Ow>0,\\
    \intertext{and that the variance is independent of the
    refinement level}
    \label{eq:var_model_MC}
    \var{\QoI_\ell} & =\vg,
    && \vg>0.
  \end{align}
\end{subequations}

\paragraph{Standard MC estimator}
For $\nrs$ i.i.d. realizations
of the parameter,  $\{\ParEar_n\}_{n=1}^\nrs$, the unbiased MC
estimator of $\E{\QoI_\ell(\ParEar)}$ is given by 
\begin{align}
  \EstMC & = \frac{1}{\nrs}\sum_{n=1}^\nrs \QoI_\ell(\ParEar_n).
  \label{eq:MC_est}
\end{align}
For $\EstMC$ to satisfy~\eqref{eq:splitting} we require 
\begin{align}
  \label{eq:bias_requirement}
  \left|\E{\QoI-\QoI_\ell}\right| & \leq (1-\splitting)\tol,
  \intertext{which according to the model~\eqref{eq:bias_model} becomes}
  \label{eq:bias_requirement_model}
  \Cw\meshparam[\ell]^\Ow  & \leq (1-\splitting)\tol.
\end{align}
For any fixed $\meshparam$ such that $\Cw\meshparam[\ell]^\Ow<\tol$,
the value of the splitting parameter, $\splitting$, is implied by
replacing the inequality in~\eqref{eq:bias_requirement_model} by equality and solving
for $\splitting$, giving
\begin{align}
  \splitting & = 1-\frac{\Cw\meshparam[\ell]^\Ow}{\tol}.
\end{align}
Thus, the model for the bias tells us how large of a statistical error
we can afford for the desired tolerance, $\tol$.
By the Central Limit Theorem, $\EstMC$ properly rescaled converges in
distribution,
\begin{align}
  \label{eq:CLT}
  \frac{\sqrt{\nrs}\left(\E{\QoI_\ell}-\EstMC\right)}{\sqrt{\vg}} 
  & \Rightarrow \mathcal{N}(0,1), 
  && \text{as $\nrs\to\infty$},
\end{align}
where $\mathcal{N}(0,1)$ is a standard normal random variable with CDF
$\Phi_{\mathcal{N}(0,1)}$. 
Hence, to satisfy the statistical error constraint~\eqref{eq:staterr},
asymptotically as $\tol\to 0$, we require 
\begin{align}
  \label{eq:var_requirement}
  \frac{\vg}{ \nrs}
  & \leq
    \left(\frac{\splitting \tol}{\confpar} \right)^2,
\end{align}
where $\confpar$ is the confidence parameter corresponding to a
$1-\confidence$ confidence interval,
i.e. $\Phi_{\mathcal{N}(0,1)}(\confpar)=1-\confidence/2$.

The computational work of generating $\EstMC$ is
\begin{align*}
  \work_{\mathrm{MC}} & \propto \nrs\work_\ell.
\end{align*}
For asymptotic analysis, assume that we can choose
$\meshparam$ by taking equality in~\eqref{eq:bias_requirement_model}
and $\nrs$ by taking equality in~\eqref{eq:var_requirement}; then we get 
the asymptotic work estimate
\begin{align}
  \label{eq:complexity_MC}
  \work_{\mathrm{MC}} & \propto 
    \frac{\tol^{-\left(2+\gamma/\Ow\right)}}
         {\splitting^2\left(1-\splitting\right)^{\gamma/\Ow}}.
\end{align}
For any fixed choice of $\splitting$, the computational complexity of
the MC method is $\tol^{-\left(2+\gamma/\Ow\right)}$. Minimizing the
right hand side in~\eqref{eq:complexity_MC} with respect to $\splitting$
gives the asymptotically optimal choice
\begin{align}
  \label{eq:optimal_splitting_MC}
  \splitting & = \left(1+\frac{\gamma}{2\Ow}\right)^{-1}\in(0,1).
\end{align}

\paragraph{MLMC estimator}

The work required to meet a given accuracy by standard MC can be
significantly improved by systematic generation of control variates
given by approximations corresponding to different mesh sizes.
In the standard MLMC approach, we use a whole hierarchy of $L+1$ 
meshes defined by decreasing mesh sizes $\{\meshparam\}_{\ell=0}^L$ 
and the telescoping representation of the expected value of the finest
approximation, $\QoI_L$, 
\begin{align*}
  \E{\QoI_L} & =
  \E{\QoI_0} + \sum_{\ell=1}^L \E{\QoI_{\ell}-\QoI_{\ell-1}},
\end{align*}
from which the MLMC estimator is obtained by approximating the
expected values in the telescoping sum by sample averages as
\begin{align}
  \label{eq:MLMC_estimator}
  \EstMLMC & = \frac{1}{\nrs_0}\sum_{n=1}^{\nrs_0}\QoI_0(\ParEar_{0,n}) +
  \sum_{\ell=1}^L\frac{1}{\nrs_\ell}\sum_{n=1}^{\nrs_\ell}
  \left(\QoI_{\ell}(\ParEar_{\ell,n})-\QoI_{\ell-1}(\ParEar_{\ell,n})\right),
\end{align}
where  $\{\ParEar_{\ell,n}\}_{n=1,\dots,\nrs_\ell}^{\ell=0,\dots,L}$
denote i.i.d. realizations of the mesh-independent random variables. 
Note that the correction terms
\begin{align}
  \label{eq:Delta}
  \Delta\QoI_{\ell}(\ParEar_{\ell,n}) & = 
  \QoI_{\ell}(\ParEar_{\ell,n})-\QoI_{\ell-1}(\ParEar_{\ell,n})
\end{align}
are evaluated with the same outcome of $\ParEar_{\ell,n}$ in both the
coarse and the fine mesh approximation. This means that
$\var{\Delta\QoI_{\ell}}\to 0$, as $\ell\to\infty$, provided that the
numerical approximation $\QoI_{\ell}$ converges strongly. Introducing
the notation 
\begin{align}
  \label{eq:Var_ell}
  V_\ell & = 
  \begin{cases} 
    \var{\QoI_0}, & \ell=0,\\
    \var{\QoI_\ell - \QoI_{\ell-1}}, & \ell>0,
  \end{cases}
\end{align}
and assuming a strong convergence rate $\Os/2$ we model
\begin{align}
  \label{eq:var_model_MLMC}
  V_\ell  &= \Cs \meshparam^{\Os}, && \text{ for } \ell > 0.
\end{align}
Note that while this holds asymptotically as $\ell\to\infty$, by the
definition of strong convergence, this model may be inaccurate for
small $\ell$, corresponding to coarse discretizations. However, it
suffices for an asymptotic work estimate.

The computational work needed to generate $\EstMLMC$ is
\begin{align}
  \label{eq:work_sum}
  \work_\mathrm{MLMC} & = \sum_{\ell=0}^{L} \nrs_\ell \work_\ell,
\end{align}
where we now assume that~\eqref{eq:work_model} also holds for the cost
of generating $\Delta\QoI_\ell$.
In order for $\EstMLMC$ to satisfy~\eqref{eq:goal}, we fix 
$\splitting\in(0,1)$ and require $\EstMLMC$ to satisfy the bias
constraint~\eqref{eq:bias_requirement} and, 
consequently~\eqref{eq:bias_requirement_model}, on the finest
discretization, $\ell=L$, leading to 
\begin{align}
  \label{eq:finest_mesh}
  \meshparam[L] & = \left(\frac{(1-\splitting)\tol}{\Cw}\right)^{1/\Ow},
\end{align}
and we also require it to satisfy the statistical error
constraint~\eqref{eq:staterr}. In the MLMC
context,~\eqref{eq:staterr} is approximated by the bound 
\begin{align}
  \label{eq:var_requirement_MLMC}
  \sum_{\ell=0}^L \frac{V_\ell}{\nrs_\ell}
  & \leq
    \left(\frac{\splitting \tol}{\confpar} \right)^2
\end{align}
on the variance of $\EstMLMC$.
Enforcing~\eqref{eq:staterr} through this bound is justified 
asymptotically, as $\tol$ converges to 0, by a Central
Limit Theorem for MLMC estimators if for example $\Os>\gamma$; see Theorem 1.1
in~\cite{Hoel_CLT_MLMC}.
Given $L$ and $\splitting$, minimizing the work~\eqref{eq:work_sum}
subject to the constraint~\eqref{eq:var_requirement_MLMC} leads to the
optimal number of samples per level in $\EstMLMC$,
\begin{align}
  \label{eq:opt_nrs}
  \nrs_\ell & = \left( \frac{\confpar}{ \splitting\tol} \right)^2
  \sqrt{\frac{V_\ell}{\work_\ell}} \sum_{\ell=0}^L \sqrt{\work_\ell V_\ell}.
\end{align}
Substituting this optimal $\nrs_\ell$ in the total work \eqref{eq:work_sum} yields:
\begin{align}
  \label{eq:total_work_MLMC}
  \work_\mathrm{MLMC} & =
  \left( \frac{\confpar}{ \splitting\tol} \right)^2
  \left( \sum_{\ell=0}^L \sqrt{\work_\ell V_\ell} \right)^2.
\end{align}
Finally, using the mesh parameter given by~\eqref{eq:finest_mesh}, 
work per sample~\eqref{eq:work_model}, and for simplicity assuming
that~\eqref{eq:var_model_MLMC} also holds for $\ell=0$, this
computational work has the asymptotic behavior
\begin{align}
  \label{eq:complexity_MLMC}
  \work_\mathrm{MLMC} & \propto
    \begin{cases}
      \tol^{-2}, & \text{if $\Os>\gamma$},\\
      \tol^{-2}\left(\log{\tol^{-1}}\right)^2, & \text{if $\Os=\gamma$},\\
      \tol^{-2\left(1+\frac{\gamma-\Os}{2\Ow}\right)}, & \text{if $\Os<\gamma$},
    \end{cases}
\end{align}
as $\tol\to 0$, assuming $\Ow\geq\min{\left(\Os,\gamma\right)}/2$; 
see e.g. Theorem~3.1 in~\cite{giles08}, or
Corollary~2.1 and Corollary~2.2 in~\cite{Haji-Ali2016}.
Similar to the standard MC case, it is possible to optimize the
choice of $\splitting$ in~\eqref{eq:splitting} for MLMC. In
particular, if $\Os>\gamma$, an asymptotic analysis gives
$\splitting\to 1$, as $\tol\to 0$, indicating an aggressive refinement of
the numerical discretization to reduce the bias.
Again, the choice of $\splitting$ does not change the rate of the
complexity, but an optimal choice may reduce the work with a constant factor. 

In all three cases in~\eqref{eq:complexity_MLMC}, the complexity is
lower than the corresponding complexity,
$\tol^{-\left(2+\gamma/\Ow\right)}$, for standard MC simulation of 
the same problem~\eqref{eq:complexity_MC}. 
This leads to very significant computational saving in complex models,
and as a result some problems that are infeasable using the standard
MC method are computationally tractable using MLMC.

\paragraph{MLMC applied to $\E{\QoI_\ell}$ of
  Section~\ref{sec:comp_weak}--\ref{sec:comp_QoI}}

The assumption on the work per sample~\eqref{eq:work_model} holds for
$\gamma=d+1$, since the degrees of freedom in the 
uniform spatial discretization are proportional to $\meshparam^{-d}$, 
and the number of time steps is  proportional to $\meshparam^{-1}$,
where work per time step of the explicit time stepping scheme is
proportional to the degrees of freedom.
In the setting described in Section~\ref{sec:prob}--\ref{sec:QoI},
the weak convergence rate, $\Ow$, is identical to the rate of convergence in
the approximation of the deterministic problem, and the strong
convergence rate, $\Os/2$, equals the weak rate.
The explicit Newmark time stepping scheme and the numerical
approximation of $\QoI_E$ and $\QoI_W$ are both of order 2, so that
asymptotically as $\meshparam\to 0$ we expect $\Ow=2$ and $\Os=4$
assuming sufficiently regular exact solution. 
Based on these observations, summarized in
Table~\ref{tab:parameters_comp}, and the complexity
estimates~\eqref{eq:complexity_MC} and~\eqref{eq:complexity_MLMC},
we expect the asymptotic complexity to improve from $\tol^{-3.5}$ to $\tol^{-2}$, 
for $d=2$, and from $\tol^{-4}$ to $\tol^{-2}(\log{(\tol^{-1})})^2$, for $d=3$, as $\tol\to 0$ and standard MC is replaced by MLMC.

\begin{table}
  \centering
  \begin{tabular}[h]{|c|c|c|c|c|c|}
    \cline{2-6}
    \multicolumn{1}{c|}{}
    & 
    \multicolumn{3}{|c|}{ Model parameters }
    &
    \multicolumn{2}{|c|}{ Asymptotic complexity } \\
    \hline
    $d$\TBstrut & \hspace{2mm} $\gamma$ \hspace{2mm} & \hspace{2mm} $\Ow$
      \hspace{2mm} & \hspace{2mm} $\Os$\hspace{2mm} & $\work_{\mathrm{MC}}$ & 
      $\work_{\mathrm{MLMC}}$ \\
    \hline
    2 & 3 & \multirow{2}{*}{2} & \multirow{2}{*}{4} & $\tol^{-3.5}$ 
      & $\tol^{-2}$\Tstrut \\
    3 & 4 & & & $\tol^{-4}$ & $\tol^{-2}(\log{(\tol^{-1})})^2$\Tstrut \\
    \hline
  \end{tabular}
  \caption{Summary of the parameters in the work and convergence
    models~\eqref{eq:MC_models} and~\eqref{eq:var_model_MLMC}, for
    the numerical approximation of Problem~\ref{prob:weak1} or
    Problem~\ref{prob:weak2}, given in
    Section~\ref{sec:comp_weak}--\ref{sec:comp_QoI}, and the
    corresponding asymptotic complexity estimates, given
    in~\eqref{eq:complexity_MC} and~\eqref{eq:complexity_MLMC}.}
  \label{tab:parameters_comp}
\end{table}

\section{Numerical Tests}
\label{sec:num}

These numerical experiments make up an initial study of the validity
of MLMC techniques as a means of accelerating the approximation of
expected values of source inversion misfit functions, where we take the
expectation with respect to random parameters modeling
uncertainties in the Earth model. 
After this initial study where the source is approximated to a point
and only synthetic data are used, our ultimate goal is to integrate
MLMC into the full source inversion problem where the finite fault
solution is to be inferred by using real seismological data. 
While the final source inversion must be based on numerical simulations
on a three-dimensional Earth model, these initial tests were made on a
two-dimensional model described in the following. Furthermore, the
misfit functions were chosen with the aim of identifying the source
location considering the source moment tensor as fixed.

We first describe the problem setup, including the
source model, computational geometry, discretization, random Earth
material parameters, and the synthetic data replacing actual
measurements in the two-dimensional test. Finally, we describe the
execution and results of MLMC computations on the given problem setup. 

\subsection{Problem Setup}

\begin{figure}
  \begin{minipage}[c]{0.5\linewidth}
    \centering
    \includegraphics[width=\linewidth]{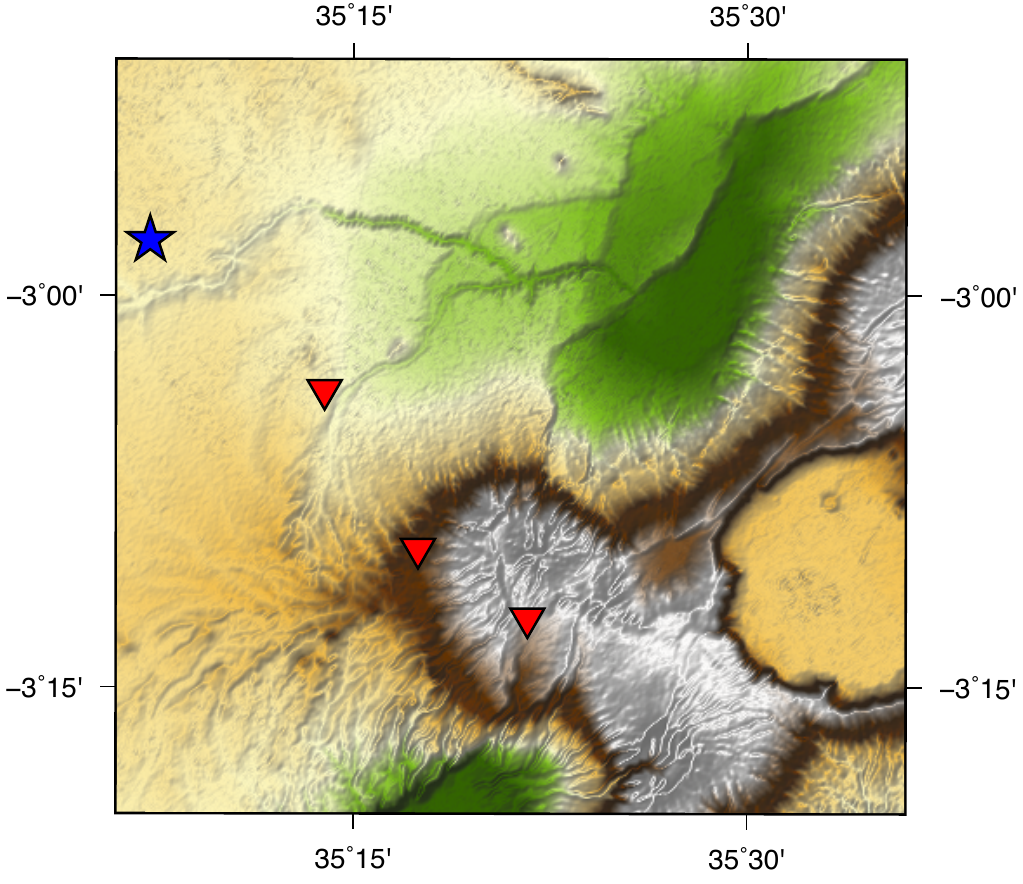}
  \end{minipage}
  \hfill
  \begin{minipage}[c]{0.4\linewidth}
    \caption{Source--receivers--geometry for the Tanzania case study,
      restricted to three receivers, marked by red triangles, which fall
      approximately along a straight line, also aligned with the
      estimated source location of a recorded seismic event (marked by
      a blue star).} 
    \label{fig:Tanzania_actual_domain}
  \end{minipage}
\end{figure}

For the numerical tests with $d=2$, we create a geometry consistent
with an actual network of receivers, belonging to a small seismic
network in the Ngorongoro Conservation Area on the East Rift, in
Tanzania. 
We do this by selecting three receivers that are approximately
aligned with the estimated epicenter of a seismic event that was
recorded.
Figure~\ref{fig:Tanzania_actual_domain}
illustrates the physical configuration.
The rough alignment of the source and the receiver locations in the
actual seismic network make this event a good opportunity to run the
tests in a two-dimensional domain. 
We describe the two-dimensional computational domains below, together with
the source and Earth parameters.

\subsubsection{Source model}

We consider a point source with a symmetric moment tensor, modeled as
a body force in the variational
equation~\eqref{eq:weak_eq_motion_attenuation} of
Problem~\ref{prob:weak2}, 
\begin{align*}
  \int_\domain \force(\xv,t)\cdot\testfun(\xv)\,d\xv & =
  -\mathbf{M}\doubledot \grad\testfun(\source)S(t),
\end{align*}
with the moment tensor 
\begin{align*}
  \mathbf{M} & = \left( 
    \begin{array}{cc}
      5.5895\e{13} & 7.9762\e{13} \\
      7.9762\e{13} & -2.5698\e{14}
    \end{array}
  \right),
\end{align*}
measured in $\mathrm{Nm}$, and a Gaussian source-time function with
corner frequency $f_{0}=2\,\mathrm{Hz}$, 
\begin{align*}
  S(t) & = 
  \dfrac{3f_{0}}{\sqrt{2\pi}} \exp\left(-\frac{9f_{0}^2(t-t_{c})^{2}}{2}\right).
\end{align*}
The time source function is centered at time $t_{c}=0\,\mathrm{s}$; the
solution time interval starts at $t_{0}=-0.6\,\mathrm{s}$ and ends at
$\mathcal{T}=25\,\mathrm{s}$, and the \QoItext~is based on $0\leq
t\leq\mathcal{T}$.

\subsubsection{Computational Domain}

The heterogeneous Earth is initially modeled with six homogeneous
layers of variable thickness, as stated in
Table~\ref{tab:thickness_layers}, and terminated by an infinite
half-space. The layers are separated by horizontal interfaces; 
topography is not included here. 
The source-receiver geometry is defined by the depth 
(vertical distance from the free surface) of the point source and 
the horizontal distances between the point source and the three 
receivers, as given in 
Table~\ref{tab:config_pars} and shown in Figure~\ref{fig:Domain}.

As described in Section~\ref{sec:comp_weak}, 
the half-plane domain is approximated by a finite domain, with
absorbing boundary conditions on the three artificial boundaries. In
the numerical approximation of Problem~\ref{prob:weak2} with seismic
attenuation, the PML boundary conditions are not perfectly absorbing,
but reflections are created at the boundary; see e.g. 
Section~3.4 in~\cite{SPECFEM2D}. The finite domain is
defined by three additional parameters $\sourcex$, $\Delta D_x$, and
$\Delta D_z$, which are chosen large enough so that no reflections
reach any of the three receivers during the simulation time interval,
$[-0.6,\mathcal{T}]\,\mathrm{s}$, given the maximal velocities in the
ranges of uncertainties. 

\begin{table}
  \centering
  \begin{minipage}{1.0\linewidth}
    \centering
    \includegraphics[width=0.8\linewidth]{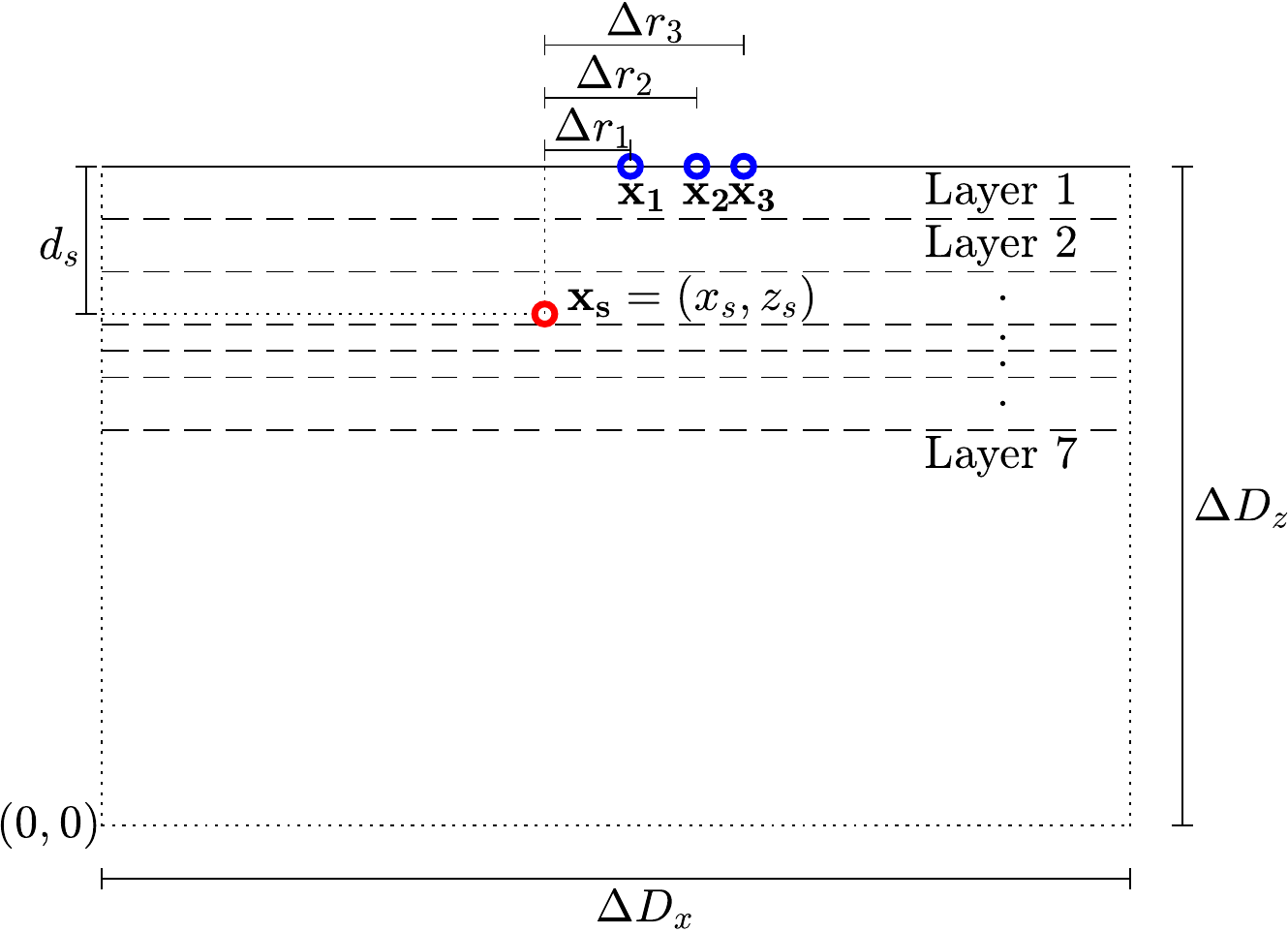}
    \captionof{figure}{Domain of the two-dimensional model. Four
      parameters define the configuration: the depth of the
      source, $d_s>0$, and the signed horizontal distances between the
      source and the three receivers, $\Delta r_1,\Delta
      r_2,\,\mathrm{and}\,\Delta r_3$. The horizontal coordinate
      direction is denoted $x$ and the vertical $z$. }
  \label{fig:Domain}
  \end{minipage}

  ~ \vspace{4mm} ~

  \begin{minipage}{0.28\linewidth}
    \begin{tabular}{ |c|c| }  
      \hline
      Layer & Thickness \\
      \hline
      1 & 10 $\mathrm{km}$ \\ 
      2 & 10 $\mathrm{km}$ \\ 
      3 & 10 $\mathrm{km}$ \\ 
      4 & 5 $\mathrm{km}$  \\ 
      5 & 5 $\mathrm{km}$  \\ 
      6 & 10 $\mathrm{km}$ \\  
      7 & -- \\  
      \hline
    \end{tabular}    
    \captionof{table}{Thickness of the layers in the half-plane two dimensional
      Earth model.}  
    \label{tab:thickness_layers}
  \end{minipage}
  \hfill
  \begin{minipage}{0.57\linewidth}
  \begin{tabular}{ |c|r|r| }
    \cline{2-3}
    \multicolumn{1}{c|}{} & Synthetic Data & \multicolumn{1}{|c|}{MLMC} \\ 
    \hline
    $d_s$        & $28.000\,\mathrm{km}$ & $28.000\,\mathrm{km}$ \\
    $\Delta r_1$ & $16.242\,\mathrm{km}$ & $11.242\,\mathrm{km}$ \\ 
    $\Delta r_2$ & $28.849\,\mathrm{km}$ & $23.849\,\mathrm{km}$ \\ 
    $\Delta r_3$ & $37.724\,\mathrm{km}$ & $32.724\,\mathrm{km}$ \\ 
    \hline
    $x_s$ &  $84.000\,\mathrm{km}$ &  $86.500\,\mathrm{km}$ \\
    $\Delta D_x$ & $195.000\,\mathrm{km}$ & $195.000\,\mathrm{km}$ \\
    $\Delta D_z$ & $125.000\,\mathrm{km}$ & $125.000\,\mathrm{km}$ \\
    \hline
  \end{tabular}
  \captionof{table}{Parameters defining the configuration in the half-plane
    geometry, $d_s$, $\Delta r_1$, $\Delta r_2$, and $\Delta r_3$, and
    additional parameters defining the truncated numerical domain,
    $x_s$, $\Delta D_x$, and $\Delta D_z$. See
    Figure~\ref{fig:Domain}.}
  \label{tab:config_pars}
  \end{minipage}
\end{table}

\paragraph{Discretization of the computational domain}

For the numerical computations using {\tt SPECFEM2D} version
$7.0$~\cite{SPECFEM2D} in double precision,
the computational domain is discretized uniformly into squares of side
$\dx$, with the coarsest mesh using $\dx[0]=2500\,\mathrm{m}$ so that
the interfaces always coincide with element boundaries.
These computations are based on spectral elements in two dimensions, 
using $n_c=9$ control points 
to define the isomorphism between the computational
element and the reference element, a basis of Lagrange polynomials of
degree $N_l=4$, and $5\times5$ Gauss-Lobatto-Legendre quadrature points.
The second order Newmark explicit time stepping scheme was used with
step size $\dt[0]=6.25\e{-3}\,\mathrm{s}$ on the coarsest
discretization, which was refined at the same rate as $\dx$ to keep
the approximate CFL condition satisfied.
A PML consisting of three elements was used on the artificial
boundaries. 

\subsubsection{Earth material properties}
\label{sec:Earth_material}

The viscoelastic property of the Earth material, as described in
Section~\ref{sec:comp_attenuation}, is approximated by a generalized
Zener model, implemented in {\tt SPECFEM2D}, with $B=3$ SLS and the
quality factor, $Q$, which is constant in each layer
(Table~\ref{tab:Earth_model_2D}). The quality factor is kept constant 
throughout the simulations.

As described in Section~\ref{sec:prob}, the triplet $(\density,
\comprv, \shearv)$, denoting the density, compression wave speed, and
shear wave speed, respectively, defines the Earth's material
properties with varying spatial position. In the particular
seven-layer domain  introduced above, a one-dimensional, piecewise
constant velocity model is used.
These three fields are then completely described by three
seven-dimensional 
random variables, $\bm{\comprv}$, $\bm{\shearv}$, and $\bm{\density}$.
Here, we detail the probability distributions we assign to these
parameters. 
To prepare the inversion of real data, $\comprv$, $\shearv$, and
$\density$ are adapted from the results of~\cite{Albaric2010}
and~\cite{Roecker2017} obtained from previous seismological
experiments in adjacent areas. 
Among the unperturbed values, denoted with a bar over the
symbols, listed in Table~\ref{tab:Earth_model_2D}, $\bar{\bm{\shearv}}$ and
$\bar{\bm{\density}}$ are treated as primary parameters, while
$\bm{\bar{\comprv}}$ is scaled from $\bar{\bm{\shearv}}$. The relation
\begin{align}
  \label{eq:DetVp}
  \nu & = \dfrac{\bar{\comprv}_{i}}{\bar{\shearv}_{i}} =1.7, 
  && i=1,2,\dots,7, 
\end{align}
is chosen because it is a common use for crustal structure and it is
in agreement with previous seismological studies in the
area~\cite{Roecker2017}.

\begin{table}[h]
  \centering
  \begin{tabular}{ |c|c|c|c|c| }  
    \hline
    layer, $i$  & $\bar{\density}_{i}$ & $\bar{\shearv}_{i}$ &
      $\bar{\comprv}_{i}$ & $Q$\Tstrut\\
    \hline
    1 & $2500$ & $3529.0$ & $6034.6$ & $300$ \\ 
    2 & $2500$ & $3705.0$ & $6335.6$ & $300$ \\ 
    3 & $2500$ & $3882.0$ & $6638.2$ & $800$ \\ 
    4 & $2500$ & $3911.0$ & $6687.8$ & $800$ \\ 
    5 & $2900$ & $4422.7$ & $7562.8$ & $800$ \\ 
    6 & $2900$ & $4506.4$ & $7705.9$ & $600$ \\
    7 & $2900$ & $4533.6$ & $7752.5$ & $600$ \\
    \hline
  \end{tabular}
  \caption{Unperturbed values of the material parameters. Here 
    $\bar\density$, $\bar{\shearv}$, and $\bar{\comprv}$ are given in
    the units
    $\mathrm{kg}/\mathrm{m}^3$, $\mathrm{m}/\mathrm{s}$, and
    $\mathrm{m}/\mathrm{s}$, respectively. The quality factor, $Q$,
    used in the seismic attenuation model is dimensionless and kept
    unperturbed in all simulations.}
  \label{tab:Earth_model_2D}
\end{table}

We model the uncertain shear wave speed, $\bm{\shearv}$, as a uniformly
distributed random variable
\begin{align*}
  \pmb{\shearv} & \sim
  \bm{\unifrnd}\left(
    \prod_{i=1}^{7}[\shearv_{i}^{lb},\shearv_{i}^{ub}]
  \right),
\end{align*}
with independent components, and where the range is a plus-minus 10\%
interval around the unperturbed value, i.e.
\begin{align*}
  \shearv_{i}^{lb} & =(1-q)\cdot\bar{\shearv}_{i},  
  &
  \shearv_{i}^{ub} & =(1+q)\cdot\bar{\shearv}_{i},
  && i=1,2,\dots,7, 
\end{align*}
where $q=0.1$. 
In keeping with~\eqref{eq:DetVp}, but assuming some variability in the
ratio, the compressional wave speed, $\bm{\comprv}$, is modeled by a 
random variable which, conditioned on $\bm{\shearv}$, is uniformly
distributed with independent components, 
\begin{align*}
  \bm{\comprv} &\sim 
  \bm{\unifrnd}\left(
    \prod_{i=1}^{7}[\nu^{lb}\cdot \shearv_{i},\nu^{ub}\cdot \shearv_{i}]
  \right),
\end{align*}
where $\nu^{lb}=1.64$ and $\nu^{ub}=1.78$, corresponding to a
range of variability of about $\pm4\%$. 
Finally, the density, $\bm{\density}$, is again uniformly distributed with
independent components
\begin{align*}
  \bm{\density} & \sim \bm{\unifrnd}\left(
    \prod_{i=1}^{7}[\density_{i}^{lb},\density_{i}^{ub}]
  \right),
\end{align*}
where
\begin{align*}
  \density_{i}^{lb} & =(1-r)\bar{\density}_{i}, &
  \density_{i}^{ub} & =(1+r)\bar{\density}_{i}, &&
  i=1,2,\dots,7, 
\end{align*}
with $r=0.1$.

\begin{figure}
  \centering
  \includegraphics[width=0.49\linewidth]{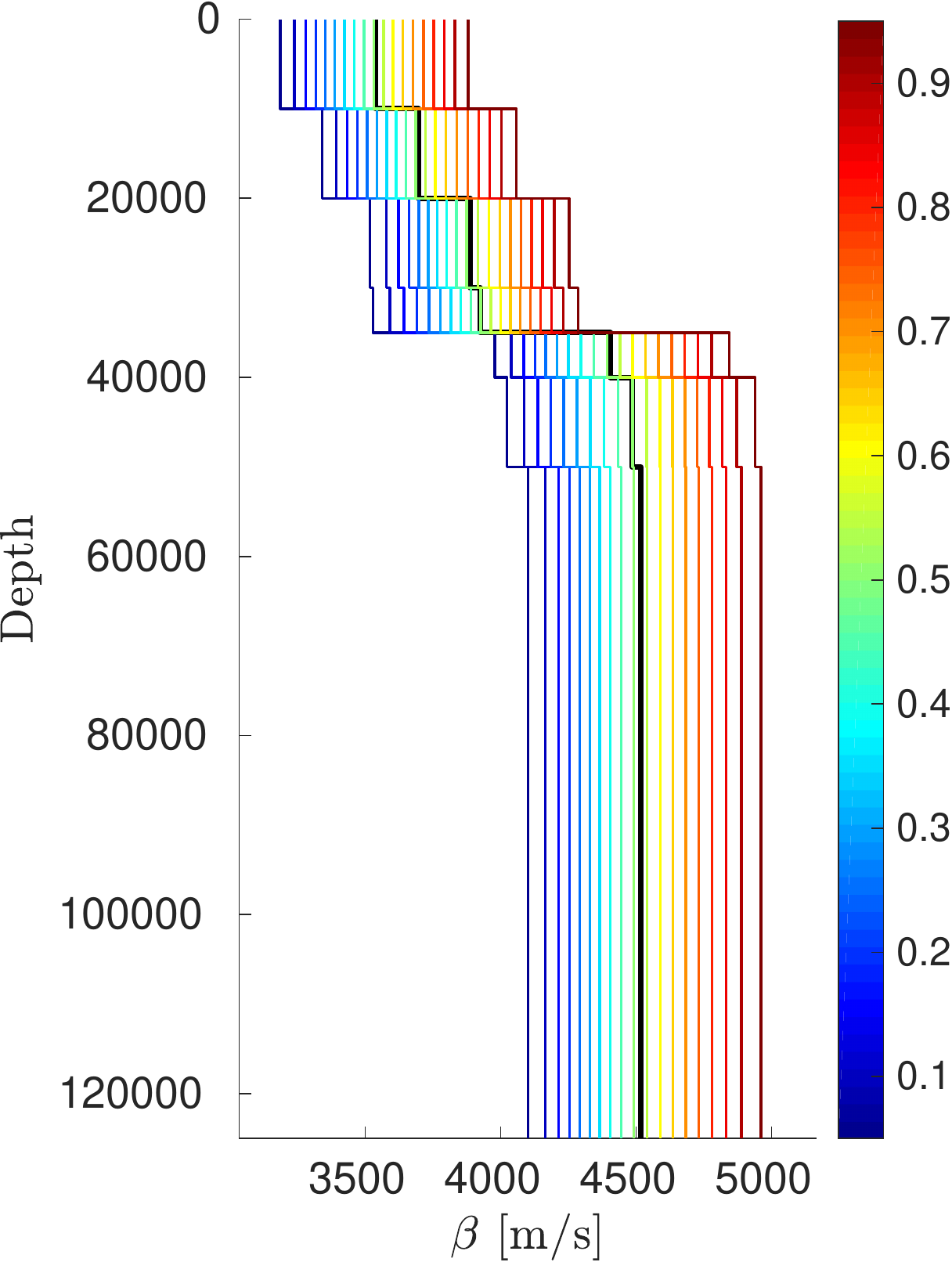}
  \includegraphics[width=0.49\linewidth]{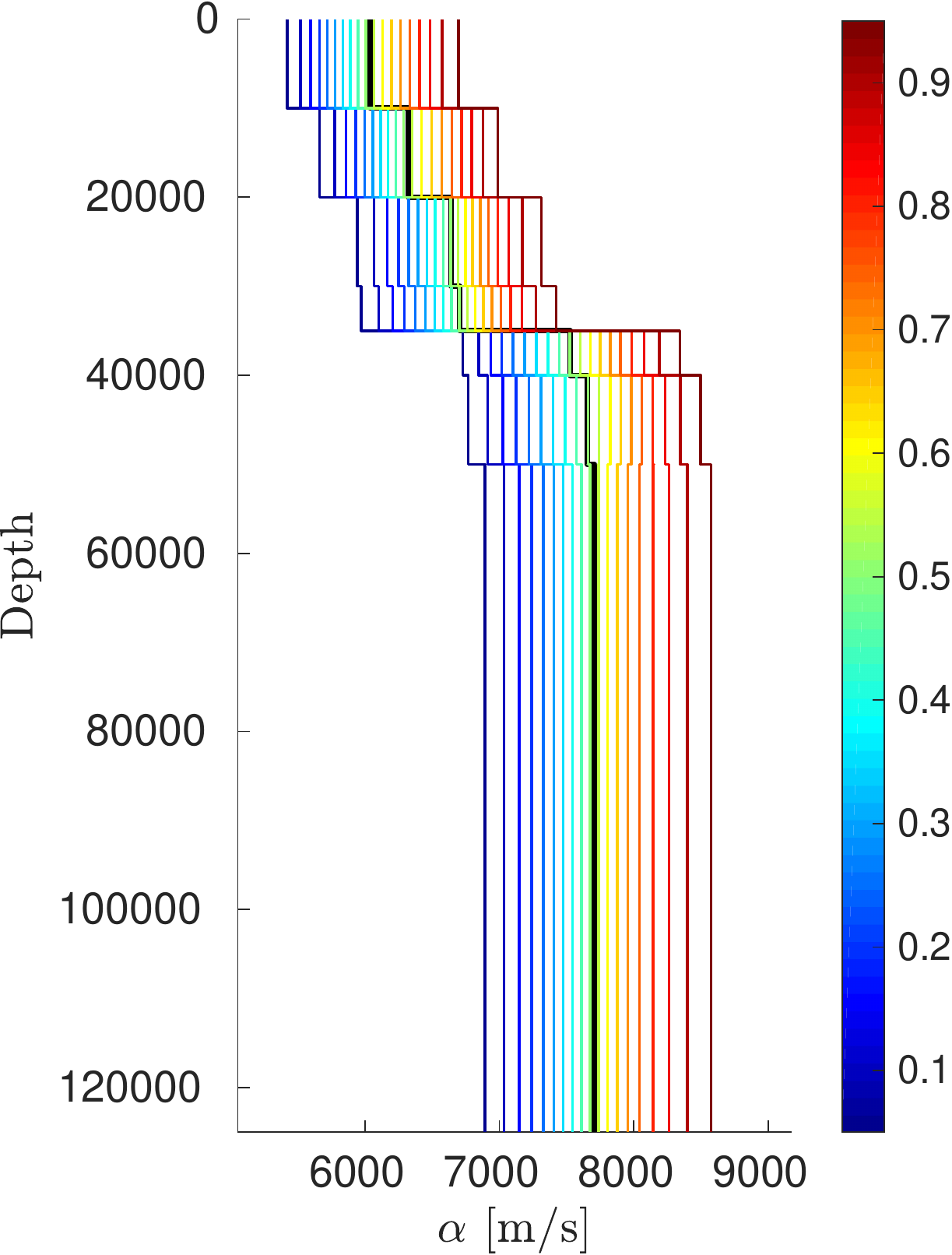}\\
  \includegraphics[width=0.49\linewidth]{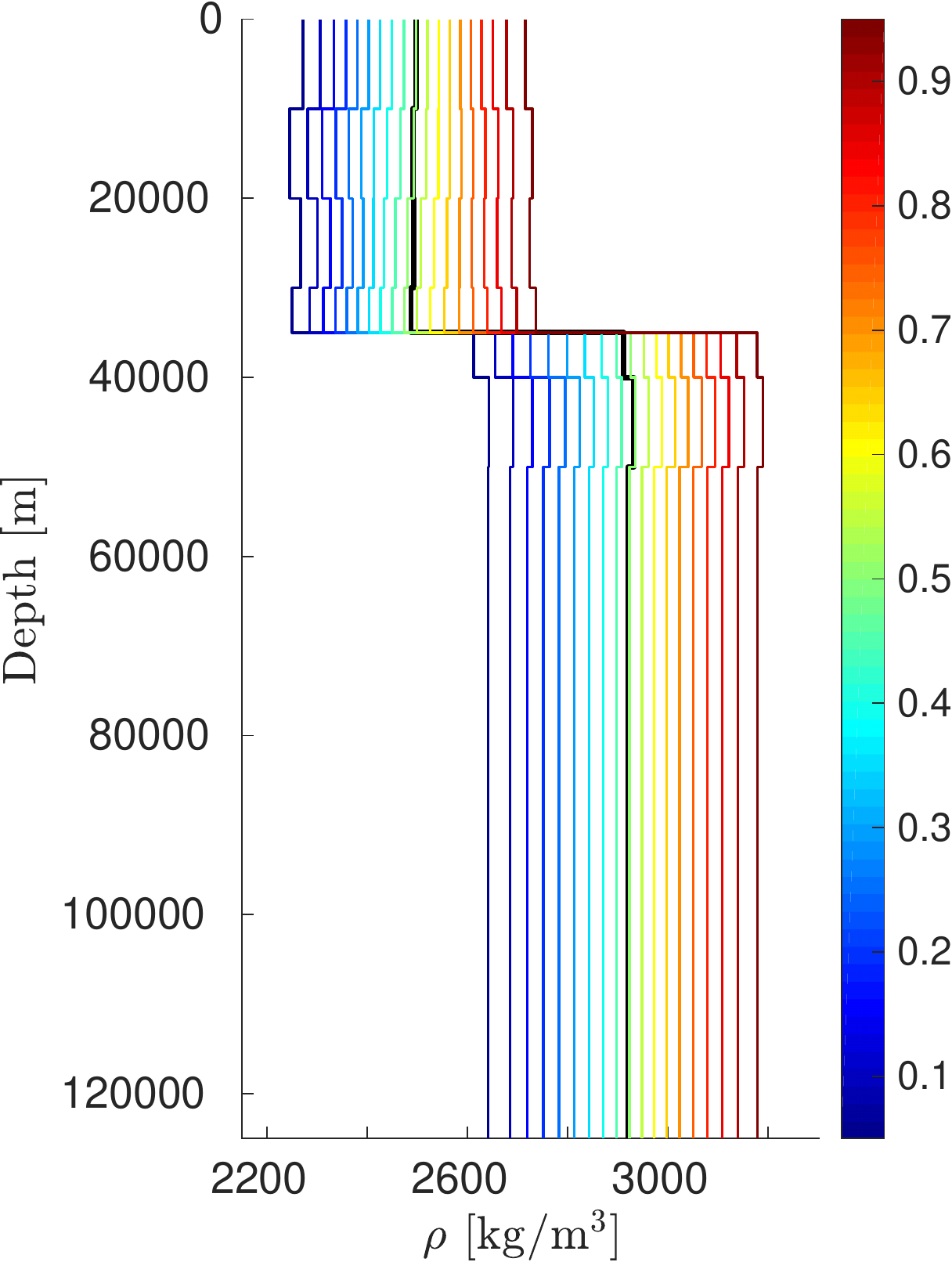}
  \includegraphics[width=0.49\linewidth]{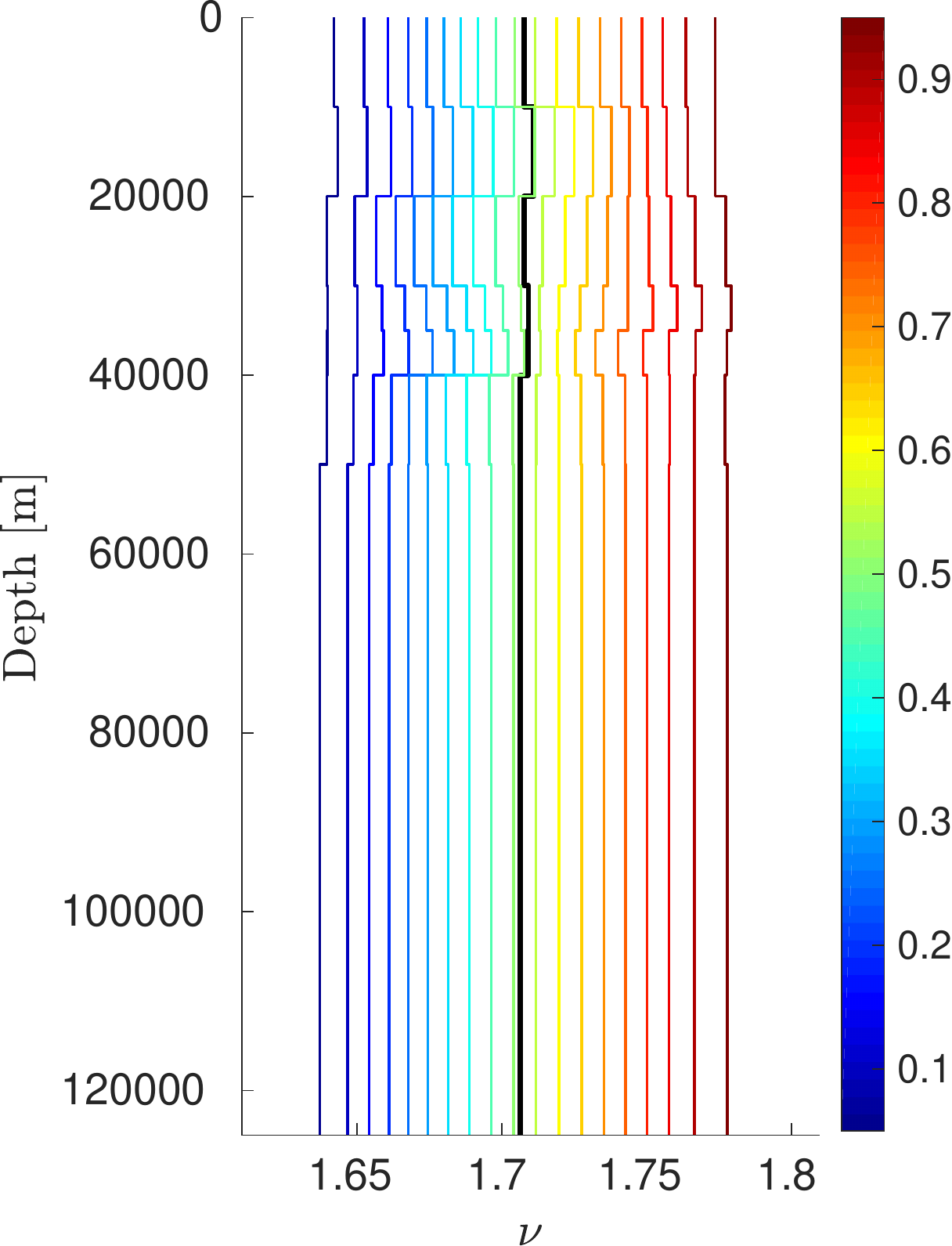}
  \caption{Sample average, marked with the thick black graph, and
    equidistributed quantiles, color-coded thinner graphs, of the Earth
    material parameters based on the 320 samples from the verification
    run.\newline 
    (Top Left) Shear wave speed, $\shearv$, (Top Right) Compressional
    wave speed, $\comprv$, (Bottom Left) Density, $\density$, and
    (Bottom Right) $\nu=\comprv/\shearv$.}
  \label{fig:VerRun_distributions}
\end{figure}

Figure~\ref{fig:VerRun_distributions} shows 
the sample mean and contour lines of the sample CDF 
for these random parameters, based on samples
used in the verification run of the problem setup.

\subsubsection{Synthetic Data}
\label{sec:num_res_synth}

Instead of actual, measured data from seismic activity, the misfit
functions for the \QoItext~in the two-dimensional computations use
synthetic data obtained from the same underlying code using a finer
discretization, $\dx[]=78.125\,\mathrm{m}$ and $\dt[]=1.953125\e{-4}\,\mathrm{s}$,
than any of the samples in the MLMC run.
The source location relative to the receivers, listed in
Table~\ref{tab:config_pars}, agrees with the independently estimated epicenter in
Figure~\ref{fig:Tanzania_actual_domain} and the fixed Earth material 
parameters, listed in Table~\ref{tab:Earth_model_2D_synthetic},
correspond to one outcome of the sampling procedure in
Section~\ref{sec:Earth_material}.
The computed displacements are illustrated in the left column of
Figure~\ref{fig:Syn_Data}. 

The resulting time series for the displacement in the three
receivers are then restricted to a much coarser time discretization, 
corresponding to a frequency of
measurements of $160\,\mathrm{Hz}$, that are in the realistic range of
frequencies for measured seismograms, and i.i.d. noise
$\WhiteNoise\sim\mathcal{N}(\mathbf{0},\sigma^{2}\Imatrix)$ with
$\sigma=2.5\e{-3}\approx 1\%$ of $\max{\displ}$, is added as
in~\eqref{eq:data}; see the right column of Figure~\ref{fig:Syn_Data}.

\begin{table}[h!]
  \centering
  \begin{minipage}{1.0\linewidth}
    \includegraphics[width=0.55\linewidth]{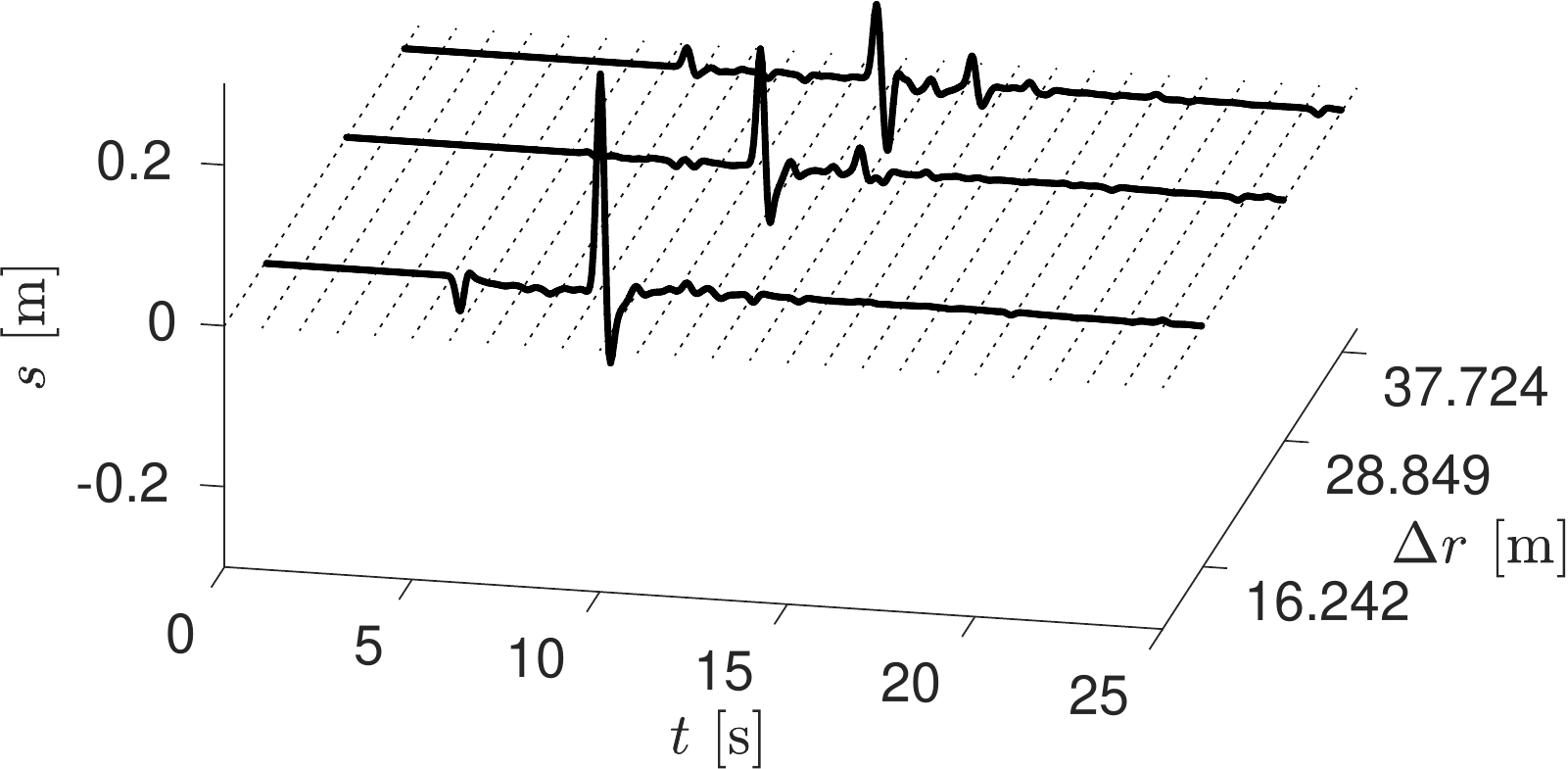}
    \hfill
    \includegraphics[width=0.4\linewidth]{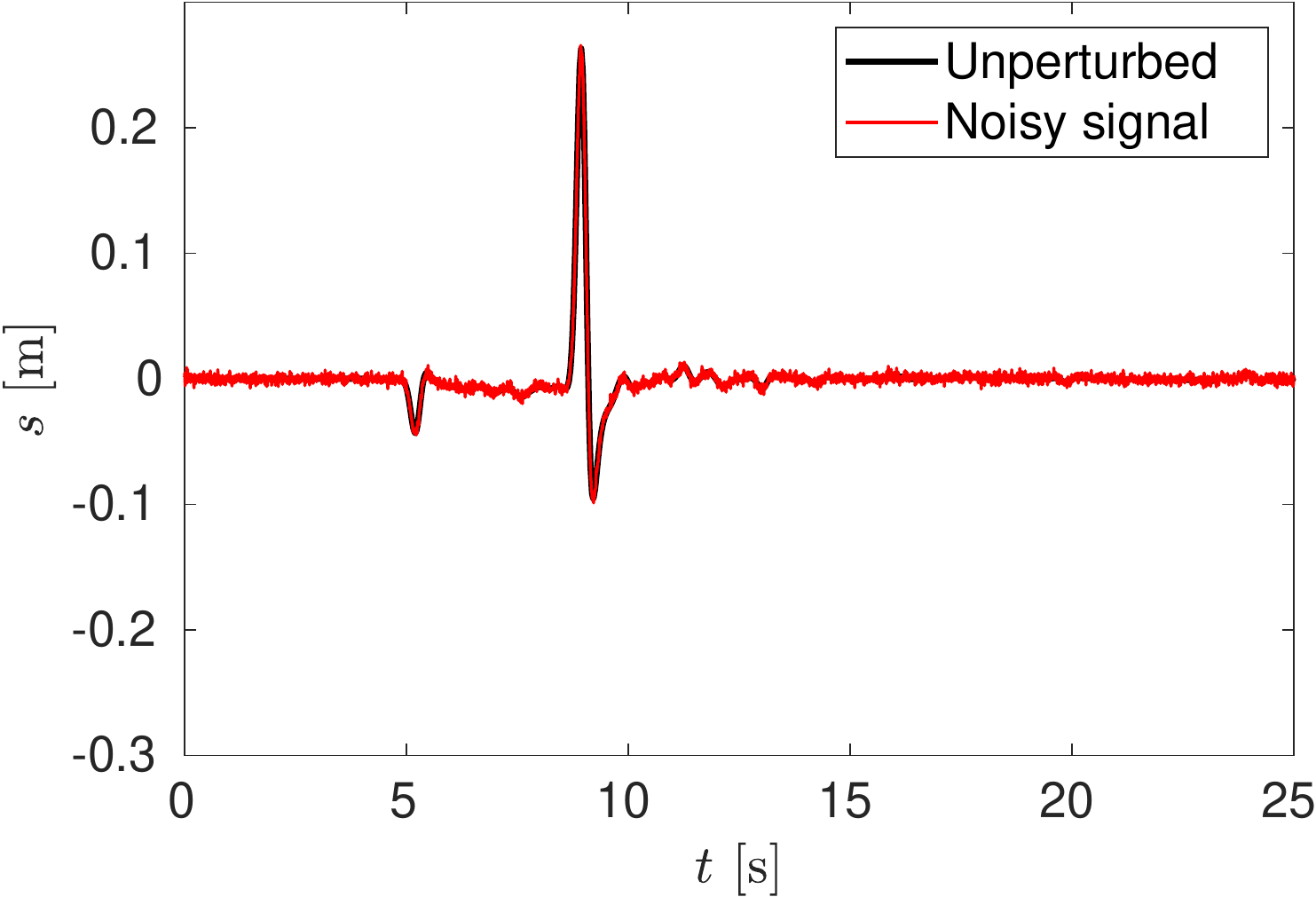}\\
    \includegraphics[width=0.55\linewidth]{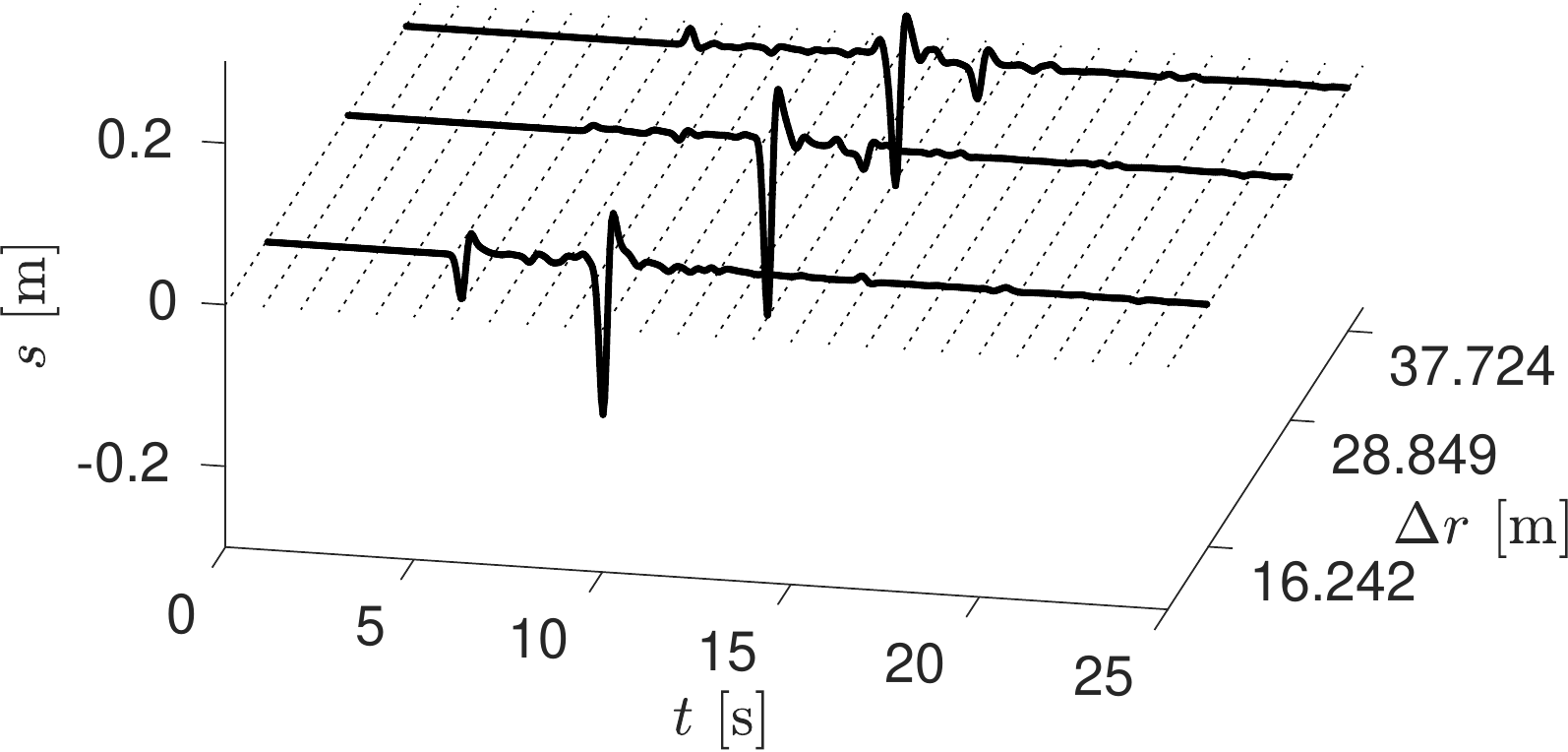}
    \hfill
    \includegraphics[width=0.4\linewidth]{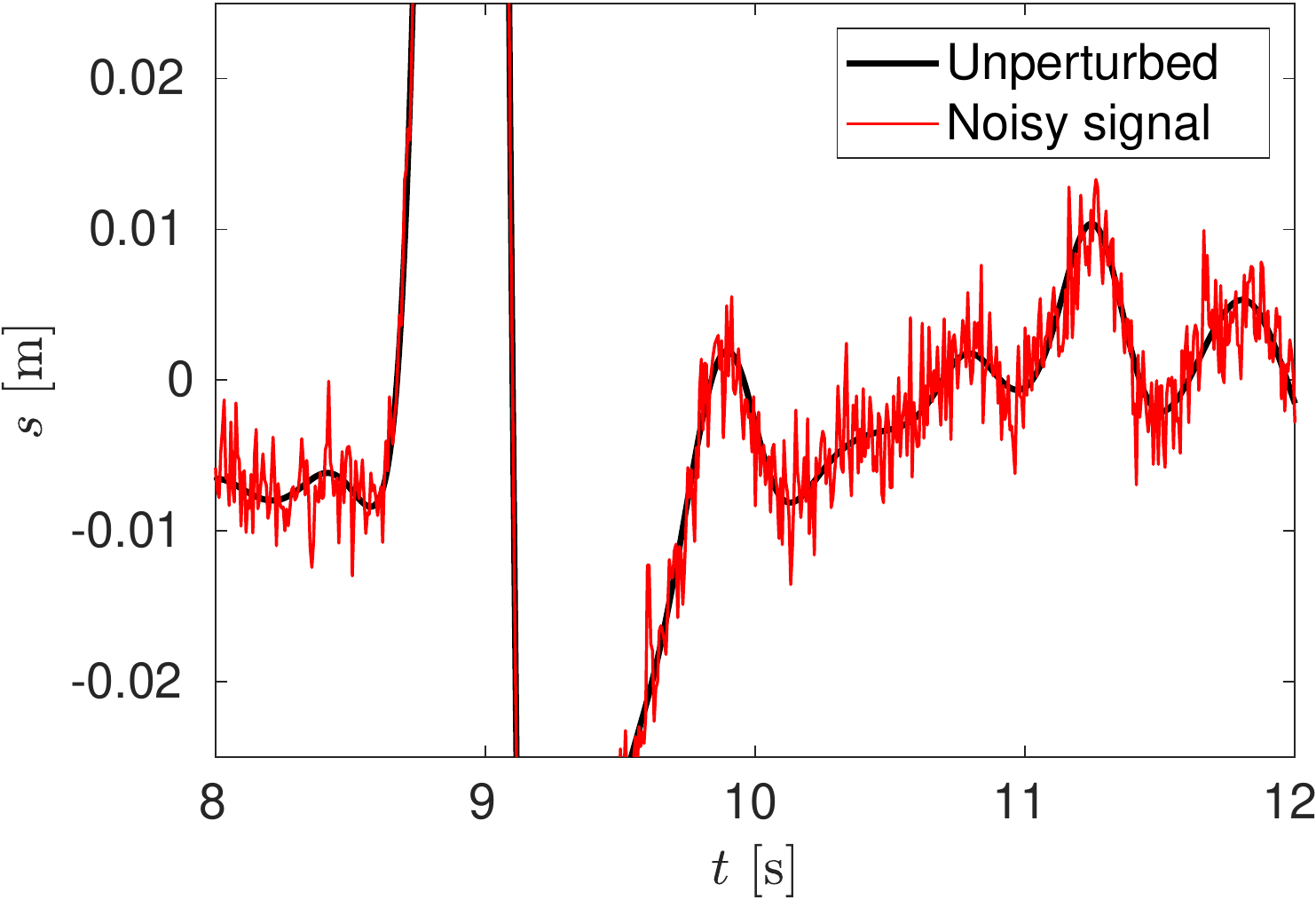}
    \captionof{figure}{Synthetic data
      $\data(\xv_{n},t)=\displ(\xv_{n},t;\source[\ast],\ParEar^\ast)$,
      where the displacement
      $\displ(\xv_{n},t;\source[\ast],\ParEar^\ast)$ is obtained from a
      computation with the source location, $\source[\ast]$, and the
      geometry parameters of Figure~\ref{fig:Domain} given in
      Table~\ref{tab:config_pars} and  
      with a fixed outcome of the random Earth parameters,
      $\ParEar^\ast$, within their ranges of uncertainty, as given in
      Table~\ref{tab:Earth_model_2D_synthetic}. To the left are the
      $x$-component (top) and $z$-component (bottom) of the synthetic
      data without added noise. To the right, the $x$-component in
      receiver location 1 is shown both with and without added noise,
      with the bottom being a detailed view of the top.}
    \label{fig:Syn_Data}
  \end{minipage}

  \vspace{4mm}

  \begin{minipage}{1.0\linewidth}
    \begin{minipage}{0.49\linewidth}
      \centering
      \begin{tabular}{ |c|c|c|c|c| }  
        \hline
        layer, $i$  & $\density_{i}$ & $\shearv_{i}$ & $\comprv_{i}$\\
        \hline
        1 & $2439.9$ & $3498.0$ & $5737.7$ \\
        2 & $2715.8$ & $3654.9$ & $6346.0$ \\
        3 & $2747.1$ & $3690.9$ & $6568.5$ \\
        4 & $2562.0$ & $4045.0$ & $7038.2$ \\
        5 & $2862.8$ & $4491.8$ & $7647.2$ \\
        6 & $2862.0$ & $4691.5$ & $7753.4$ \\
        7 & $2809.7$ & $4969.8$ & $8790.9$ \\
        \hline
      \end{tabular}
      \captionof{table}{Values of the material parameters in the simulation which
        generated the synthetic data. Here 
        $\density$, $\shearv$, and $\comprv$ are given in
        the units
        $\mathrm{kg}/\mathrm{m}^3$, $\mathrm{m}/\mathrm{s}$, and
        $\mathrm{m}/\mathrm{s}$, respectively, and rounded to
        five digits.}
      \label{tab:Earth_model_2D_synthetic}
    \end{minipage}
    \hfill
    \begin{minipage}{0.49\linewidth}
      \centering
      \begin{tabular}{ |c|c|c|c|c| }  
        \hline
        layer, $i$  & $\density_{i}$ & $\shearv_{i}$ & $\comprv_{i}$\\
        \hline
        1 & $2333.8$ & $3789.8$ & $6503.4$ \\
        2 & $2539.0$ & $3562.3$ & $6117.5$ \\
        3 & $2416.8$ & $3726.8$ & $6157.8$ \\
        4 & $2521.5$ & $3724.4$ & $6176.2$ \\
        5 & $2793.1$ & $4062.9$ & $7228.6$ \\
        6 & $2822.0$ & $4422.7$ & $7602.2$ \\
        7 & $2811.9$ & $4612.2$ & $7733.8$ \\
        \hline
      \end{tabular}
      \captionof{table}{Values of the material parameters in the
        simulations used to assess the effect of including attenuation
        in the model. Units and rounding of 
        $\density$, $\shearv$, and $\comprv$ as in
        Table~\ref{tab:Earth_model_2D_synthetic}.\\~}
      \label{tab:Earth_model_elastic}
    \end{minipage}
  \end{minipage}
\end{table}

\subsubsection{The impact of attenuation on the \QoItexts}
\label{sec:comp_elastic}

To assess the effect of attenuation in the problem described above, we
compute the \QoItexts~obtained with and without attenuation in the
model, for one outcome of the Earth material parameters, given in
Table~\ref{tab:Earth_model_elastic}, and with the geometry of the MLMC
runs in Table~\ref{tab:config_pars}, using discretization levels 1
and 2 in Table~\ref{tab:mesh_hierarchy_2D}. Including attenuation 
changed both \QoItexts~several percent; see
Table~\ref{tab:comp_elastic}. 

The significantly reduced computing time of MLMC, compared to
standard MC with corresponding accuracy, allows for simulation with
an error tolerance small enough to evaluate the usefulness of
including attenuation in the model in the presence of uncertainty in
the Earth material parameters. 

An alternative variant of the MLMC approach in this paper, would be to
use \QoItexts~sampled using the elastic model as control variates for
\QoItexts~sampled using the model with attenuation. Since the work
associated with the elastic model is smaller, we expect an MLMC method
where coarse grid samples are based on the elastic model to further
reduce the computational cost of achieving a desired accuracy in the
expected value of the \QoItext. 

\begin{table}
  \centering
  \begin{tabular}{|c|c|c|c|c|}
    \hline
    & Level, $\ell$ & Elast. & Atten. & \textbf{Change} \\ 
    \hline
    \multirow{2}{*}{$\QoI_E$} 
    & 1 & 4.20\e{-3} & 3.82\e{-3} & \textbf{-9.1\%}\\
    \cline{2-5} 
    & 2 & 4.24\e{-3} & 3.85\e{-3} & \textbf{-9.2\%}\\ 
    \hline
    \multirow{2}{*}{$\QoI_W$} 
    & 1 & 8.94\e{-2} & 1.11\e{-1} & \textbf{23.6\%}\\
    \cline{2-5} 
    & 2 & 1.33\e{-1} & 1.37\e{-1} & \textbf{3.0\%}\\ 
    \hline
  \end{tabular}
  \caption{Effect of including attenuation in the Earth material model.}
  \label{tab:comp_elastic}
\end{table}

\subsection{MLMC Tests}
\label{sec:num_res_mlmc}

In this section, we describe how we apply the MLMC algorithm to the 
test problem introduced above and present results showing a
significant decrease in cost in order to achieve a given accuracy,
compared to standard MC estimates.

\subsubsection{Verification and parameter estimation}
\label{sec:ver_run}

In a verification step, we compute a smaller number of samples on
four discretization levels corresponding to a repeated halving of
$\dx[]$ and $\dt[]$, as specified in
Table~\ref{tab:mesh_hierarchy_2D}. Statistics of the underlying Earth
material samples on the coarsest discretization are illustrated in
Figure~\ref{fig:VerRun_distributions}.

This verification step is necessary to verify that our problem
configuration works with the underlying code as expected. At the 
same time, the assumptions~\eqref{eq:MC_models}
and~\eqref{eq:var_model_MLMC} are tested, by experimentally observing
the computation time per sample on the different levels, as well as sample
averages~\eqref{eq:MC_est} and sample variances, 
\begin{align}
  \varMC(\QoI(\ParEar)) & = 
  \frac{1}{\nrs-1}\sum_{n=1}^\nrs
    \left(
      \QoI(\ParEar_n)-\frac{1}{\nrs}\sum_{m=1}^\nrs\QoI(\ParEar_m)
    \right)^2,
  \label{eq:varMC}
\end{align}
 of either \QoItext~in
Section~\ref{sec:comp_QoI}, $\QoI_{\ast,\ell}$, $\ell=0,1,2,3$, and the
corresponding two-grid correction terms $\Delta\QoI_{\ast,\ell}$,
$\ell=1,2,3$.

\paragraph{Work estimates}

The cost per sample is taken to be the cost of generating one sample
of the displacement time series
$\{\displ(\rec[r,n],t_j;\cdot,\cdot)\}_{j=0,n=1}^{J,\quad\Nrec}$
using \specfem. It is measured as the reported elapsed time from the
job scheduler on the supercomputer, multiplied by the number of cores
used, as given in Table~\ref{tab:mesh_hierarchy_2D}. The
post-processing of the time series to approximate 
the~\QoItext~and compute the MLMC estimators, as outlined in 
Section~\ref{sec:comp_QoI}--\ref{sec:MLMC} without any additional filters 
etc., is performed on laptops and
workstations at a negligible cost, compared to the reported time.

The time per sample on a given level varied very little, and its
average over the samples in the verification run, shown in
Figure~\ref{fig:VerRun_work}, verifies the expectation from
Section~\ref{sec:comp_weak} that $\work_\ell\propto\meshparam^{-3}$,
corresponding to $\gamma=3$ in~\eqref{eq:work_model}.

\begin{table}[h]
  \centering
  \begin{minipage}{1.0\linewidth}
    \centering
    \begin{tabular}{ |c|r|l|c||c| }  
      \hline
      Level, $\ell$ & \multicolumn{1}{|c|}{$\dx$} & 
        \multicolumn{1}{|c|}{$\dt$} & \# Cores & $\nrs_\ell$ (Ver.Run) \\
      \hline
      0 & $2500.0\,\mathrm{m}$ & $6.2500\e{-3}\,\mathrm{s}$  & 4  & 160\\ 
      1 & $1250.0\,\mathrm{m}$ & $3.1250\e{-3}\,\mathrm{s}$  & 4  & 160\\ 
      2 & $625.0\,\mathrm{m}$  & $1.5625\e{-3}\,\mathrm{s}$  & 16 &  40 \\ 
      3 & $312.5\,\mathrm{m}$  & $7.8125\e{-4}\,\mathrm{s}$  & 64 &  10 \\ 
      \hline
    \end{tabular}
    \captionof{table}{Discretization and computer parameters:
      $\dx$, side of square elements in the uniform spatial
      discretizations, 
      $\dt$, uniform time step size,
      \#~Cores, the number of cores per sample used by \specfem, and
      $\nrs_\ell$, number of samples per level in the MLMC
      discretization hierarchy for the verification run.}
    \label{tab:mesh_hierarchy_2D} 
  \end{minipage}

  ~ \vspace{4mm} ~

  \begin{minipage}{0.5\linewidth}
    \centering
    \includegraphics[width=\linewidth]{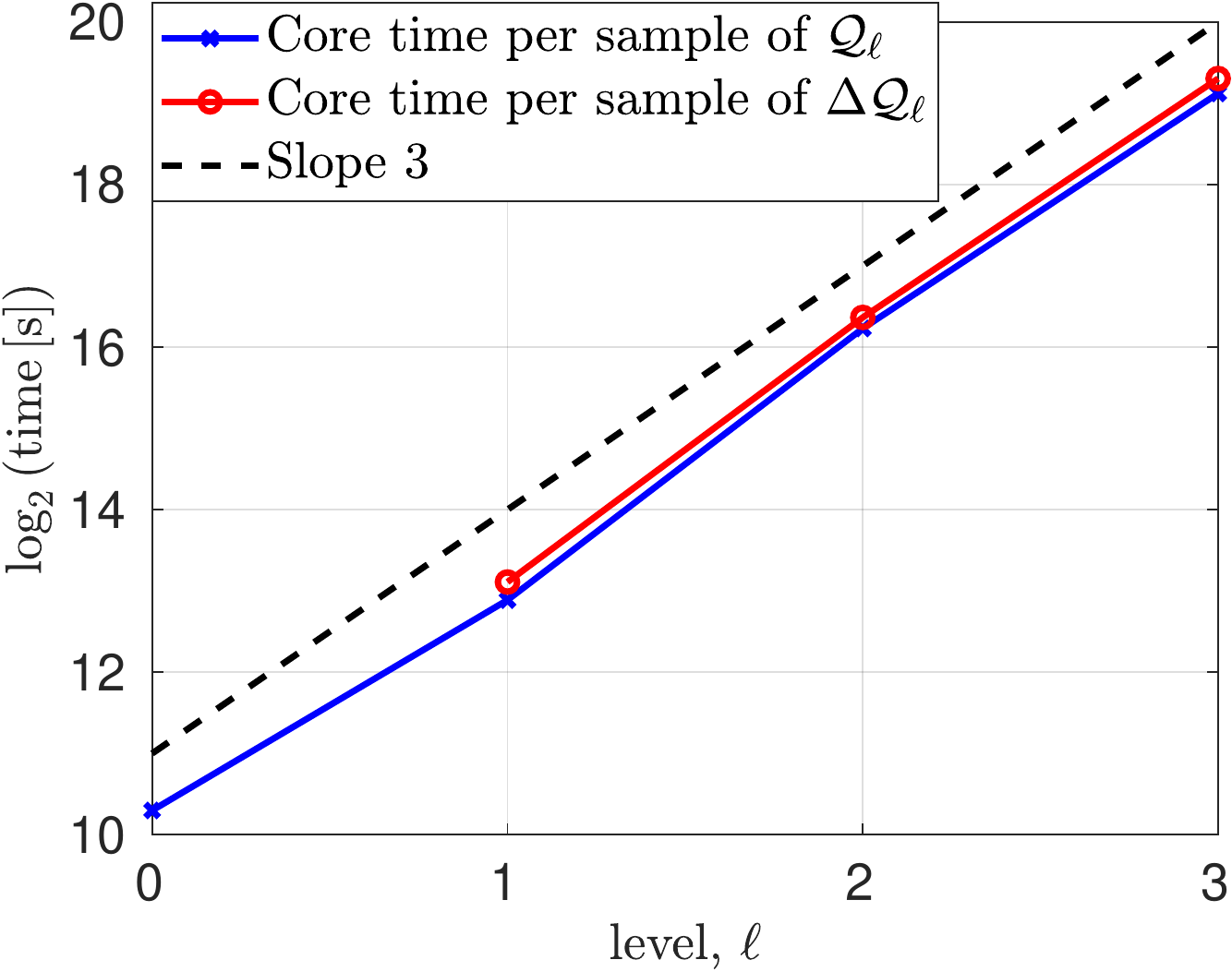}
  \end{minipage}
  \hfill
  \begin{minipage}[c]{0.4\linewidth}
    \captionof{figure}{Work per sample from the verification run measured as the
      core time, i.e. the elapsed time reported from the
      supercomputer's job scheduler multiplied by the number of cores
      per job, given in Table~\ref{tab:mesh_hierarchy_2D}.}
    \label{fig:VerRun_work}
  \end{minipage}
\end{table}

\paragraph{\QoItext~based on the $L^2$-misfit}

For $\QoI_{E,\ell}$ defined in~\eqref{eq:QoI_E}, both
$\E{\Delta\QoI_{E,\ell}}$ and $\var{\Delta\QoI_{E,\ell}}$ appear to
decrease faster in this range of discretizations than the
asymptotically optimal rates, which are $\Ow=2$ and $\Os=4$ given the
underlying numerical approximation methods; see
Figure~\ref{fig:VerRun_conv_L2}. The increased convergence rates 
indicate that we are in a pre-asymptotic regime where, through the 
stability constraint, the time step is taken so small that the time
discretization error does not yet dominate the error in $\QoI_E$, the
way it eventually will as $\tol\to 0$. In the context of the intended
applications for inverse problems, it can be justified to solve the
forward problem to a higher accuracy than one would demand in a
free-standing solution to the forward problem, but it may still be
unrealistic to use smaller relative tolerances than those used in
this study. This means we can expect to remain in the
pre-asymptotic regime, and simple extrapolation-based estimates of the
bias will be less reliable. 

We note from the sample variances, $\varMC(\QoI_{E,\ell})$, that the
standard deviation of $\QoI_E$ is of the order $5\e{-4}$, whereas
$\E{\QoI_E}$ itself is of the order $4\e{-3}$ (see the reference value
in Table~\ref{tab:ref_hierarchy}).

Given that $\varMC\left(\Delta\QoI_{E,1}\right)$ is significantly
smaller than $\varMC\left(\QoI_{E,0}\right)$, and that the cost per
sample of $\QoI_{E,0}$ is significantly smaller than the corresponding
cost of $\QoI_{E,1}$, it is intuitively clear that the optimal MLMC
approximation should include samples starting at the discretization
level labeled $\ell=0$ in Table~\ref{tab:mesh_hierarchy_2D}, and that 
the finest level, $L$, will depend on the tolerance, $\tol$. 

\paragraph{\QoItext~based on $W_2^2$-distances}

For $\QoI_{W,\ell}$ defined in~\eqref{eq:QoI_W}, we again seem to be
in the pre-asymptotic regime;
see Figure~\ref{fig:VerRun_conv_W}. 
Here, unlike for $\QoI_E$, it is clear that samples on level $\ell=0$
are of no use, since $\varMC\left(\Delta\QoI_{W,\ell}\right)$ only
become smaller than $\varMC\left(\QoI_{W,\ell}\right)$ for
$\ell\geq2$. This shows that typically the optimal MLMC approximation
should start with the discretization labeled $\ell=1$ in
Table~\ref{tab:mesh_hierarchy_2D}. 

We note from the sample variances, $\varMC(\QoI_{W,\ell})$, that the
standard deviation of $\QoI_W$ is of the order $3\e{-2}$, while
$\E{\QoI_W}$ is of the order $1\e{-1}$; see the reference value
in Table~\ref{tab:ref_hierarchy}. Thus the uncertainty in the Earth material
parameters contribute significantly to $\QoI_W$ when, as in this case,
the source location, $\source$, in the MLMC simulation is not too far
from the source location, $\source[\ast]$, used when generating the
synthetic data.

\begin{table}[!ht]
  \centering
  \begin{minipage}{1.0\linewidth}
    \centering
    \includegraphics[width=0.49\linewidth]{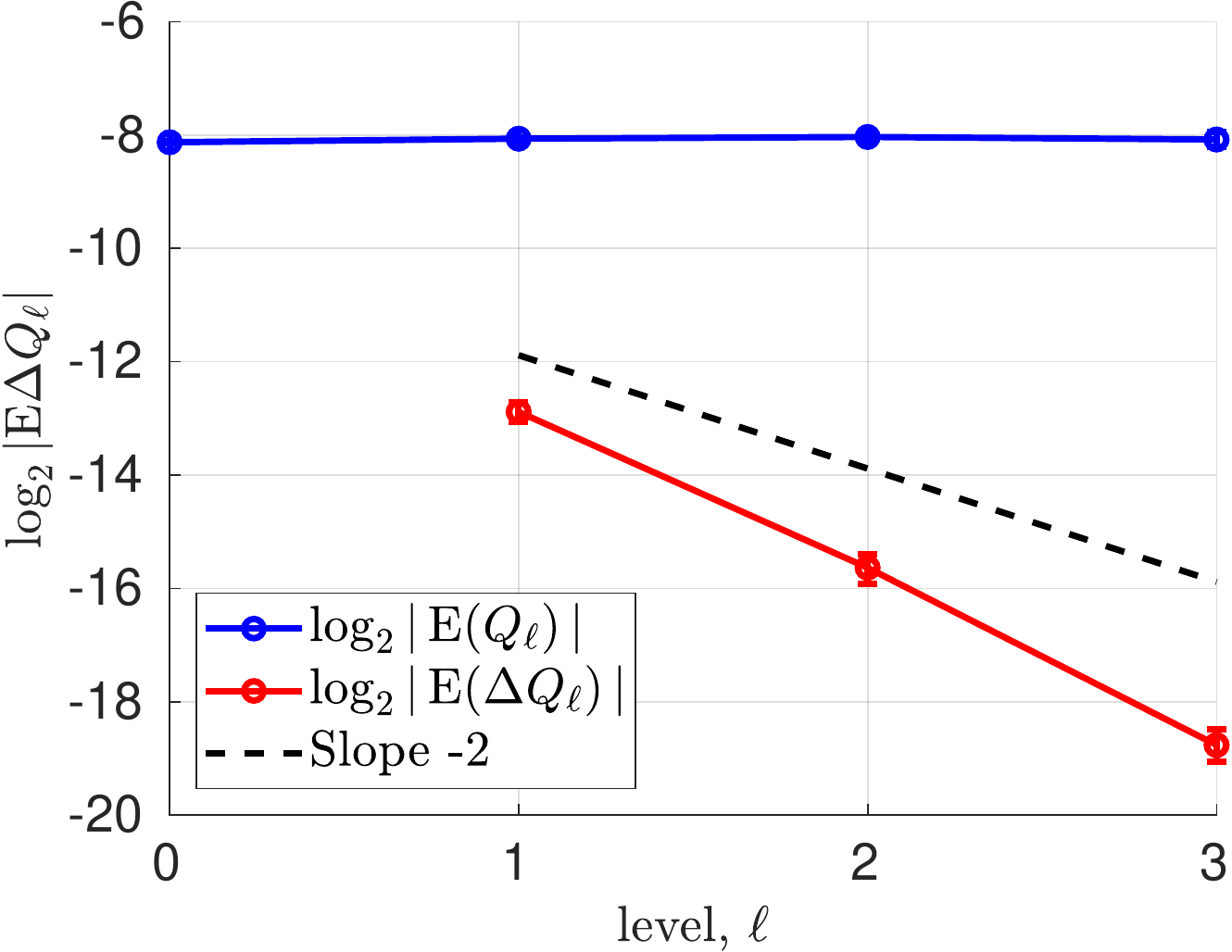}
    \includegraphics[width=0.49\linewidth]{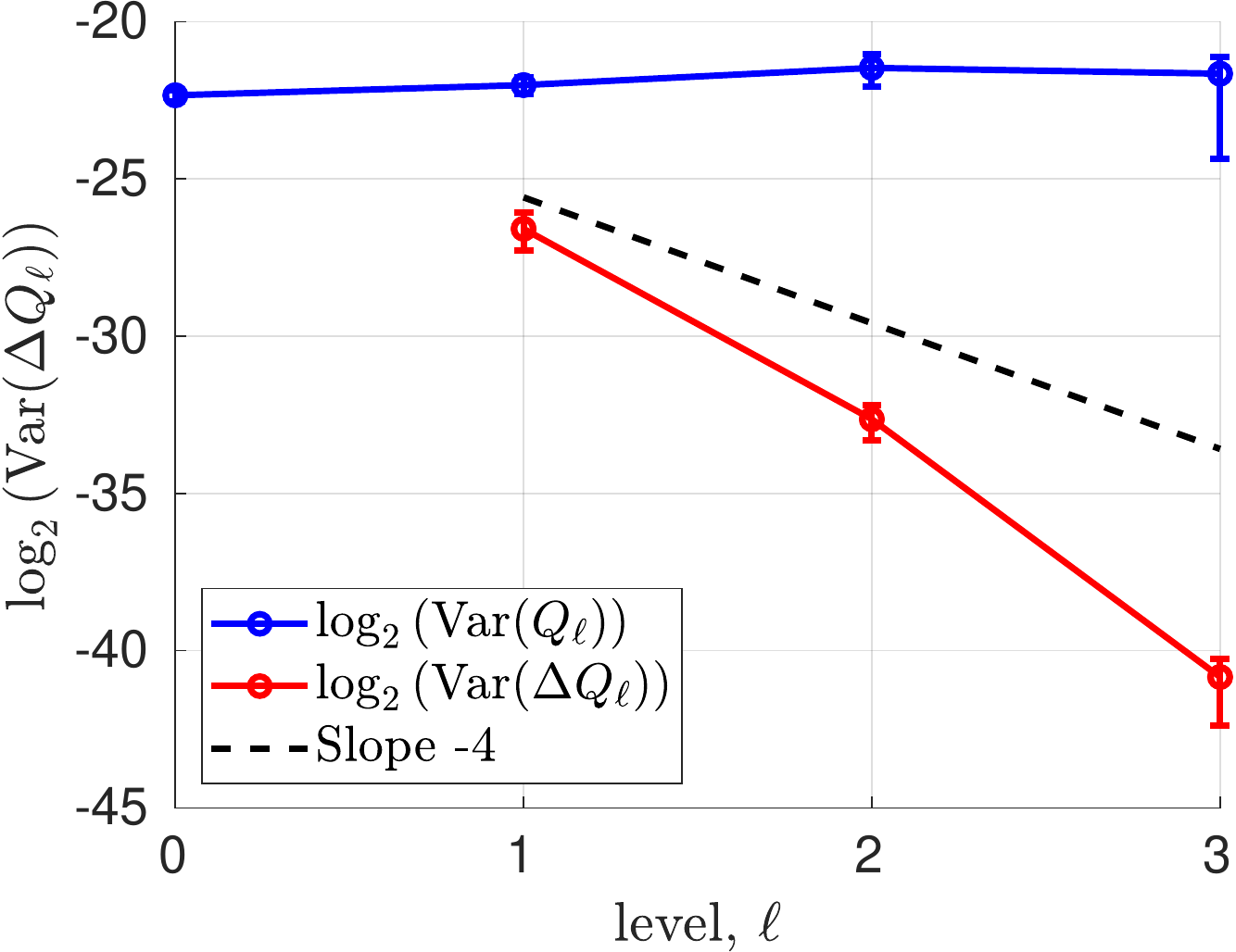}
    \captionof{figure}{Sample averages (left) and sample variances (right) of
      $\QoI_{E,\ell}$ and  $\Delta\QoI_{E,\ell}$ based on the
      verification run with the number of samples, $\nrs_\ell$, given in
      Table~\ref{tab:mesh_hierarchy_2D}. Error bars show bootstrapped
      $95\%$ confidence intervals.
    } 
    \label{fig:VerRun_conv_L2}
  \end{minipage}

  ~ \vspace{4mm} ~

  \begin{minipage}{1.0\linewidth}
    \centering
    \includegraphics[width=0.49\linewidth]{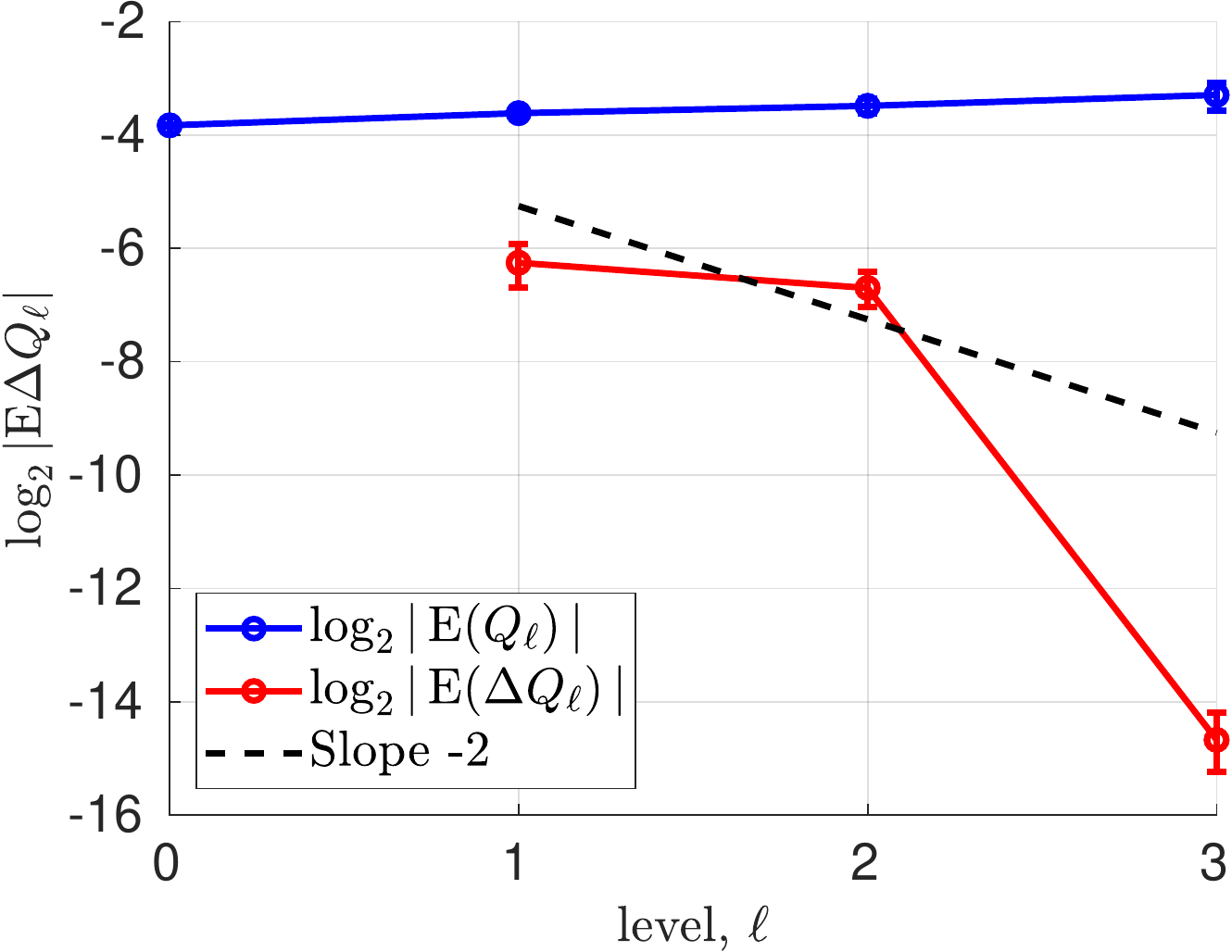}
    \includegraphics[width=0.49\linewidth]{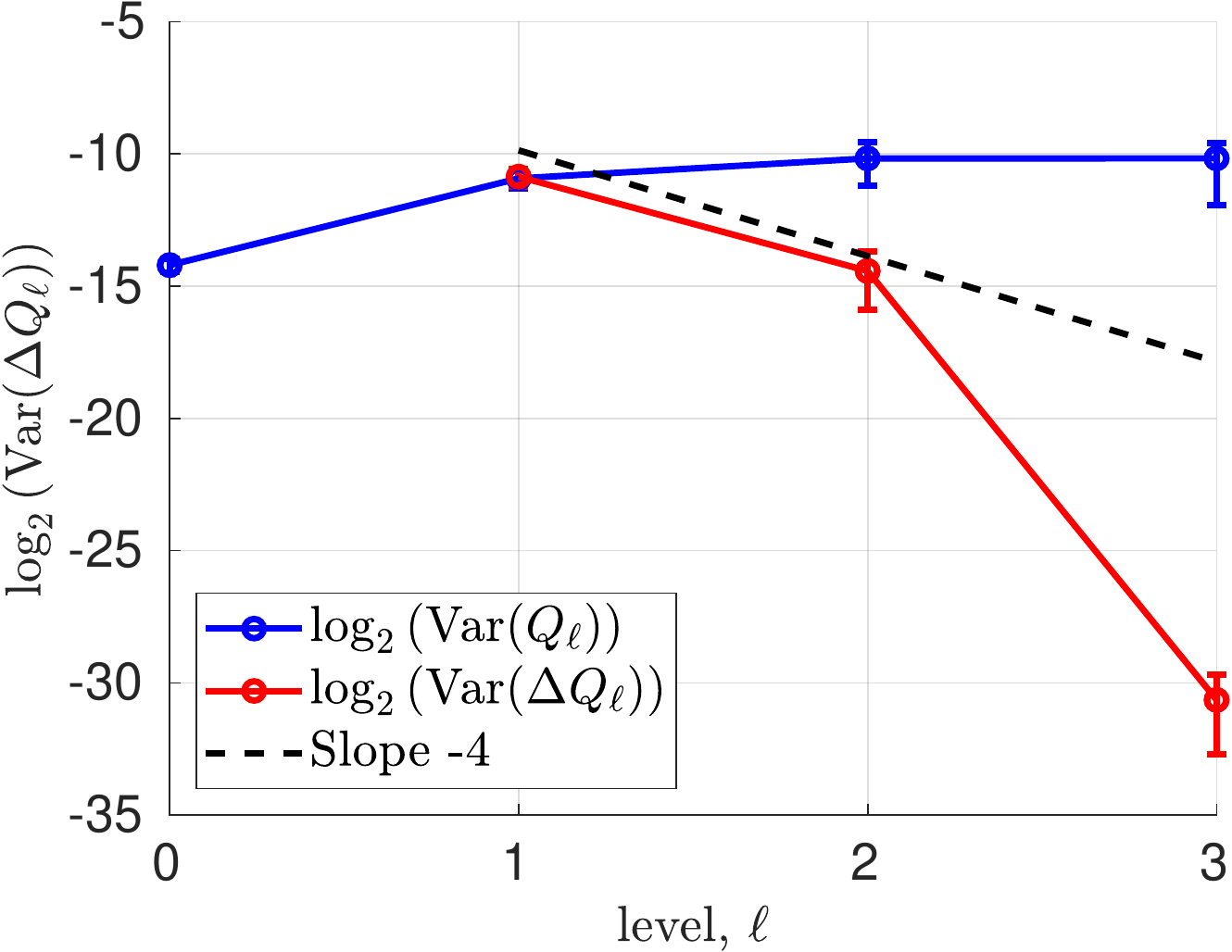}
    \captionof{figure}{Sample averages (left) and sample variances (right) of
      $\QoI_{W,\ell}$ and  $\Delta\QoI_{W,\ell}$ based on the
      verification run with the number of samples, $\nrs_\ell$, given in
      Table~\ref{tab:mesh_hierarchy_2D}. Error bars show bootstrapped
      $95\%$ confidence intervals.} 
    \label{fig:VerRun_conv_W}
  \end{minipage}
\end{table}

\subsubsection{Generation of MLMC and MC runs}
\label{sec:num_res_generate_Mx}

To test the computational complexity of generating MC and MLMC
estimators with a given tolerance, $\tol$, in $\E{\QoI_E}$ and
$\E{\QoI_W}$, we take a sequence of tolerances, $\{\tol_i\}_{i=1}^I$,
and predict which refinement levels to use, and how many samples to use
on each level to achieve an error within a given tolerance, as follows:

\paragraph{Parameters in models of work and convergence}

We take the cost per sample, cf.~\eqref{eq:work_model},
\begin{align}
  \work_\ell & =
  \begin{cases}
    \overline{\work_\ell}, & \text{for $\ell=0,\dots,3$,}\\
    \overline{\work_3}\,2^{\gamma(\ell-3)}, & \text{for $\ell>3$,}
  \end{cases}
  \label{eq:use_work}
\end{align}
where $\overline{\work_\ell}$ denotes the average core time in the
verification run and $\gamma=3$.
For the bias estimate, cf.~\eqref{eq:bias_model}, we make the
assumption that the asymptotic weak convergence rate holds for
$\ell\geq3$ and we approximate the bias on level $\ell<3$ using the
correction to level $\ell+1$. More precisely,
\begin{align}
  \left|\E{\QoI_\ast-\QoI_{\ast,\ell}}\right| & =
    \begin{cases}
      \EstMC^{95\%}(\Delta\QoI_{\ast,\ell+1}), & \text{for $\ell=0,1,2$,}\\
      \EstMC^{95\%}(\Delta\QoI_{\ast,3})\,2^{-\Ow(\ell-2)}, & \text{for $\ell>2$,}
    \end{cases}
  \label{eq:use_bias}
\end{align}
where $\EstMC^{95\%}(\Delta\QoI_{\ast,\ell})$ denotes the maximum
absolute value in the bootstrapped 95\% confidence interval of
$\EstMC(\Delta\QoI_{\ast,\ell})$. 
Similarly, for the estimate of the variances,
cf.~\eqref{eq:var_model_MC} and~\eqref{eq:var_model_MLMC}, we assume 
\begin{align}
  \vg & =
    \begin{cases}
      \varMC^{95\%}(\QoI_{E,0}), & \text{for $\QoI_E$,}\\ 
      \varMC^{95\%}(\QoI_{W,1}), & \text{for $\QoI_W$,} 
    \end{cases} 
  \label{eq:use_var0}
  \\
  \intertext{and}
  V_{\ast,\ell} & =
    \begin{cases}
      \varMC^{95\%}(\Delta\QoI_{\ast,\ell}), & \text{for $\ell=1,2,3$,}\\
      \varMC^{95\%}(\Delta\QoI_{\ast,3})\,2^{-\Os(\ell-3)}, & \text{for $\ell>3$,}
    \end{cases}
  \label{eq:use_varD}
\end{align}
where $\varMC^{95\%}(\QoI_{\ast,\ell})$ and
$\varMC^{95\%}(\Delta\QoI_{\ast,\ell})$ denote the maximum value in
the bootstrapped 95\% confidence interval of
$\varMC(\QoI_{\ast,\ell})$ and $\varMC(\Delta\QoI_{\ast,\ell})$ respectively. 

\paragraph{Tolerances used and corresponding estimators}

For the convergence tests, we estimate the scale of $\E{\QoI_\ast}$
from the verification run and choose sequences of decreasing
tolerances
\begin{align*}
  \tol_k & = \tol_1\left(\frac{1}{\sqrt{2}}\right)^k, && \text{for k=1,\dots,K,}
\end{align*}
where, for $\E{\QoI_E}$, $\tol_1 = 4.650\e{-4}\approx 12.5\%$ of
$\E{\QoI_E}$ and, for $\E{\QoI_W}$, $\tol_1 = 1.920\e{-2}\approx 21\%$
of $\E{\QoI_W}$.

To determine which levels to include and how many samples to use on
each level in the MLMC, we proceed as
follows: Given $\tol$, $\confpar=2$, $\ell_{max}$, and the
models~\eqref{eq:use_work}--\eqref{eq:use_varD}, we use the brute force
optimization described in Algorithm~\ref{alg:MLMC_hierarchy} to
determine the optimal choices 
$\mathbf{H} =
\left(\ell_0,L,\{\nrs_\ell\}_{\ell=\ell_0}^L\right)$. Here, 
$\ell_{max}$ is a relatively small positive integer, since the number
of levels grows at most logarithmically in $\tol^{-1}$. 
The resulting choices are shown in Table~\ref{tab:nr_samples_MLMC_L2}
on page~\pageref{tab:nr_samples_MLMC_L2}
and Table~\ref{tab:nr_samples_MLMC_W22}
on page~\pageref{tab:nr_samples_MLMC_W22}, for $\QoI_E$ and $\QoI_W$, 
respectively. 

For standard MC estimators, we make the analogous brute force
optimization to determine on which level to sample and how many
samples to use; see Table~\ref{tab:nr_samples_MC_L2}
and Table~\ref{tab:nr_samples_MC_W22}.

\begin{algorithm}[ht]
  \begin{algorithmic}[1]
    \State $W=\infty$
    \For{$\ell_0=0:\ell_{max}$}
      \For{$L=\ell_0:\ell_{max}$}
        \State $b \Leftarrow$ bias estimated by~\eqref{eq:use_bias}
          with $\ell=L$
        \If{$b<\tol$}
          \State $\varphi \Leftarrow 1-b/\tol$,
          cf.~\eqref{eq:splitting}
          \State $\{\nrs_\ell^\ast\}_{\ell=\ell_0}^L \Leftarrow$
            optimal samples in~\eqref{eq:opt_nrs}, summing from
            $\ell_0$ to $L$, given $\tol$,~$\confpar$,~$\varphi$  and\\
            \hspace{2.0cm}~work estimates, $\work_\ell$, 
            in~\eqref{eq:use_work}, variance estimates,
            $V_{\ell_0}$, in~\eqref{eq:use_var0}, and
            $V_{\ell}$ in~\eqref{eq:use_varD}, for $\ell>\ell_0$, 
          \State $\{\nrs_\ell\}_{\ell=\ell_0}^L \Leftarrow
            \max{\left\{2,\lceil\nrs_\ell^\ast\rceil\right\}}$ 
          \State $W^\ast \Leftarrow $ work
            estimate~\eqref{eq:work_sum}, summing from $\ell_0$
            to $L$
          \If{$W^\ast<W$}
            \State $W \Leftarrow W^\ast$
            \State $\mathbf{H} \Leftarrow \left(\ell_0,L,\{\nrs_\ell\}_{\ell=\ell_0}^L\right)$
          \EndIf  
        \EndIf  
      \EndFor
    \EndFor
  \end{algorithmic}
  \caption{Selection of optimal MLMC hierarchy}
  \label{alg:MLMC_hierarchy}
\end{algorithm}

\paragraph{Observations regarding the suggested MLMC and MC
  parameters}

Recalling that we expect, $\gamma=3$, $\Ow=2$, and $\Os=4$,
asymptotically for both $\QoI_E$ and $\QoI_W$, and that this should
lead to an asymptotic complexity
$\work_{\mathrm{MC}}\propto\tol^{-3.5}$ and
$\work_{\mathrm{MLMC}}\propto\tol^{-2}$, as $\tol\to0$, we show the
predicted work for both MC and MLMC together with these asymptotic
work rate estimates in Figure~\ref{fig:VerRun_pred_work}. In both
cases, it is clear that the predicted MLMC work grows with the
asymptotically expected rate, which is the optimal rate for Monte
Carlo type methods, as it is the same rate obtained for MC sampling
when samples can be generated at unit cost, independently of $\tol$. 

It is clear, by comparing Figure~\ref{fig:VerRun_pred_L} (showing the
refinement level) with Figure~\ref{fig:VerRun_pred_work}, that there
are ranges of values of $\tol$ for which the predicted work for MC 
first grows approximately as $\tol^{-2}$, and then faster towards the 
end of the stage. These ranges correspond to values of $\tol$
resulting in the same refinement level, so that the cost per sample
and the bias estimate, within each range, are independent of $\tol$. 
As noted above, in this pre-asymptotic regime, the apparent
convergence $\E{\QoI_{\ast,\ell}}$, with respect to $\ell$, is faster
than the asymptotically expected rate, $\Ow=2$. Therefore, the MC work
will also grow at a slower rate than the asymptotic estimate. In
particular, $\dx[]$ and $\dt[]$ are decreased at a lower rate with
decreasing $\tol$.

This faster apparent weak convergence rate is also reflected in the
value of the splitting parameter, $\splitting$
in~\eqref{eq:splitting}, implicitly obtained through the brute force
optimization in Algorithm~\ref{alg:MLMC_hierarchy} 
(Figure~\ref{fig:VerRun_pred_splitting}). 
The optimal splitting for MC, given the asymptotic rates of work per
sample, $\gamma=3$, and weak convergence, $\Ow=2$, is
$\splitting_\mathrm{MC}=4/7$ according
to~\eqref{eq:optimal_splitting_MC}, while for $\tol$ in the given
range the observed $\splitting$ is typically closer to 1 due to the
fast decay of the bias estimate.
For MLMC, in contrast, we expect $\splitting\to 1$, as $\tol\to 0$,
with , $\gamma=3$, $\Ow=2$, and $\Os=4$.

\paragraph{Predicted savings of MLMC compared to MC}

We recall from Figure~\ref{fig:VerRun_pred_work} that MLMC still
provides significant savings, compared to MC, even in the range of
tolerances where the work of MC grows at a slower rate than we can
predict that it will do asymptotically, as $\tol\to0$. For example, for
the finest tolerance in Table~\ref{tab:nr_samples_MLMC_L2}, MLMC is
predicted to reduce the work of MC by about $97\%$, and for the
finest tolerance in Table~\ref{tab:nr_samples_MLMC_W22}, by about
$78\%$.  
This is also illustrated in Figure~\ref{fig:VerRun_pred_ratio}.

\begin{table}[!ht]
  \centering
  \begin{minipage}{1.0\linewidth}
    \centering
    \includegraphics[width=0.49\linewidth]{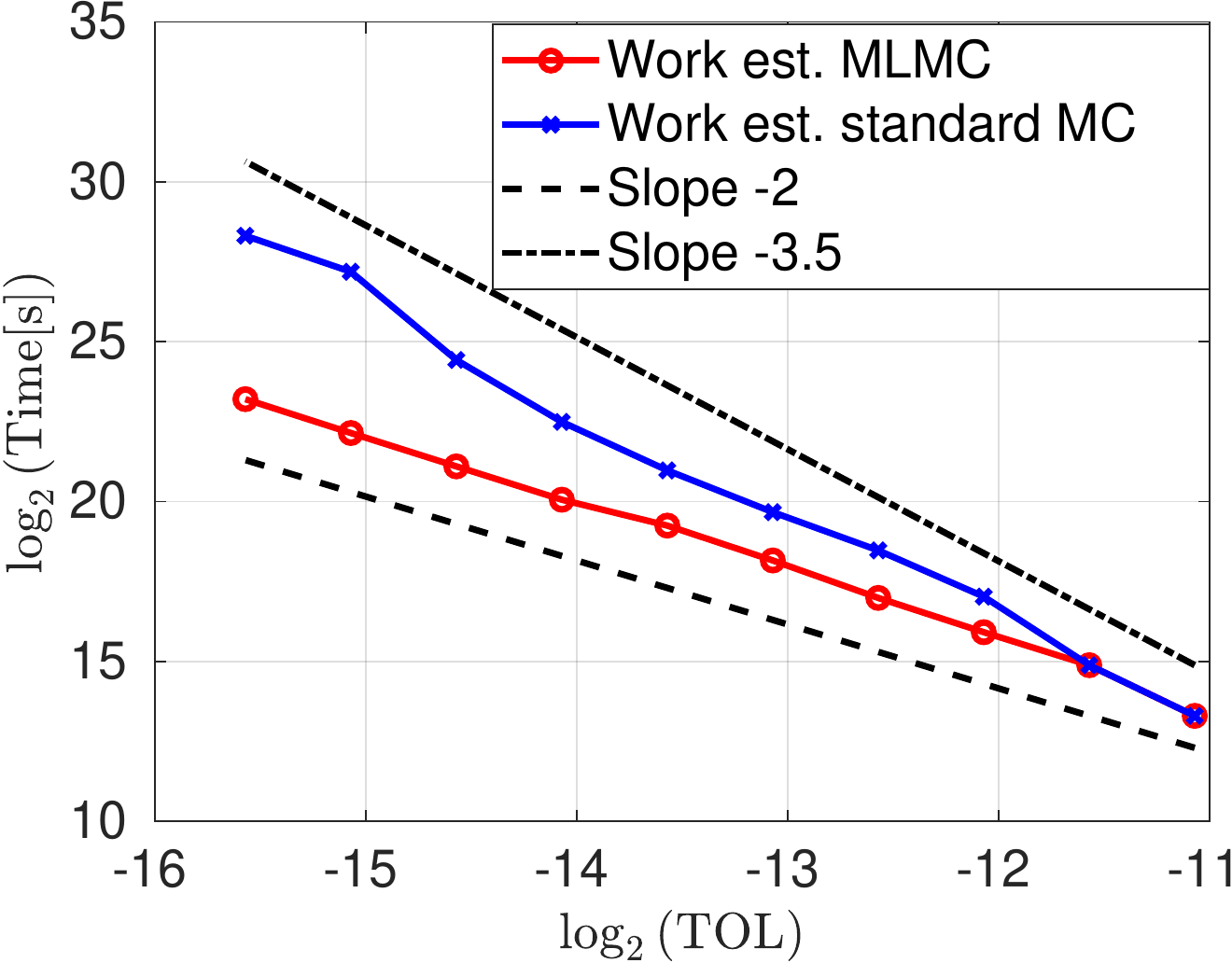}
    \includegraphics[width=0.49\linewidth]{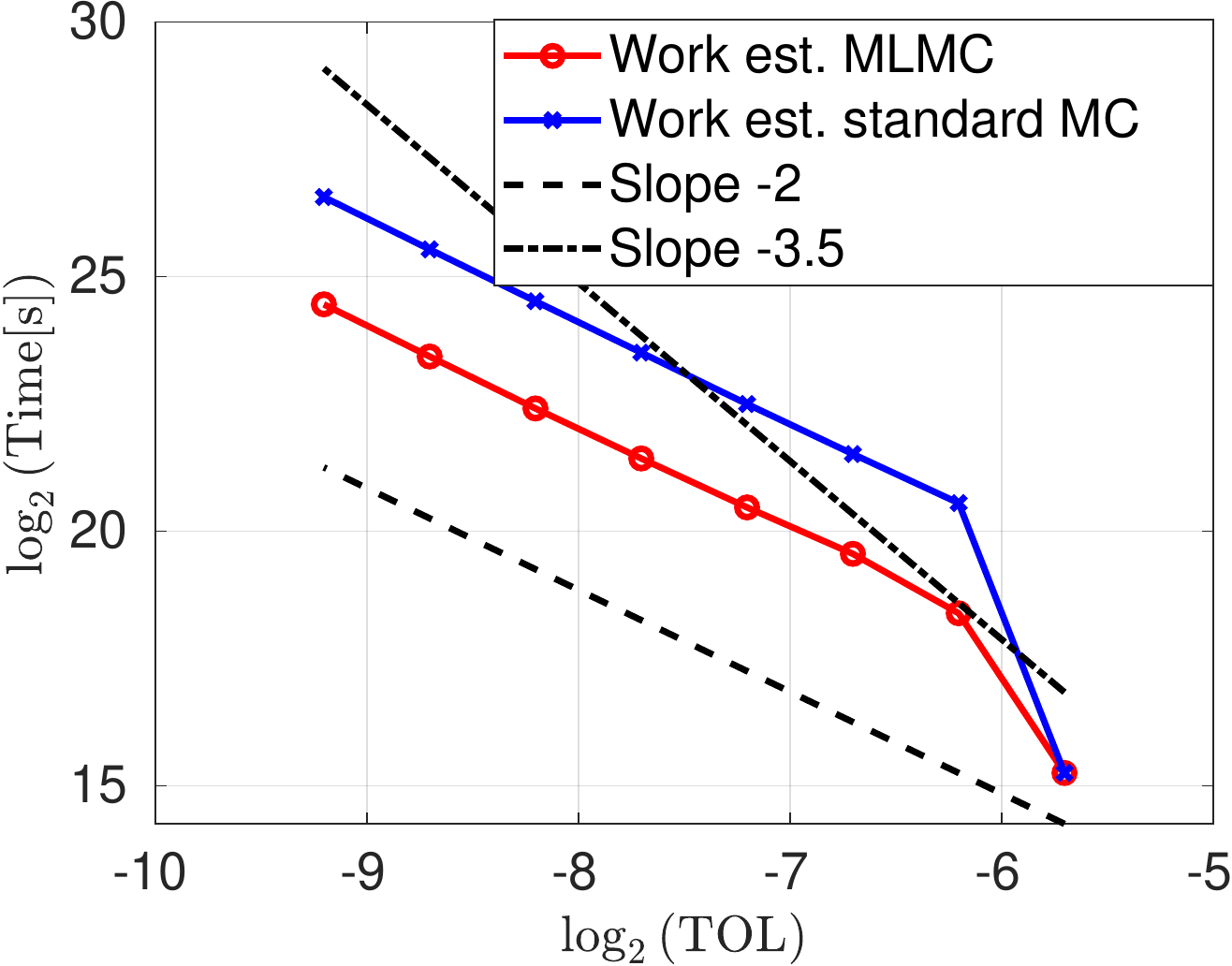}
    \captionof{figure}{Predicted work of MLMC and MC based on the verification
      run which resulted in
      Table~\ref{tab:nr_samples_MLMC_L2} for $\QoI_E$,  (Left), 
      and Table~\ref{tab:nr_samples_MLMC_W22} for $\QoI_W$ (Right).}
    \label{fig:VerRun_pred_work}
  \end{minipage}

  ~ \vspace{4mm} ~

  \begin{minipage}{1.0\linewidth}
    \centering
    \includegraphics[width=0.49\linewidth]{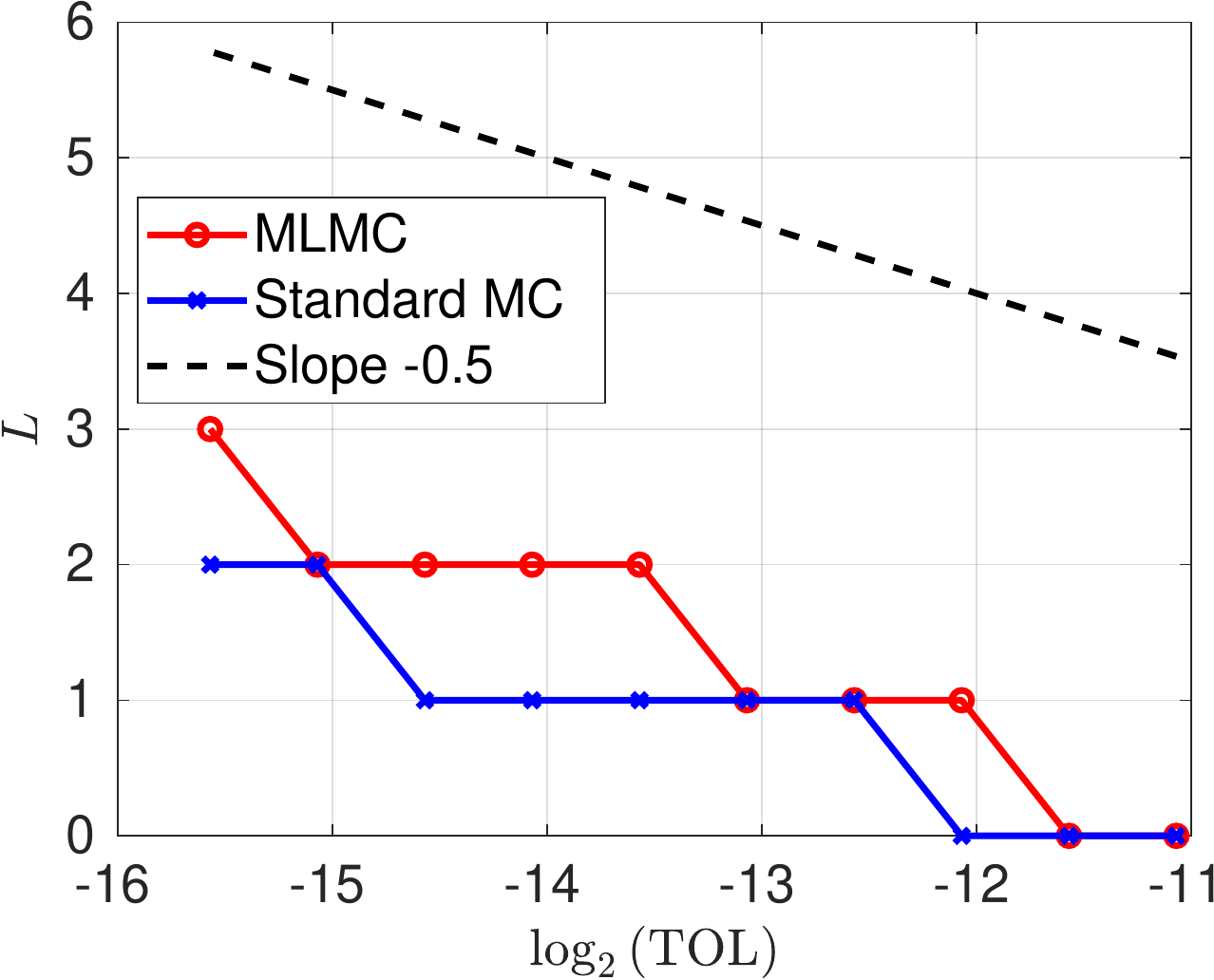}
    \includegraphics[width=0.49\linewidth]{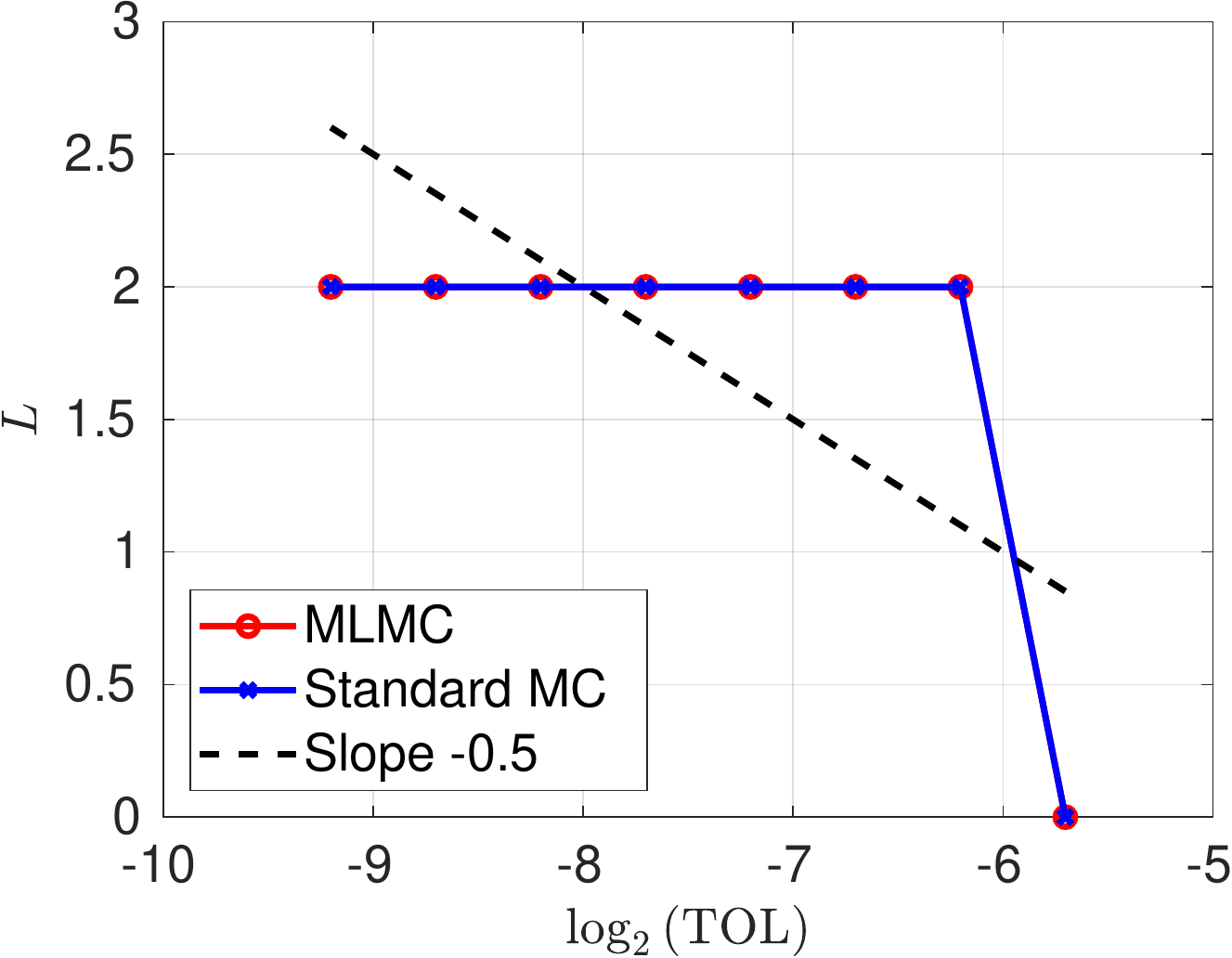}
    \captionof{figure}{The refinement level, $L$, determining the bias of the
      estimators, i.e., the maximum refinement level for MLMC and single
      refinement level for MC, for the sequence of tolerances used for
      $\QoI_E$, (Left), and  $\QoI_W$, (Right).}
    \label{fig:VerRun_pred_L}
  \end{minipage}
\end{table}

\begin{table}[!ht]
  \centering
  \begin{minipage}{1.0\linewidth}
    \centering
    \includegraphics[width=0.49\linewidth]{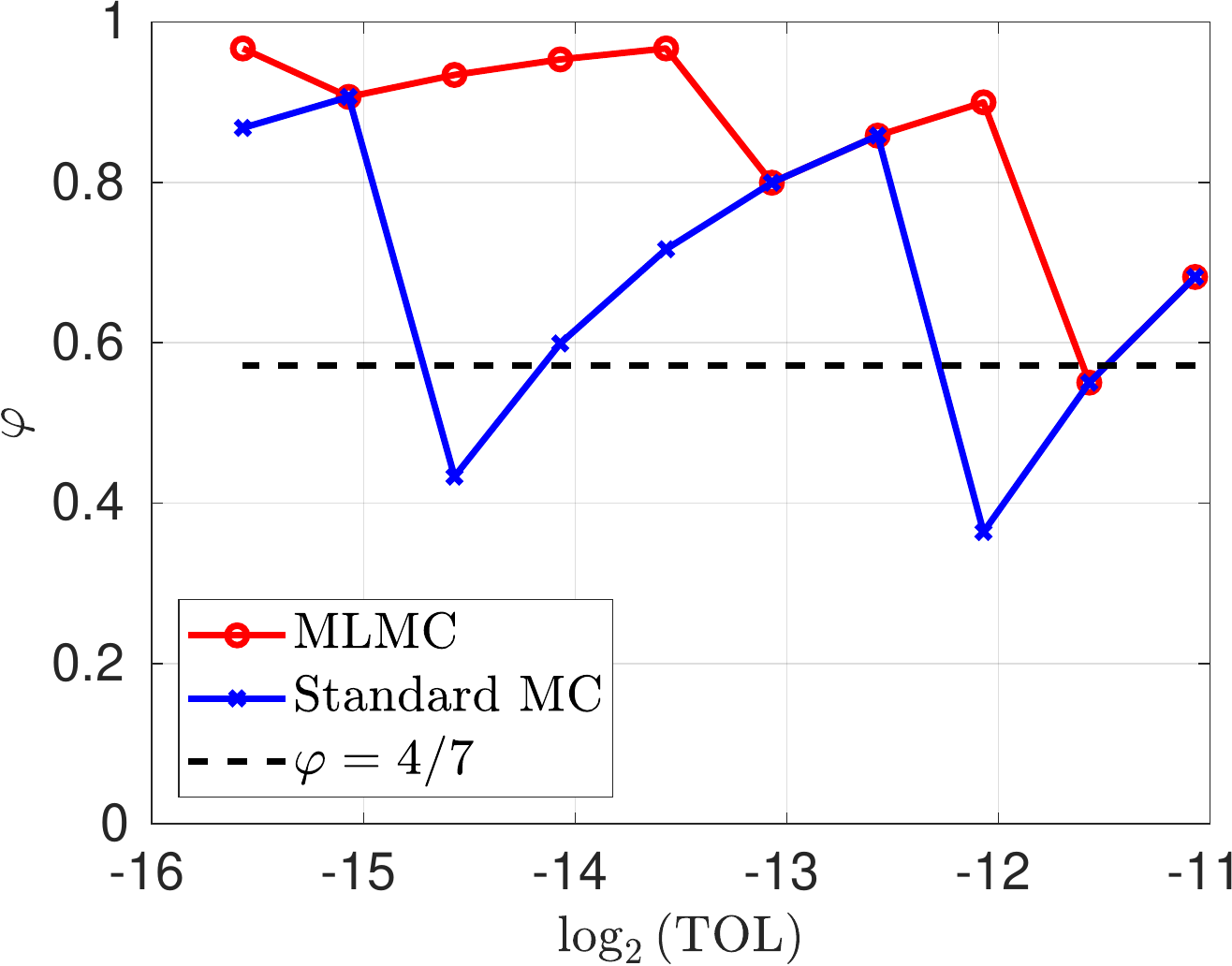}
    \includegraphics[width=0.49\linewidth]{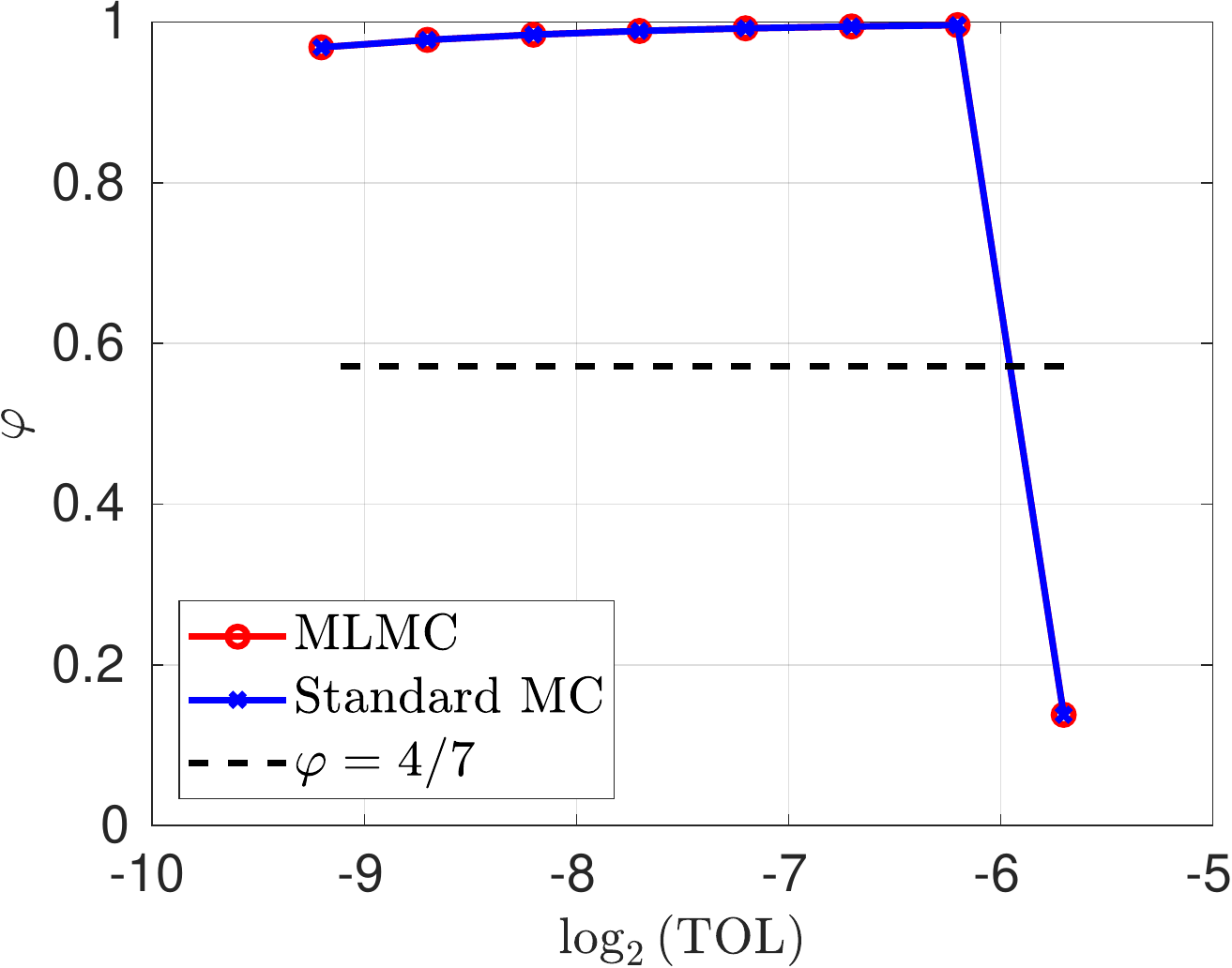}
    \captionof{figure}{Splitting parameter, $\splitting$, in~\eqref{eq:splitting},
      implicit in the MLMC and MC estimators listed in
      Table~\ref{tab:nr_samples_MLMC_L2} for  $\QoI_E$,  (Left), 
      and Table~\ref{tab:nr_samples_MLMC_W22} for $\QoI_W$, (Right).
      The dashed line denotes the asymptotically optimal value
      $\splitting=4/7$ for MC, as $\tol\to 0$, given the asymptotic work
      and convergence rates.} 
  \label{fig:VerRun_pred_splitting}
  \end{minipage}

  ~ \vspace{4mm} ~

  \begin{minipage}{1.0\linewidth}
    \centering
    \includegraphics[width=0.49\linewidth]{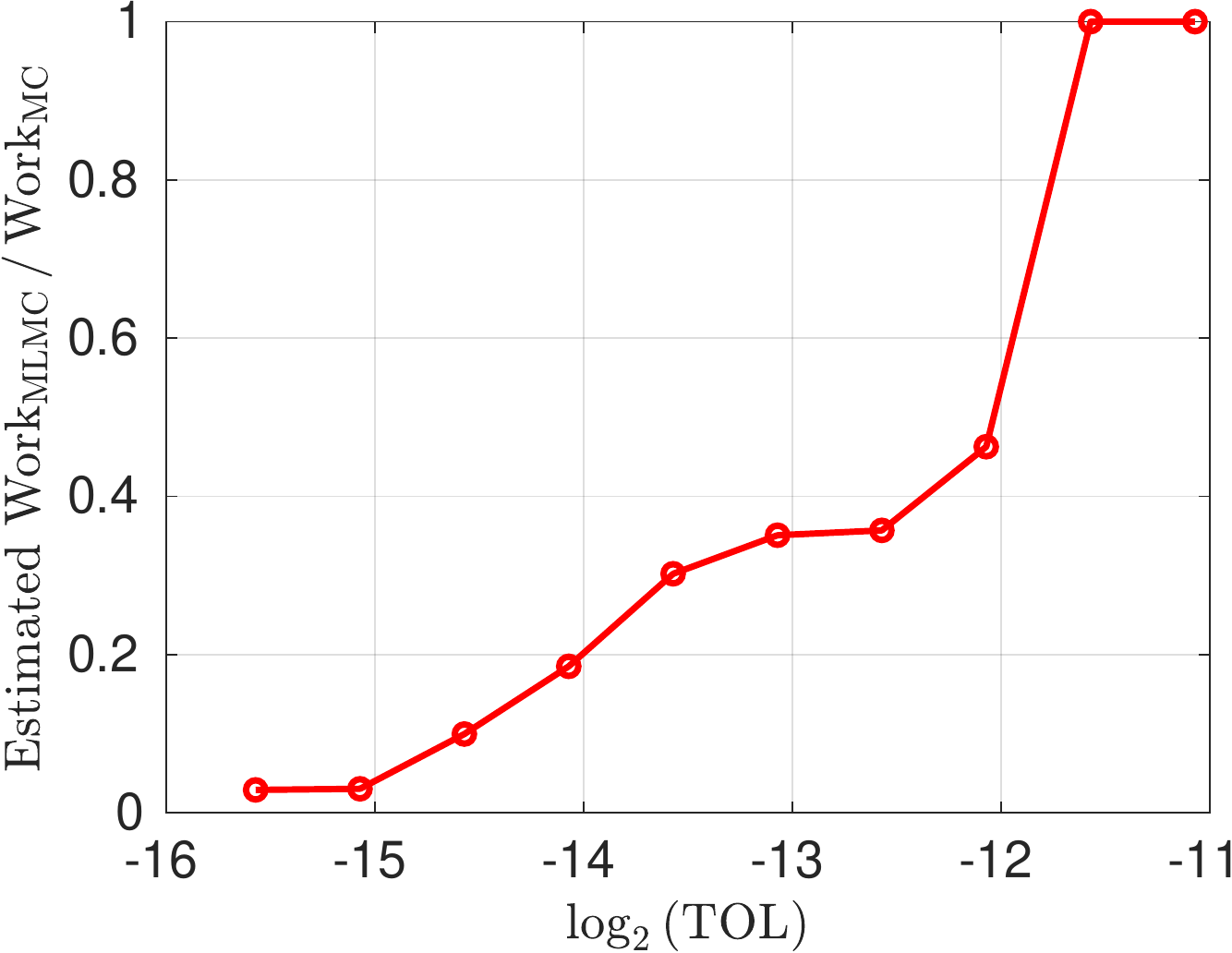}
    \includegraphics[width=0.49\linewidth]{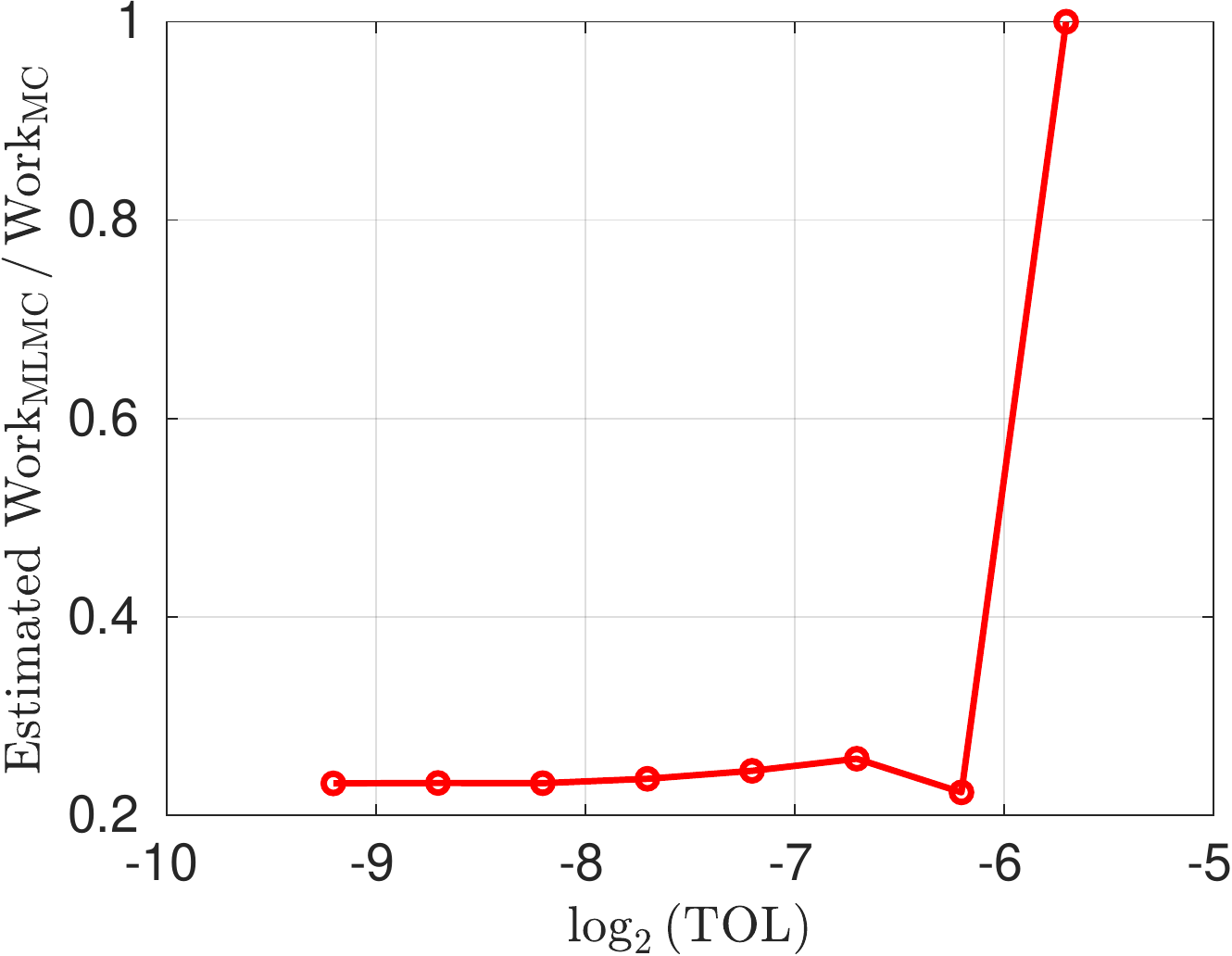}
    \captionof{figure}{Predicted ratio of MLMC to MC work based on the verification
      run which resulted in
      Table~\ref{tab:nr_samples_MLMC_L2} for $\QoI_E$,  (Left), 
      and Table~\ref{tab:nr_samples_MLMC_W22} for $\QoI_W$ (Right).}
    \label{fig:VerRun_pred_ratio}
  \end{minipage}
\end{table}

\begin{table}[ht]
  \centering
  \begin{minipage}{1.0\linewidth}
    \centering
    \begin{tabular}{ |c|r|r|r|r|c|c| }
      \hline
      \multirow{2}{*}{}
      & \multicolumn{4}{|c|}{Samples per level, $\nrs_\ell$\Tstrut}
      & \multirow{2}{*}{$\tol$} & $\work_\mathrm{MLMC}$ \\ 
      \cline{2-5}
      & $\ell=0$ & $\ell=1$ & $\ell=2$ & 
      $\ell=3$\Tstrut & & Core Time [s]\\
      \hline
      $\tol_{1}$   &    8 &   -- &  -- & --  &  4.650\e{-4}\Tstrut & 8.730\e{3}\\
      $\tol_{2}$   &   24 &   -- &  -- & --  &  3.288\e{-4}\Tstrut & 2.995\e{4}\\
      $\tol_{3}$   &   28 &    3 &  -- & --  &  2.325\e{-4}\Tstrut & 5.708\e{4}\\
      $\tol_{4}$   &   61 &    6 &  -- & --  &  1.644\e{-4}\Tstrut & 1.182\e{5}\\
      $\tol_{5}$   &  141 &   13 &  -- & --  &  1.163\e{-4}\Tstrut & 2.688\e{5}\\
      $\tol_{6}$   &  220 &   20 &   2 & --  &  8.220\e{-5}\Tstrut & 6.017\e{5}\\
      $\tol_{7}$   &  452 &   40 &   2 & --  &  5.813\e{-5}\Tstrut & 1.131\e{6}\\
      $\tol_{8}$   &  941 &   82 &   4 & --  &  4.110\e{-5}\Tstrut & 2.230\e{6}\\
      $\tol_{9}$   & 1996 &  173 &   7 & --  &  2.906\e{-5}\Tstrut & 4.619\e{6}\\
      $\tol_{10}$  & 3580 &  311 &  13 &  2  &  2.055\e{-5}\Tstrut & 9.498\e{6}\\
      \hline
    \end{tabular}
    \captionof{table}{Parameters defining the MLMC estimator for different
      tolerances and the resulting work, measured in core time, for the
      approximation of $\E{\QoI_E}$. The largest tolerance, $\tol_1$,
      corresponds to approximately 12.5\% of the reference value in 
      Table~\ref{tab:ref_hierarchy}.}
    \label{tab:nr_samples_MLMC_L2}
  \end{minipage}

  ~ \vspace{4mm} ~

  \begin{minipage}{1.0\linewidth}
    \centering
    \begin{tabular}{ |c|c|r|c|c| }
      \hline
      & $\ell$ & $\nrs_\ell$ & $\tol$ & $\work_\mathrm{MC}$, Core Time [s] \\ 
      \hline
      $\tol_{1}$   & 0 &    8 & 4.650\e{-4}\Tstrut & 8.730\e{3}\\
      $\tol_{2}$   & 0 &   24 & 3.288\e{-4}\Tstrut & 2.995\e{4}\\
      $\tol_{3}$   & 0 &  106 & 2.325\e{-4}\Tstrut & 1.182\e{5}\\
      $\tol_{4}$   & 1 &   48 & 1.644\e{-4}\Tstrut & 3.588\e{5}\\
      $\tol_{5}$   & 1 &  110 & 1.163\e{-4}\Tstrut & 8.260\e{5}\\
      $\tol_{6}$   & 1 &  272 & 8.220\e{-5}\Tstrut & 2.053\e{6}\\
      $\tol_{7}$   & 1 &  778 & 5.813\e{-5}\Tstrut & -- \\
      $\tol_{8}$   & 1 &  2978 & 4.110\e{-5}\Tstrut & -- \\
      $\tol_{9}$   & 2 &  1986 & 2.906\e{-5}\Tstrut & -- \\
      $\tol_{10}$   & 2 &  4333 & 2.055\e{-5}\Tstrut & -- \\
      \hline
    \end{tabular}
    \captionof{table}{Parameters defining the MC estimator for the approximation of
      $\E{\QoI_E}$ for different tolerances and the resulting work,
      measured in core time, in the cases where the MC estimate has been
      computed.}
    \label{tab:nr_samples_MC_L2}
  \end{minipage}
\end{table}

\begin{table}[ht]
  \centering
  \begin{minipage}{1.0\linewidth}
    \centering
    \begin{tabular}{ |c|r|r|r|r|c|c| }
      \hline
      \multirow{2}{*}{}
      & \multicolumn{4}{|c|}{Number of samples per level, $\nrs_\ell$\Tstrut}
      & \multirow{2}{*}{$\tol$} & $\work_\mathrm{MLMC}$ \\ 
      \cline{2-5}
      & \hspace{2mm}$\ell=0$ & \hspace{2mm}$\ell=1$ & \hspace{2mm}$\ell=2$ & 
      \hspace{2mm}$\ell=3$\Tstrut & & Core Time [s]\\
      \hline
      $\tol_1$ & 31 &   -- &  -- & --  & 1.920\e{-2}\Tstrut & 3.396\e{4}\\
      $\tol_2$ & -- &   23 &   2 & --  & 1.357\e{-2}\Tstrut & 3.373\e{5}\\
      $\tol_3$ & -- &   46 &   5 & --  & 9.598\e{-3}\Tstrut & 7.652\e{5}\\
      $\tol_4$ & -- &   91 &   9 & --  & 6.787\e{-3}\Tstrut & 1.491\e{6}\\
      $\tol_5$ & -- &  183 &  17 & --  & 4.799\e{-3}\Tstrut & 2.790\e{6}\\
      $\tol_6$ & -- &  368 &  33 & --  & 3.393\e{-3}\Tstrut & 5.521\e{6}\\
      $\tol_7$ & -- &  744 &  67 & --  & 2.399\e{-3}\Tstrut & 1.120\e{7}\\
      $\tol_8$ & -- & 1519 & 136 & --  & 1.697\e{-3}\Tstrut & 2.302\e{7}\\
      \hline
      $\tol_9$ & -- & 2928 & 261 &  2  & 1.200\e{-3}\Tstrut & --\\
      \hline
    \end{tabular}
    \captionof{table}{Parameters defining the MLMC estimator for different
      tolerances and the resulting work, measured in core time, for the
      approximation of $\E{\QoI_W}$. The smallest tolerance, $\tol_9$,
      is included among the bootstrapped MLMC estimators, but not the
      primary realizations. The largest tolerance, $\tol_1$,
      corresponds to approximately 21\% of the reference value in 
      Table~\ref{tab:ref_hierarchy}.}
    \label{tab:nr_samples_MLMC_W22}
  \end{minipage}

  ~ \vspace{4mm} ~

  \begin{minipage}{1.0\linewidth}
    \centering
    \begin{tabular}{ |c|c|r|c|c| }
      \hline
      & $\ell$ & $\nrs_\ell$ & $\tol$ & $\work_\mathrm{MC}$, Core Time [s] \\ 
      \hline
      $\tol_1$   & 0 &    31 & 1.920\e{-2}\Tstrut & 3.396\e{4}\\
      $\tol_2$   & 2 &    20 & 1.357\e{-2}\Tstrut & 1.513\e{6}\\
      $\tol_3$   & 2 &    39 & 9.598\e{-3}\Tstrut & 2.968\e{6}\\
      $\tol_4$   & 2 &    77 & 6.787\e{-3}\Tstrut & 6.012\e{6}\\
      $\tol_5$   & 2 &   155 & 4.799\e{-3}\Tstrut & -- \\
      $\tol_6$   & 2 &   312 & 3.393\e{-3}\Tstrut & -- \\
      $\tol_7$   & 2 &   632 & 2.399\e{-3}\Tstrut & -- \\
      $\tol_8$   & 2 &  1288 & 1.697\e{-3}\Tstrut & -- \\
      \hline     
    \end{tabular}
    \captionof{table}{Parameters defining the MC estimator for the approximation of
      $\E{\QoI_W}$ for different tolerances and the resulting work,
      measured in core time, in the cases where the MC estimate has been
      computed.}
    \label{tab:nr_samples_MC_W22}
  \end{minipage}
\end{table}

\subsubsection{MLMC and MC Runs}
\label{sec:MLMC_tests}

Here, we present computational results based on the actual MLMC and MC
runs performed with the parameters listed in
Table~\ref{tab:nr_samples_MLMC_L2}, for $\QoI_E$, and
Table~\ref{tab:nr_samples_MLMC_W22}, for $\QoI_W$.

\paragraph{On the use of parameters estimated in the verification run}

Note that we use information from the verification run when we set up
the convergence tests. This is in line with the intended use of MLMC
in the inverse problem setting that involves repeatedly computing
approximate solutions to the underlying forward problem with different
parameter values in the course of solving the inverse problem, so that
prior information about parameters from earlier runs becomes available. 
Additionally, a continuation type algorithm~\cite{haji_CMLMC}, can be
used in the inverse problem setting.

\paragraph{Computational results}

For these tests, one sample of $\EstMLMC$ was computed for each
tolerance. Note that $\EstMLMC$ by itself is a random variable and
that here the samples corresponding to different tolerances are
independent. 

\begin{figure}[!ht]
  \centering
  \includegraphics[width=0.49\linewidth]{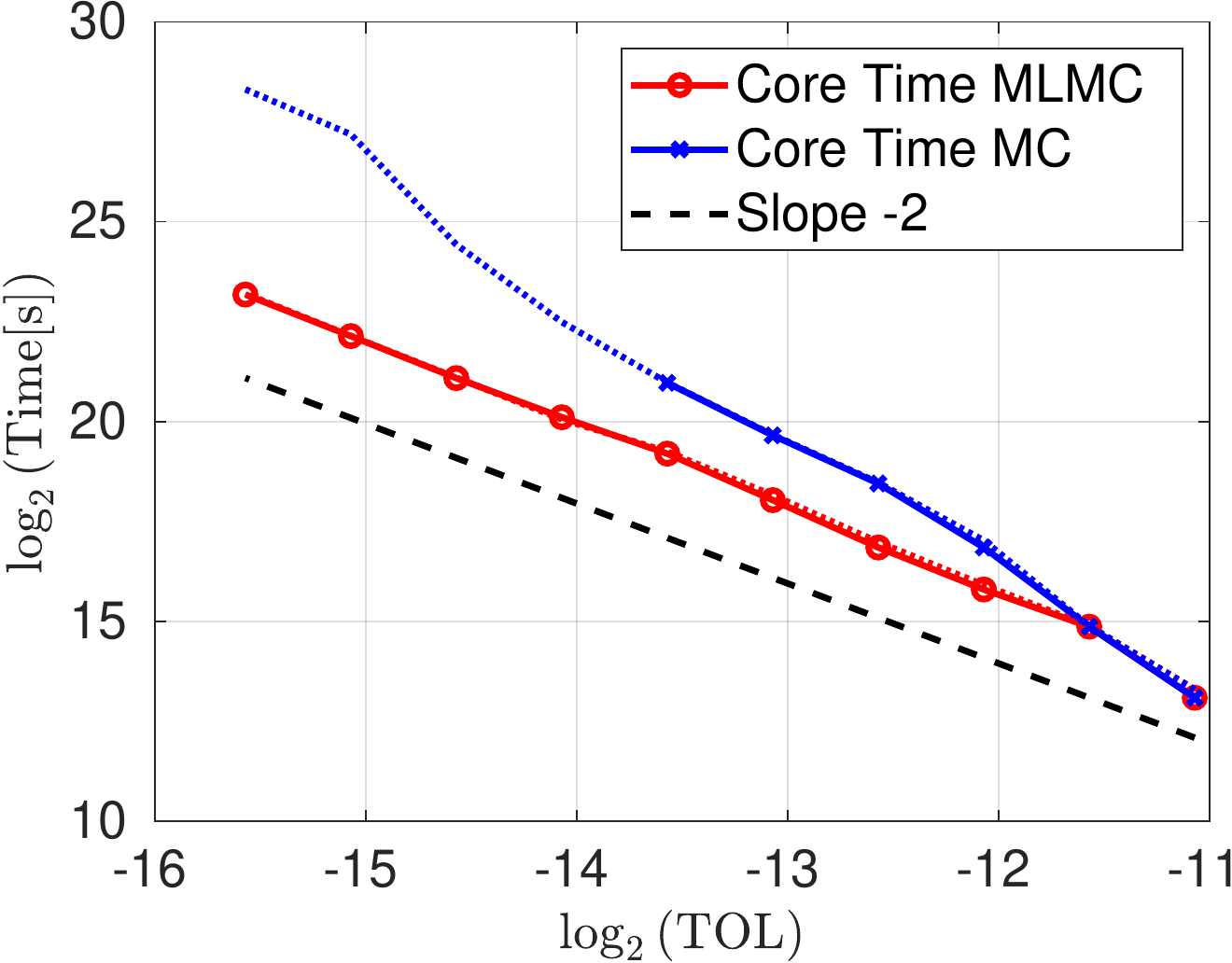}
  \includegraphics[width=0.49\linewidth]{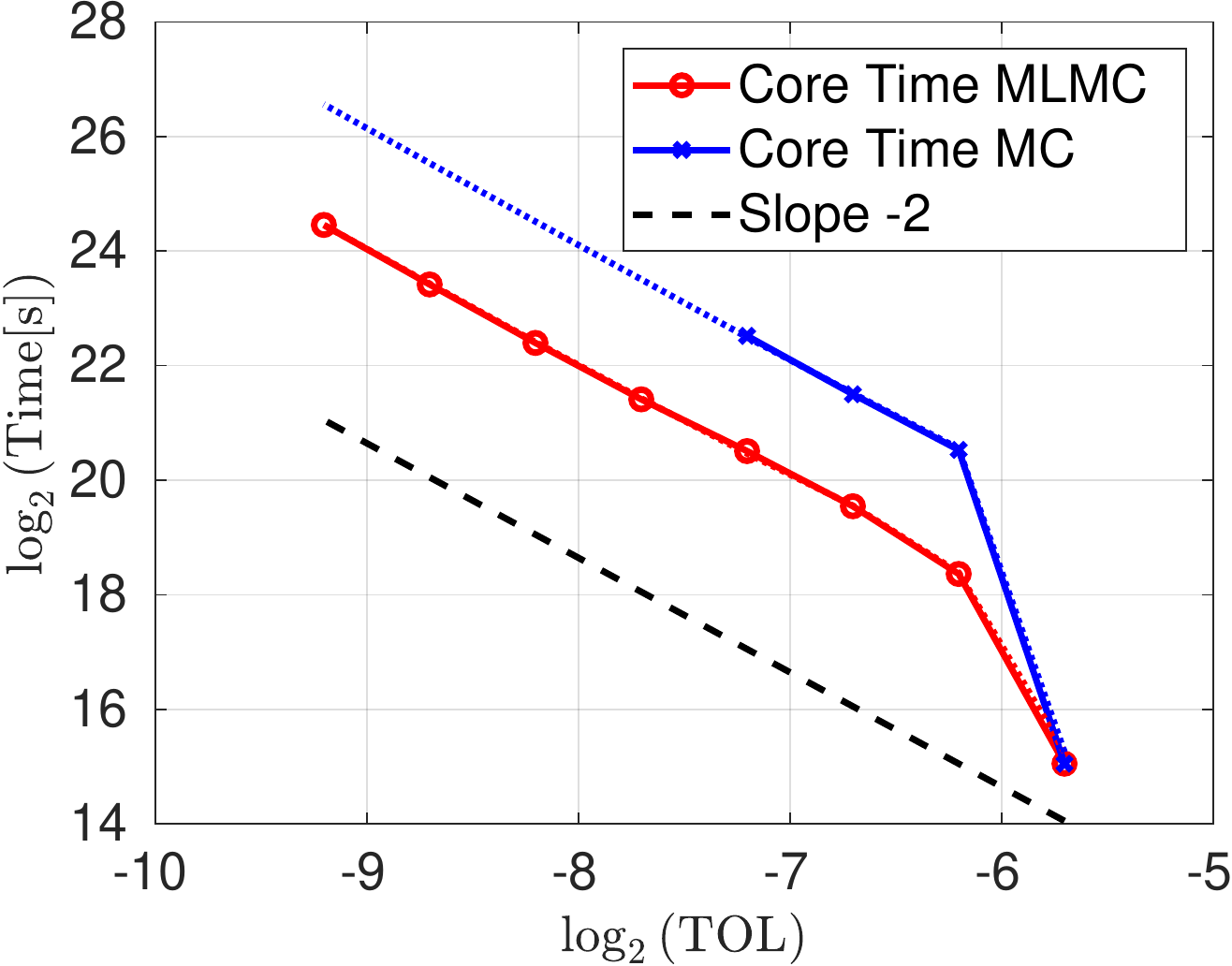}
  \caption{MLMC and MC work as a function of $\tol$ for the
    approximation of $\E{\QoI_E}$ (left) and $\E{\QoI_W}$
    (right). Solid lines show actual computational time, as described
    in Section~\ref{sec:ver_run}, and the dotted line shows the
    predicted work based on the verification run. The dashed line
    shows the slope of the optimal complexity for Monte Carlo type
    methods, $\tol^{-2}$.} 
  \label{fig:MLMC_work}
\end{figure}

The computational work, shown in Figure~\ref{fig:MLMC_work}, agrees
very well with the work predicted in
Section~\ref{sec:num_res_generate_Mx}, due to the highly consistent
execution time of \specfem and the fact that the number of samples on
each level was fixed beforehand, based on the verification run
results. 
For those tolerances where both MC and MLMC estimates were computed,
significant savings of computational time for MLMC relative to MC was
observed, as discussed in Section~\ref{sec:num_res_generate_Mx}.

In the absence of an a priori known exact solution to the test
problem, we estimate the accuracy of the MLMC results by comparing
them to a reference solution obtained by pooling a larger number of
samples on each level, \emph{including} all samples used to
generate the MLMC estimators for varying tolerances; see
Table~\ref{tab:ref_hierarchy} on page~\pageref{tab:ref_hierarchy}. 
Thus, while the
samples of these estimators are mutually independent, they are not
independent of the reference value. On the other hand, the number of
samples used to obtain the reference value vastly exceeds the number of
samples for larger tolerances and significantly exceeds the number of
samples used for the smaller tolerances. The errors compared to this
reference solution are shown as red circles in
Figure~\ref{fig:MLMC_conv}. 

Additionally, 100 samples of $\EstMLMC$ for each value of $\tol$ were
obtained by bootstrapping from the same pool of samples used to
generate the reference solution. The corresponding errors, marked
with black crosses in Figure~\ref{fig:MLMC_conv}, indicate the
variability of the error. 

\begin{table}[!ht]
  \centering
  \begin{minipage}{1.0\linewidth}
    \centering
    \includegraphics[width=0.49\linewidth]{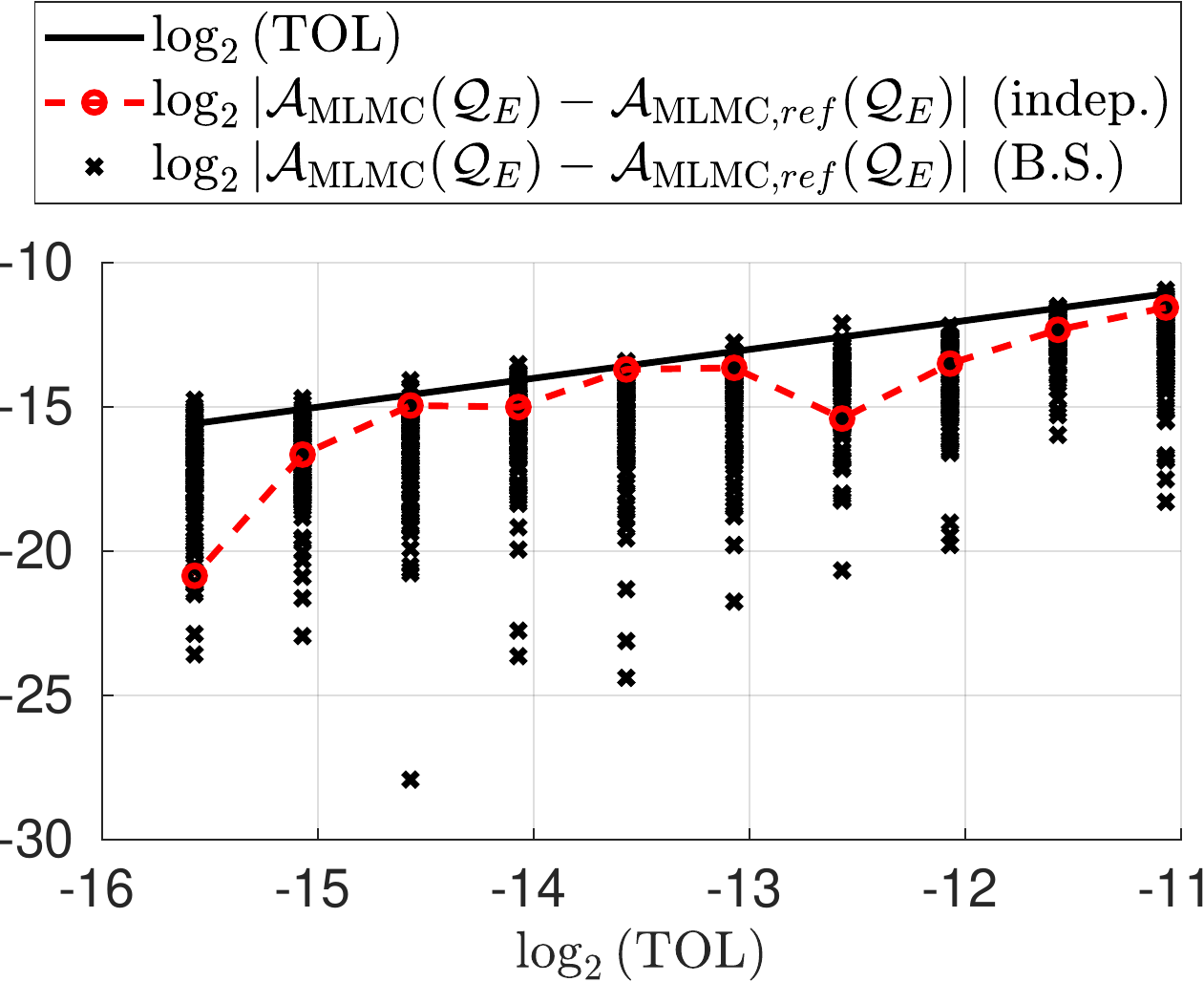}
    \includegraphics[width=0.49\linewidth]{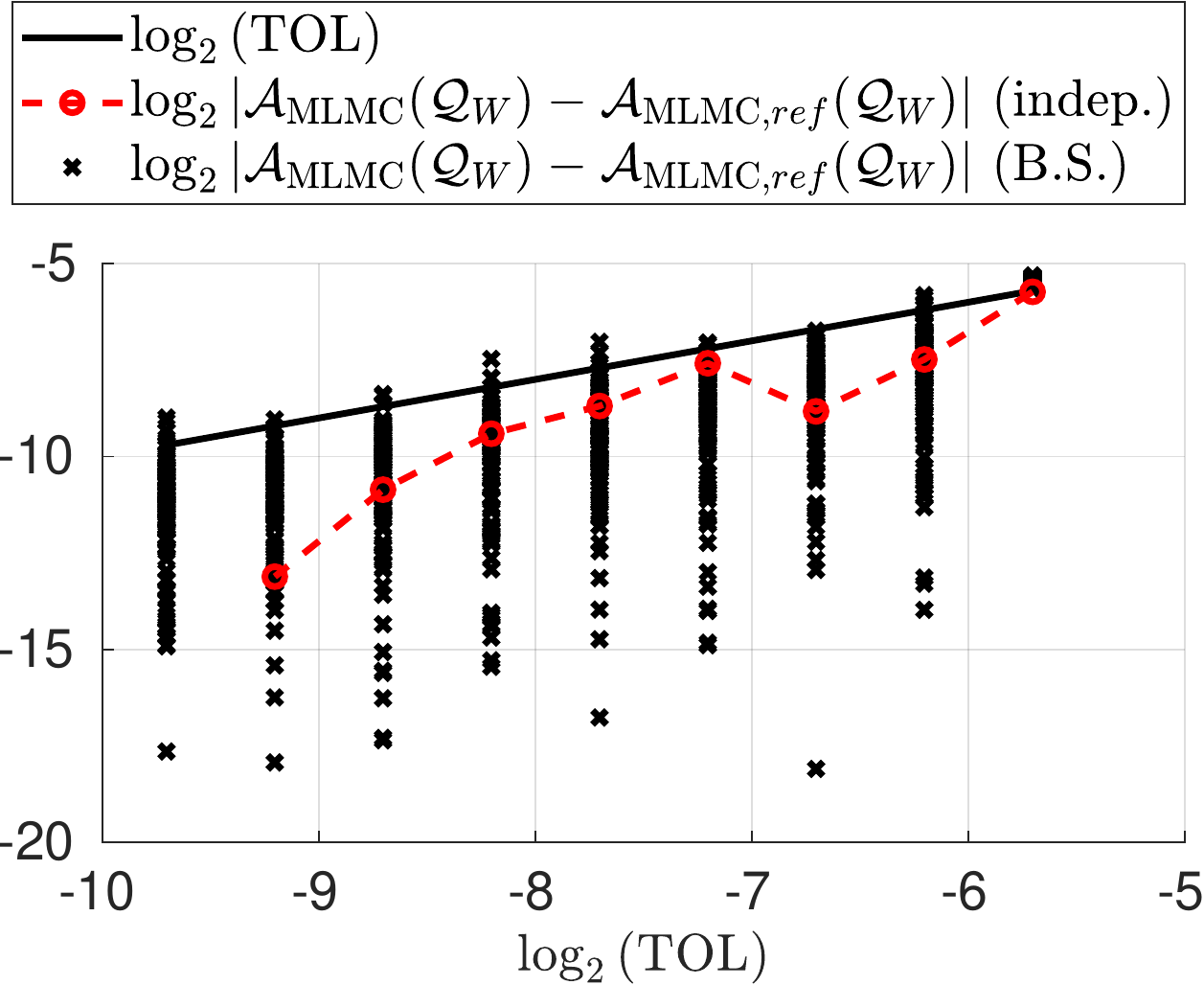}
    \captionof{figure}{Convergence of the MLMC estimators of
      Table~\ref{tab:nr_samples_MLMC_L2} for $\E{\QoI_E}$ (left) and
      Table~\ref{tab:nr_samples_MLMC_W22} for $\E{\QoI_W}$ (right). One
      realization of 
      $\EstMLMC(\QoI_\ast)$ per value of the tolerance was computed,
      based on samples independent of those used for all other
      tolerances, and the error was approximated using the reference value
      of $\EstMLMC(\QoI_\ast)$ in Table~\ref{tab:ref_hierarchy}; this
      error estimate is labeled (indep.). 
      In addition, 100 statistically dependent realizations of
      $\EstMLMC(\QoI_\ast)$ for all tolerances were obtained by
      bootstrapping with replacement from the whole pool of samples; see
      Table~\ref{tab:pool}. These bootstrapped error estimates are
      labeled (B.S.).
    }
  \label{fig:MLMC_conv}
  \end{minipage}

  ~ \vspace{4mm} ~

  \begin{minipage}{1.0\linewidth}
    \centering
    \includegraphics[width=0.49\linewidth]{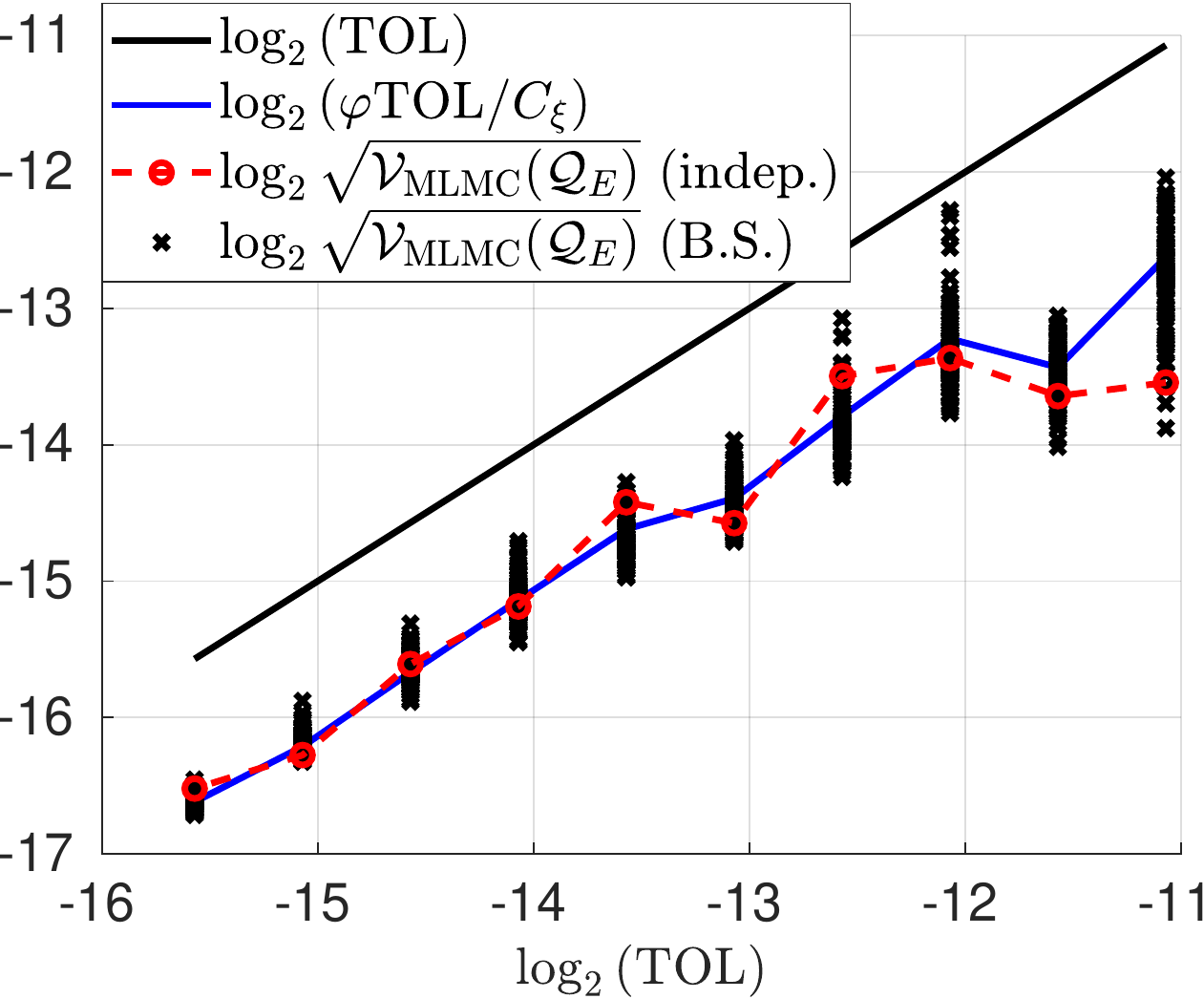}
    \includegraphics[width=0.49\linewidth]{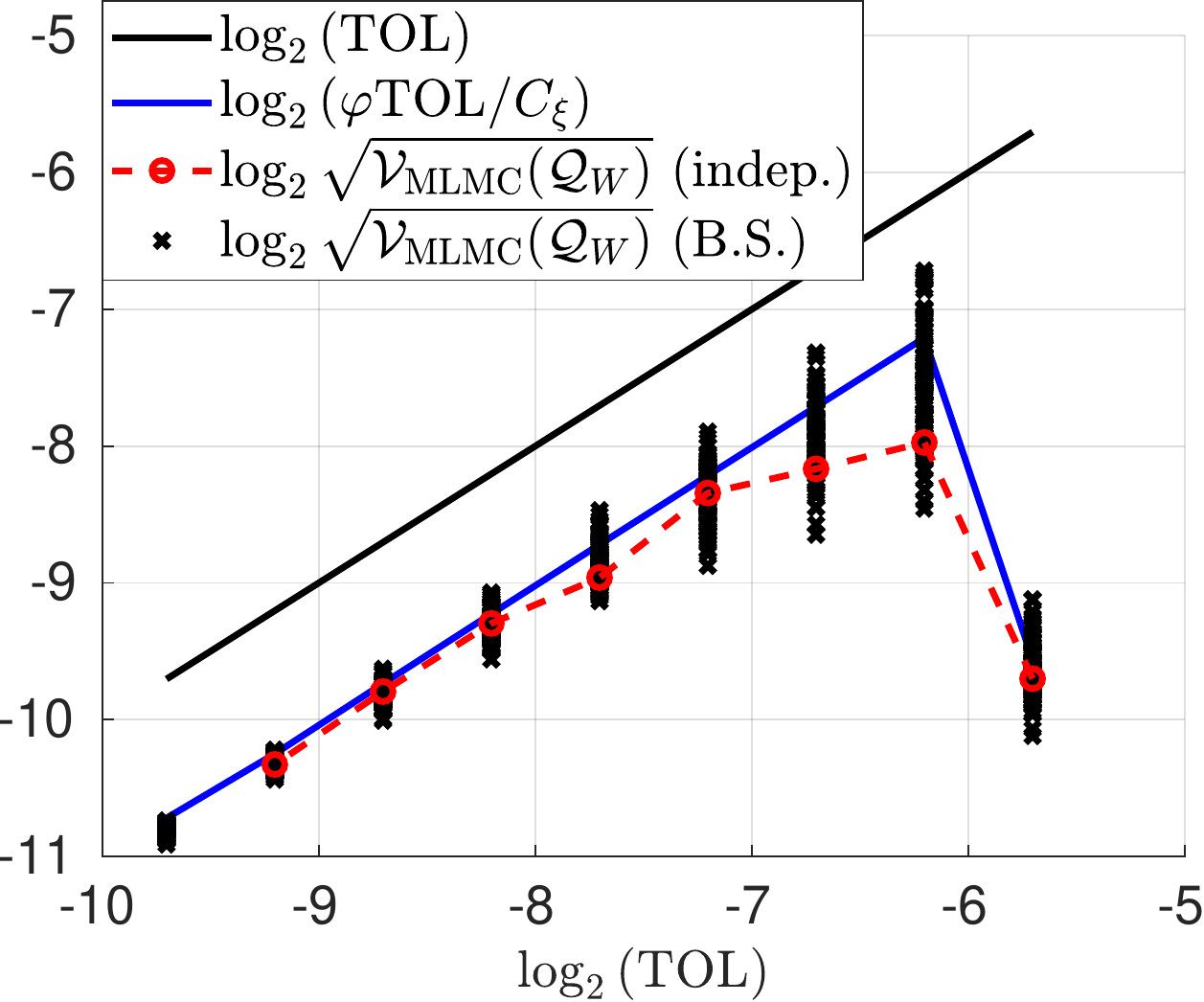}
    \captionof{figure}{Convergence of the statistical error estimate corresponding
      to convergence study in Figure~\ref{fig:MLMC_conv}. Here
      $\mathcal{V}_{MLMC}$ refers to the estimator of the variance of 
      $\EstMLMC(\QoI_\ast)$ in~\eqref{eq:var_MLMC_Est}, obtained by replacing
      the true variances by their unbiased estimators~\eqref{eq:varMC}.}
    \label{fig:Stat_Err_MLMC}
  \end{minipage}
\end{table}

The variance of the MLMC estimator
\begin{align}
  \var{\EstMLMC(\QoI)} & = \frac{1}{\nrs_0}\var{\QoI_0} + 
     \sum_{\ell=1}^L\frac{1}{\nrs_\ell}\var{\Delta\QoI_\ell},
  \label{eq:var_MLMC_Est}
\end{align}
is approximated by sample variances to verify
that~\eqref{eq:var_requirement_MLMC} is satisfied. As shown in
Figure~\ref{fig:Stat_Err_MLMC},
$\sqrt{\mathcal{V}_{MLMC}}\approx\frac{\splitting\tol}{\confpar}$,
with $\confpar=2$, consistently with the required constraint on
the statistical error, as to be expected based on the sample variance
estimates from the verification phase.

For comparison, the corresponding tests for standard MC are shown in
Figure~\ref{fig:MC_conv} and Figure~\ref{fig:Stat_Err_MC}. Note in
particular, that though the computation here were significantly more
expensive than the corresponding MLMC computations, on smaller
tolerances, the statistical error was not over-resolved.

\begin{table}[ht]
  \centering
  \begin{minipage}{1.0\linewidth}
    \centering
    \includegraphics[width=0.49\linewidth]{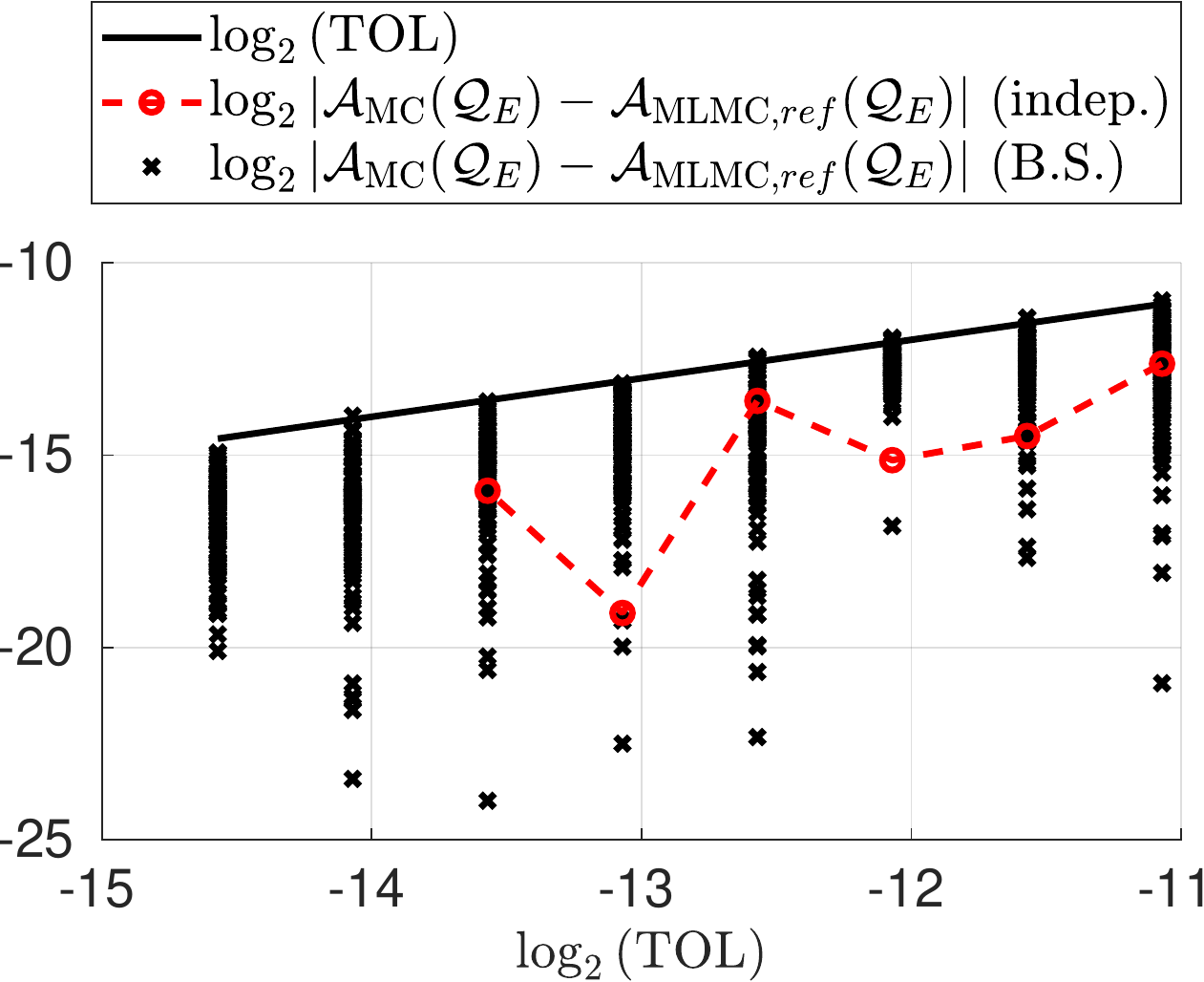}
    \includegraphics[width=0.49\linewidth]{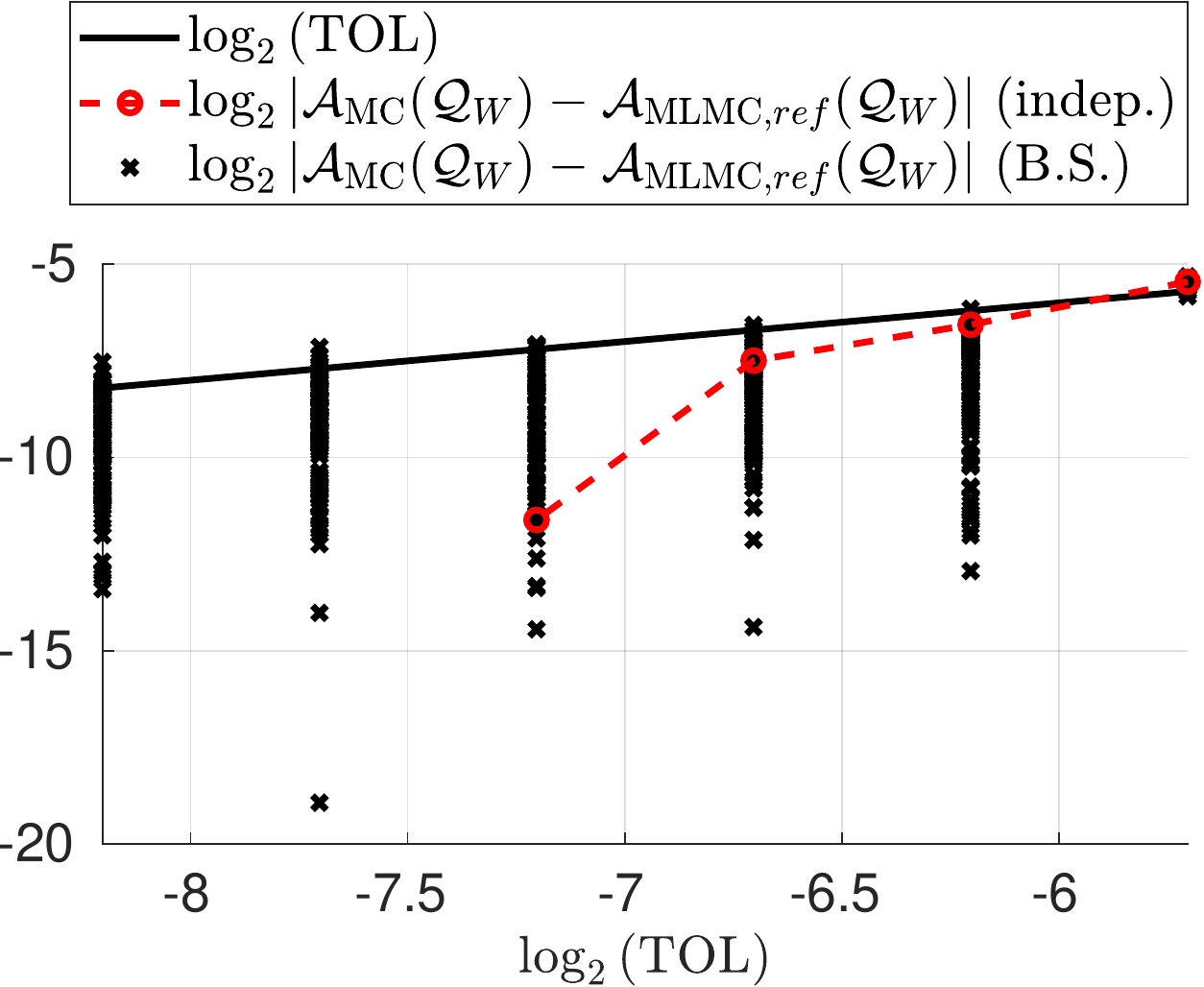}
    \captionof{figure}{Convergence of the MC estimators of
      Table~\ref{tab:nr_samples_MC_L2} for $\E{\QoI_E}$ (left) and
      Table~\ref{tab:nr_samples_MC_W22} for $\E{\QoI_W}$ (right). One
      realization of 
      $\EstMC(\QoI_\ast)$ per value of the tolerance was computed,
      based on samples independent of those used for all other
      tolerances, and the error was approximated using the reference value
      of $\EstMLMC(\QoI_\ast)$ in Table~\ref{tab:ref_hierarchy}; this
      error estimate is labeled (indep.). 
      In addition, 100 statistically dependent realizations of
      $\EstMC(\QoI_\ast)$ for all tolerances were obtained by
      bootstrapping with replacement from the whole pool of samples; see 
      Table~\ref{tab:pool}. These bootstrapped error estimates are
      labeled (B.S.).}
    \label{fig:MC_conv}
  \end{minipage}

  ~ \vspace{4mm} ~

  \begin{minipage}{1.0\linewidth}
    \centering
    \includegraphics[width=0.49\linewidth]{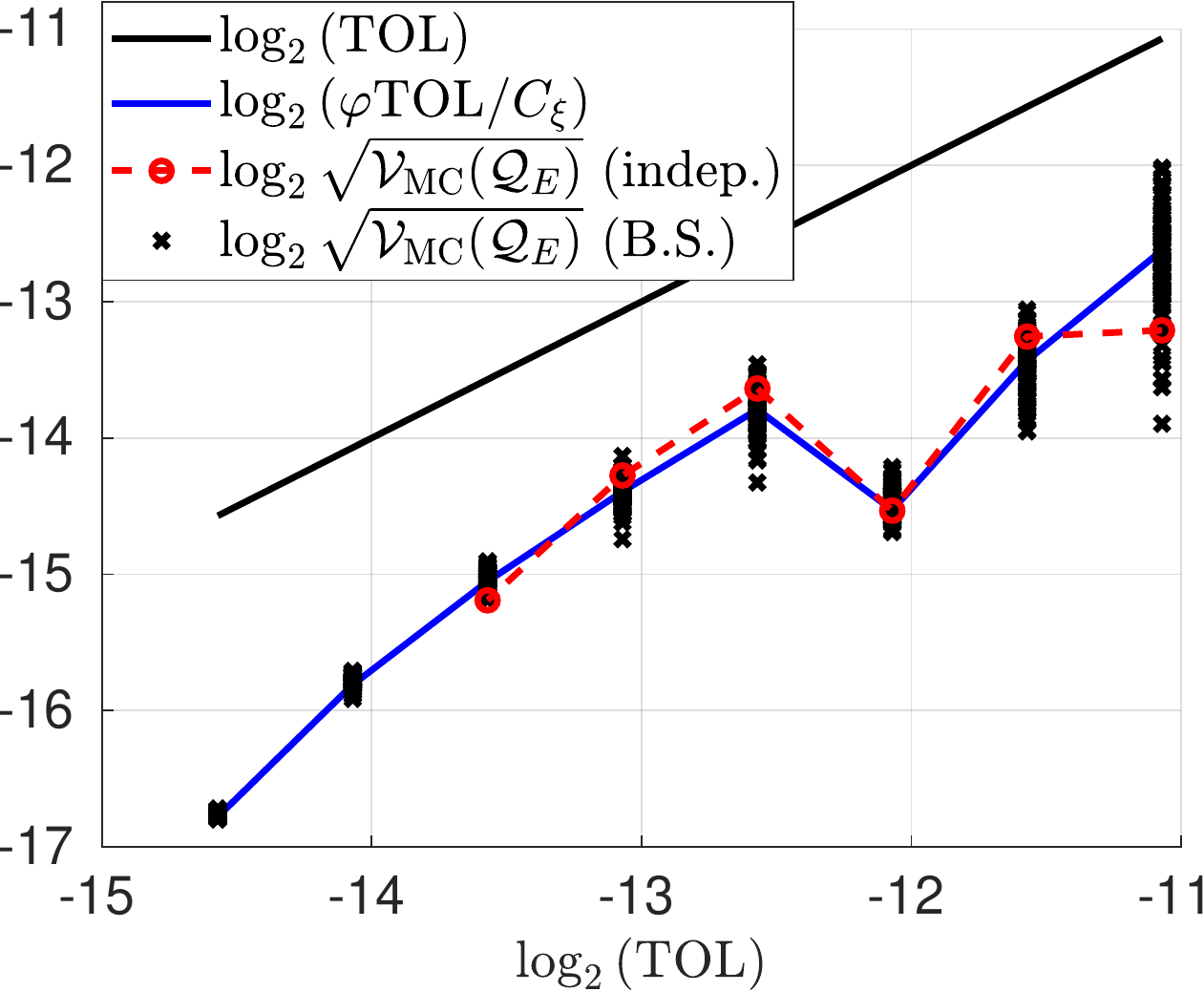}
    \includegraphics[width=0.49\linewidth]{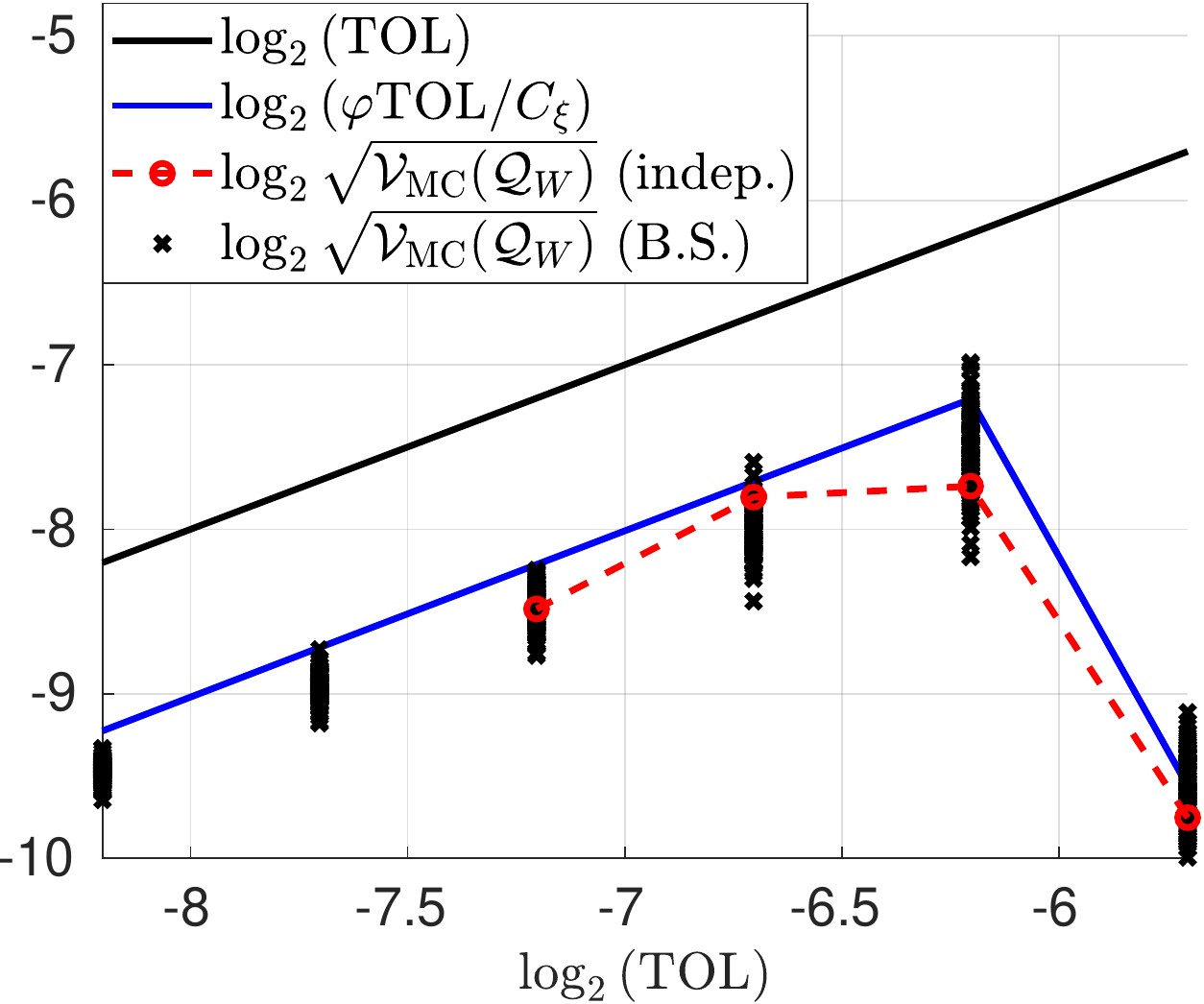}
    \captionof{figure}{Convergence of the statistical error estimate corresponding
      to convergence study in Figure~\ref{fig:MC_conv}. Here
      $\mathcal{V}_{MC}$ refers to the estimator of
      $\var{\EstMC(\QoI_{\ast,L})}=\frac{\var{(\QoI_{\ast,L})}}{\nrs_L}$,
      obtained by unbiased sample variance estimators~\eqref{eq:varMC}.
      To control the statistical error, we chose the number of samples
      so that $\sqrt{\mathcal{V}_{MC}}=\splitting\tol/\confpar$, for
      $\mathcal{V}_{MC}$ \emph{predicted} based on the parameters
      estimated from the verification run.}
    \label{fig:Stat_Err_MC}
  \end{minipage}
\end{table}

\begin{table}[ht]
  \centering
  \begin{minipage}{1.0\linewidth}
    \centering
    \begin{tabular}{ |c|r|r|r|r|c|c| }
      \hline
      \multirow{2}{*}{} & 
      \multicolumn{4}{|c|}{Number of samples per level, $\nrs_\ell$\Tstrut}
      & \multirow{2}{*}{$\EstMLMC(\QoI_{\ast,ref})$} &
      \multirow{2}{*}{$\varMC_{\mathrm{MLMC}}\left(\QoI_{\ast,ref}\right)$}
      \\ 
      \cline{2-5}
      & \hspace{2mm}$\ell=0$ & \hspace{2mm}$\ell=1$ & \hspace{2mm}$\ell=2$ & 
      \hspace{2mm}$\ell=3$\Tstrut & 
      & \\
      \hline
      $\QoI_E$ & 19045 & 5608 & 351 & 4 & 3.729\e{-3}\Tstrut & 1.31\e{-11} \\
      $\QoI_W$ & --    & 5608 & 351 & 4 & 9.227\e{-2}\Tstrut & 1.98\e{-7} \\
      \hline
    \end{tabular}
    \captionof{table}{Reference values of $\E{\QoI_E}$ and $\E{\QoI_W}$ together
      with the samples per level in the MLMC estimators used to obtain
      them.}
    \label{tab:ref_hierarchy}
  \end{minipage}

  ~ \vspace{4mm} ~

  \begin{minipage}{1.0\linewidth}
    \centering
    \begin{tabular}{ |c|r|r|r|r| }
      \hline
      \multicolumn{5}{|c|}{Number of samples per level, $\nrs_\ell$,
      in pool\Tstrut}\\ 
      \hline
      & \hspace{2mm}$\ell=0$ & \hspace{2mm}$\ell=1$ & \hspace{2mm}$\ell=2$ & 
      \hspace{2mm}$\ell=3$\Tstrut\\
      \hline
      MLMC $\QoI_E$ & 19045 & 5608 & 351 & 4 \\
      MLMC $\QoI_W$ & -- & 5608 & 351 & 4 \\
      MC & 24653 & 5959 & 355 & 4 \\
      \hline
    \end{tabular}
    \captionof{table}{The number of i.i.d. samples per level in the pool of
      samples used when bootstrapping the estimators $\EstMLMC(\QoI_E)$
      and $\EstMLMC(\QoI_W)$, where the used number of samples per level
      are given in Table~\ref{tab:nr_samples_MLMC_L2} and
      Table~\ref{tab:nr_samples_MLMC_W22} for $\QoI_E$ and $\QoI_W$
      respectively. }
    \label{tab:pool}
  \end{minipage}

  ~ \vspace{4mm} ~

  \begin{minipage}{1.0\linewidth}
    \centering
    \includegraphics[width=0.49\linewidth]{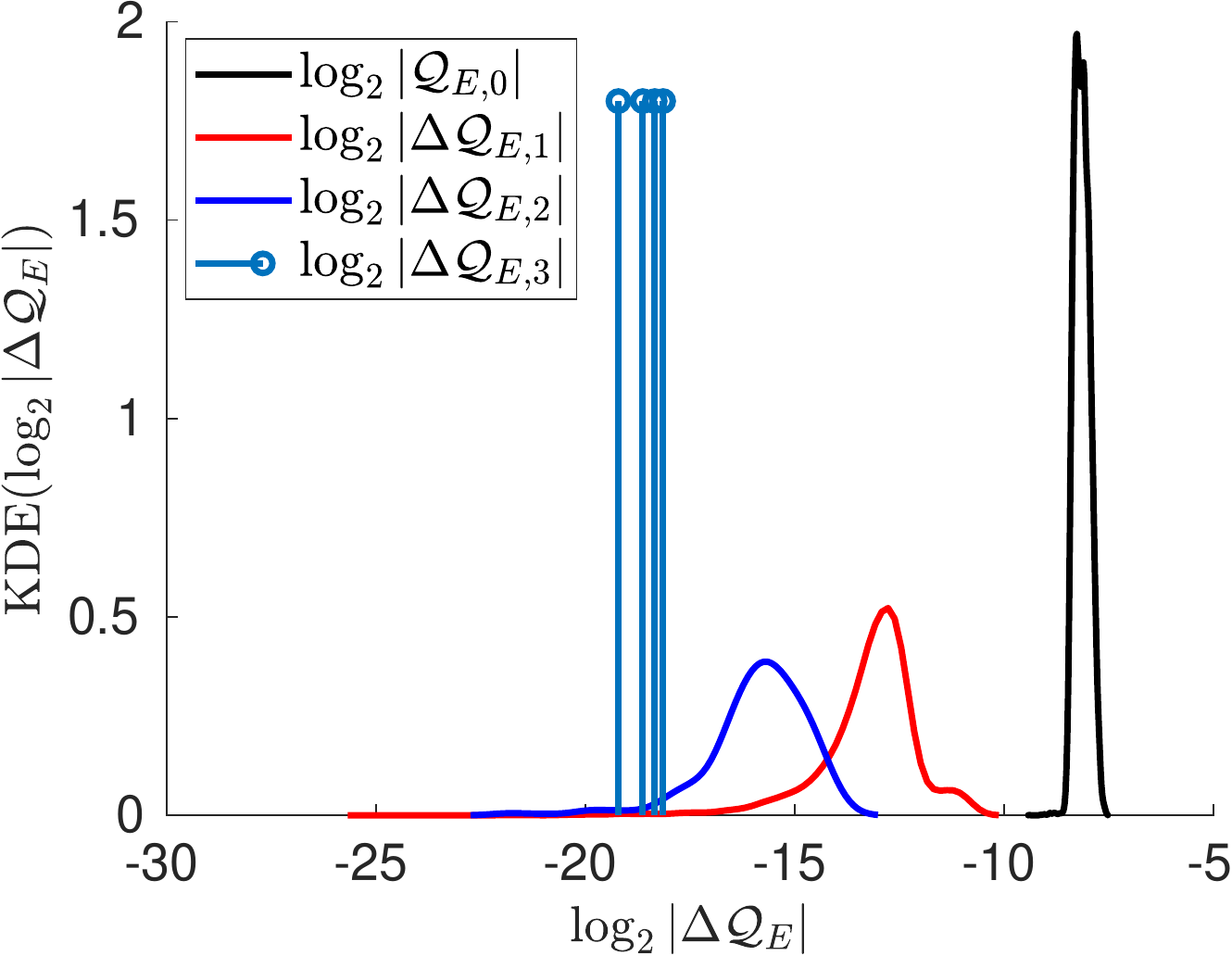}
    \includegraphics[width=0.49\linewidth]{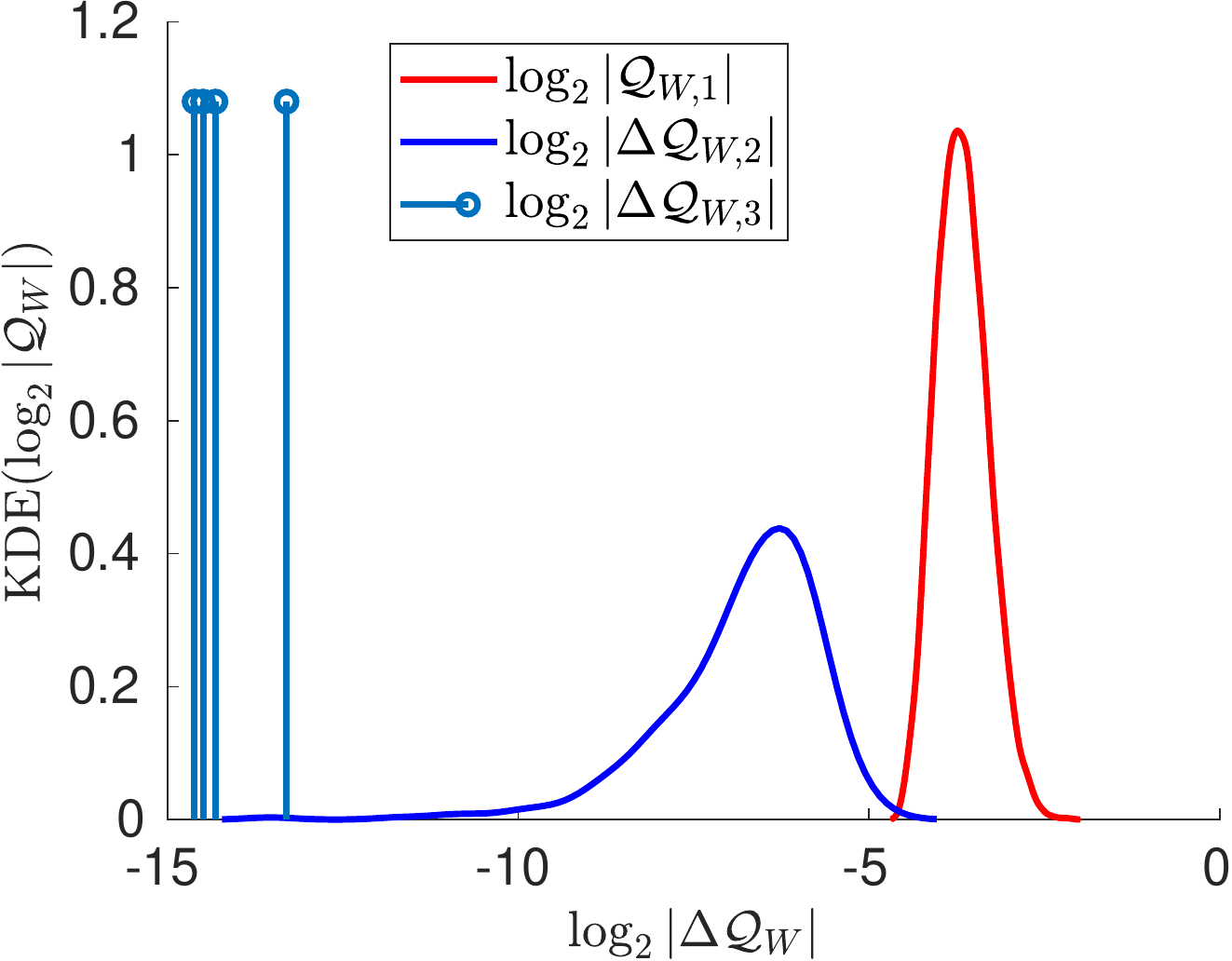}
    \captionof{figure}{Distribution of the samples of $\Delta\QoI_E$ and $\Delta\QoI_W$
      used for the bootstrapping of MLMC estimators and for the
      reference solutions; compare Table~\ref{tab:ref_hierarchy} and
      Table~\ref{tab:pool}.}  
    \label{fig:MLMC_KDE}
  \end{minipage}
\end{table}

\section{Conclusions and Future Work}
\label{sec:conc}

We have verified experimentally that MLMC techniques can significantly
reduce the computational cost of approximating expected values of
selected quantities of interest, defined in terms of misfit functions
between simulated waveforms and synthetic data with added noise, and
where the expected values are taken with respect to random parameters
modeling the uncertainties of the Earth's material properties. 
The numerical experiments conducted in this work were performed on
two-dimensional physical domains, but the extension to
three-dimensional physical domains does not create any additional
difficulties other than a higher computational cost per
sample, due to the numerical approximation of the underlying wave 
propagation model in higher spatial dimension.
Furthermore, the asymptotic complexity of the MLMC method for these
particular underlying approximation methods remains the same in the
three-dimensional case up to logarithmic factors in the user-specified
error tolerance. 

Future work includes defining the misfit function between computed
waveforms from three-dimensional simulations and actual measurement
data obtained in field studies, instead of synthetic data with added
noise, thus addressing the associated seismic inversion problem of
inferring the source location.  
Replacing the coarse level samples in the MLMC hierarchies with
samples computed using an elastic model will likely further improve
the computational gains of MLMC compared to standard MC.
Other future work is related to considering alternative ways to define
the misfit function between computed and measured seismic signals. In
this context, the normalized integration method (NIM), proposed
in~\cite{Liu2012}, and other recently proposed optimal transport-based
approaches~\cite{Metivier2018} will be considered.  

\section*{Acknowledgments}
This work is supported by the KAUST Office of 
Sponsored Research (OSR) under Award No. URF/1/2584-01-01 in 
the KAUST Competitive Research Grants Program-Round 4 (CRG2015)
and the Alexander von Humboldt Foundation. For computer time, this research used the resources of the Supercomputing Laboratory at KAUST, under the development project k1275.
The authors are grateful to Prof. Martin Mai and Dr. Olaf Zielke,
Dr. Luis F.R. Espath and Dr. Håkon Hoel, Prof. Mohammad Motamed, Prof. Daniel Appelö, and Prof. Jesper Oppelstrup for valuable
discussions and comments. 
We are grateful for the support provided by Dr.~Samuel~Kortas,
Computational Scientist, High Performance Computing, KAUST. In
particular, we are using the job-scheduler extension \emph{decimate
  0.9.5}~\cite{decimate}, developed by Dr.~Kortas. 
We would also like to acknowledge the use of the open
source software package {\tt SPECFEM2D}~\cite{SPECFEM2D}, 
provided by Computational Infrastructure for Geodynamics
\texttt{(http://geodynamics.org)} which is funded by the National
Science Foundation under awards EAR-0949446 and EAR-1550901.

M.~Ballesio, J.~Beck, A.~Pandey, E.~von~Schwerin, and R.~Tempone are 
members of the KAUST SRI Center for Uncertainty Quantification in 
Computational Science and Engineering.

\bibliographystyle{plain}
\bibliography{./Seismology_MLMC.bib}

\end{document}